\magnification=\magstep1
\input amstex
\UseAMSsymbols
\input pictex
\NoBlackBoxes


	    \def\op{{\text{op}}}
	       \def\Rep{\operatorname{Rep}}
	       \def\m{\operatorname{m}}
	       \def\n{\operatorname{n}}
	       \def\mo{\operatorname{mod}}
	          \def\Conic{\operatorname{Conic}}
		     \def\rep{\operatorname{rep}}

			      \def\Hom{\operatorname{Hom}}
			         \def\End{\operatorname{End}}
				    \def\Ext{\operatorname{Ext}}
				       \def\Tor{\operatorname{Tor}}
				          
					     \def\add{\operatorname{add}}

							    \def\bdim{\operatorname{\bold{dim}}}
							         
								      \def\arr#1#2{\arrow <1.5mm> [0.25,0.75] from #1 to #2}

									        

\vglue2truecm

\centerline{\bf The Representation Theory of Dynkin Quivers.}
		\medskip 
\centerline{\bf Three Contributions.}
		\bigskip
\centerline{Claus Michael Ringel}
		  	    \bigskip\bigskip
The representations of the Dynkin quivers and the corresponding Euclidean quivers
are treated in many books. It is well-known that the Dynkin quivers are
representation-finite, that the dimension vectors of the indecomposable representations
are just the corresponding positive roots and that any (not necessarily finite-dimensional)
representation is the direct sum of indecomposable representations. The Euclidean quivers
are the minimal tame quivers and also here, the category of all representations is well
understood and various methods of proof are known.

These notes present three building blocks for dealing with 
representations of Dynkin (and Euclidean) quivers.
They should be helpful as part of a direct approach to study representations of quivers,
and they shed some new light on properties of Dynkin and Euclidean quivers.
The presentations is based on
lectures given at SJTU (Shanghai) in 2011 and 2015 and also at KAU (Jeddah) in 2012 and at IYET
(Izmir) in 2014. 
	\bigskip
	
Part 1 considers the type $\Bbb A_n$. As one knows, any (not necessarily finite-dimensional)
representation is the direct sum of thin representations. We provide
a straight-forward arrangement of arguments
in order to avoid indices and clumsy inductive considerations, but also avoiding somewhat
fancy tools such as the recursive use of the 
Bernstein-Gelfand-Ponomarev reflection functors or bilinear forms
and root systems. As we will see, the essential case to be considered is the case $\Bbb A_3$
as studied in any first year course in Linear Algebra.

In Part 2 we deal with the case $\Bbb D_n$. If $Q$ is a quiver of type $\Bbb D_n$ with
$n\ge 4$, the category of all representations of $Q$ contains a full subcategory 
which is equivalent to the category of representations of a quiver of type $\Bbb D_{n-1}$ such that the
remaining indecomposable representations are thin. This provides an inductive procedure starting with
$\Bbb D_3 = \Bbb A_3.$ No further linear algebra knowledge is required for dealing with
the cases $\Bbb D_n$. The proof relies on some easily established equivalences of
categories. As an afterthought, we give an interpretation in terms of
Auslander-Reiten quivers.

Part 3 is devoted to the quivers $\Delta$ of type $\Bbb E_6, \Bbb E_7, \Bbb E_8$
and the categories $\rep\Delta$ of finite-dimensional representations. 
The cases $\Bbb E_6, \Bbb E_7, \Bbb E_8$ can be
reduced to $\Bbb A_5, \Bbb D_6$ and $\Bbb E_7$, respectively, and there is a further
reduction to $\Bbb A_2\sqcup\Bbb A_2, \Bbb A_5, \Bbb E_6$, respectively.   
It seems to be of interest that this reduction scheme follows the rules
of the magic Freudenthal-Tits square. 
In this way, we obtain a uniform way to construct 
the maximal indecomposable representation $M$ in 
each of the cases $\Bbb E_6, \Bbb E_7, \Bbb E_8$.
Also, this provides a unified method to deal with 
the corresponding Euclidean quivers of type 
$\widetilde{\Bbb E}_6,\widetilde{\Bbb E}_7$ and $\widetilde{\Bbb E}_8$, namely 
to construct the one-parameter families, as well as
representatives of the $\tau$-orbits of the simple regular representations.  

Actually, a similar procedure can also be used for the cases 
$\Bbb D_n$ and $\widetilde {\Bbb D}_n$. 
For all the Dynkin quivers $\Delta$ of type $\Bbb D_n$ and $\Bbb E_m$
and for the corresponding Euclidean quivers $\widetilde \Delta$, 
we obtain a reduction to $\Bbb D_4$ and $\widetilde{\Bbb D}_4$,
respectively. 

Our approach in Part 3 may be formulated also as follows. Let $\Delta$ be a
Dynkin quiver of type $\Bbb D_n$ or $\Bbb E_n$ and $y$ its exceptional vertex
(this is the uniquely determined vertex such that
the corresponding Euclidean quiver $\widetilde \Delta$
is obtained from $\Delta$ by adding 
a vertex $z$ as well as an arrow between $y$ and $z$).
We look at the quiver $\Delta'$ 
obtained from $\Delta$ by deleting $y$. 
The essential observation concerns the restriction $M|\Delta'$.
As we will see, $M|\Delta'$ is the direct sum of three
indecomposable representations; we call these representations  $A(1),A(2),A(3)$ 
the special antichain triple of $\Delta$. The thick subcategory generated by
$A(1),A(2),A(3)$ and the simple representation $S(y)$ will be said to be 
the core of $\rep\Delta$. In case $\Delta$ is of type $\Bbb E_m$, it is also of interest to
look at the quiver $\Delta''$ obtained from $\Delta$ by deleting $y$ and the
neighbors of $y$. We will show that the restriction $M|\Delta''$ is the direct sum
of four indecomposable representations and that its endomorphism ring 
is 6-dimensional. 

    \medskip 
The level of the presentation varies considerably and increases throughout the discussion.
Whereas Part 1 is completely elementary, just based on some results in linear algebra,
the further text is less self-contained. The inductive proof of Theorem 2.1 requires no 
prerequisites, but for a proper understanding one needs to have some knowledge about
Auslander-Reiten quivers, the use of antichains and simplification as well as 
about perpendicular categories.
In Part 3 we will use hammocks and one-point extensions. 
We have tried to incorporate most of the relevant definitions, but a neophyte may 
have to consult a standard reference such as [ARS], or also [R3]. Comments 
concerning the possible use of the material presented here can be found in 
remarks at the end of each part. The Panoramic View at the end of the paper 
shows how the three Parts can be completed for a full discussion of the borderline
between representation-finite and representation-infinite quivers.
	   \bigskip
{\bf Acknowledgment.} The author is indebted to Markus Schmidmeier for many helpful comments.
In particular, he suggested to use the word conical. In first versions of part 3, we had
restricted the attention to quivers of type $\Bbb E_6,\ \Bbb E_7,\ \Bbb E_8.$ On the basis
of discussions with him, also the cases $\Bbb D_n$ are now treated in a parallel way.  
     \bigskip\medskip 
\vfill
{\bf Preliminaries.}
     \medskip
We fix some field $k$. All vector spaces to be considered will be $k$-spaces, 
the algebras will be $k$-algebras, the categories will be $k$-categories. Note that in Part 1 and at the
beginning of Part 2, the vector spaces to be considered may be infinite-dimensional. Later, we will
restrict the discussion to finite-dimensional representations. 
	 \medskip 
{\bf 0.1. Quivers.} We denote by $Q = (Q_0,Q_1,s,t)$ a finite quiver, thus $Q_0,Q_1$ are finite sets and
$s,t\:Q_1\to Q_0$ are (set) maps. The elements of $Q_0$ are called vertices, the 
elements in $Q_1$ arrows. Usually, an arrow $\alpha\in Q_1$ will be drawn in the following way
$s(\alpha) \overset{\alpha}\to{\rightarrow} t(\alpha)$ or $\alpha\:s(\alpha)\to t(\alpha)$. 

If $Q$ is a quiver, its underlying graph $\overline Q$ is a graph (possibly with loops and multiple edges),
it is obtained from $Q$ by replacing the maps $s,t$ by the map $\{s,t\}$
which sends $\alpha\in Q_0$ to the subset $\{s(\alpha),t(\alpha)\}$ of $Q_1$. The main concern of these
notes are the {\it Dynkin quivers:} these are the quivers $Q$ such that the underlying graph $\overline Q$
is one of the following Dynkin graphs $\Bbb A_n, \Bbb D_n, \Bbb E_6, \Bbb E_7, \Bbb E_8$:
$$
{\beginpicture
\setcoordinatesystem units <0.6cm,0.6cm>

\put{\beginpicture
\multiput{$\circ$} at  1 0  2 0  4 0 /
\multiput{$\star$} at 0 0  5 0 /
\plot 0.2 0  0.8 0 /
\plot 1.2 0  1.8 0 /
\plot 2.2 0  2.5 0 /
\plot 3.5 0  3.8 0 /
\plot 4.2 0  4.8 0 /
\put{$\cdots$} at 3 0 
\put{$\Bbb A_n$} at -2 0
\put{} at 6 0
\endpicture} at 0 0

\put{\beginpicture
\setcoordinatesystem units <0.6cm,0.4cm>
\multiput{$\circ$} at 0 1  0 -1  1 0  3 0  5 0 /
\put{$\star$} at   4 0 
\plot 0.2 0.8  0.8 0.2 /
\plot 0.2 -.8  0.8 -.2 /
\plot 1.2 0  1.5 0 /
\plot 2.5 0  2.8 0 /
\plot 3.2 0  3.8 0 /
\plot 4.2 0  4.8 0 /
\put{} at 6 0
\put{$\cdots$} at 2 0 
\put{$\Bbb D_n$} at -2 0

\endpicture} at 0 -2
\put{\beginpicture
\setcoordinatesystem units <0.6cm,0.5cm>
\multiput{$\circ$} at 0 0  1 0  2 0  3 0  4 0 /
\put{$\star$} at 2 1 
\plot 0.2 0  0.8 0 /
\plot 1.2 0  1.8 0 /
\plot 2.2 0  2.8 0 /
\plot 3.2 0  3.8 0 /
\plot 2 0.2  2 0.8 /
\put{$\Bbb E_6$} at -2 0
\put{} at 6 0
\endpicture} at 0 -4
\put{\beginpicture
\setcoordinatesystem units <0.6cm,0.5cm>
\multiput{$\circ$} at  1 0  2 0  3 0  4 0  5 0  2 1 /
\put{$\star$} at   0 0 
\plot 0.2 0  0.8 0 /
\plot 1.2 0  1.8 0 /
\plot 2.2 0  2.8 0 /
\plot 3.2 0  3.8 0 /
\plot 4.2 0  4.8 0 /
\plot 2 0.2  2 0.8 /
\put{$\Bbb E_7$} at -2 0
\put{} at 6 0
\endpicture} at 0 -6
\put{\beginpicture
\setcoordinatesystem units <0.6cm,0.5cm>
\multiput{$\circ$} at 0 0  1 0  2 0  3 0  4 0  5 0   2 1 /
\put{$\star$} at   6 0 
\plot 0.2 0  0.8 0 /
\plot 1.2 0  1.8 0 /
\plot 2.2 0  2.8 0 /
\plot 3.2 0  3.8 0 /
\plot 4.2 0  4.8 0 /
\plot 5.2 0  5.8 0 /
\plot 2 0.2  2 0.8 /
\put{$\Bbb E_8$} at -2 0
\endpicture} at 0 -8
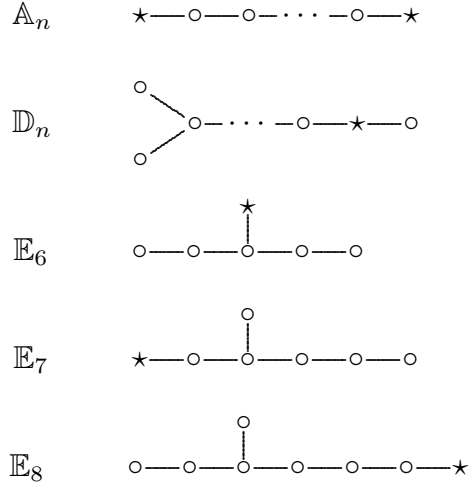
\endpicture}
$$
with vertices drawn as $\circ$ or $\star$. 
The index $n$ in $\Bbb A_n, \Bbb D_n, \Bbb E_n$ is just the number of vertices.
We call the vertices $\star$ 
{\it exceptional}. The exceptional vertices of a Dynkin graph 
(see [R2], p.6) can be
characterized in many different ways, using for example the root system or the 
quadratic form attached to the graph. 
Note that the graphs of type $\Bbb D_n$ and
$\Bbb E_m$ have a unique exceptional vertex. 

If $\Delta$ is a Dynkin graph different from $\Bbb A_1$, the
corresponding Euclidean graph $\widetilde\Delta$ is  obtained from $\Delta$ by adding a vertex $z$
and one edge between $z$ and any exceptional vertex of $\Delta$. If $\Delta = \Bbb A_1$, say with vertex $y$,
we obtain  the corresponding Euclidean graph $\widetilde{\Bbb A}_1$ by adding a vertex $z$
and two edges between $y$ and $z$. In addition, the graph with one vertex and one loop may also be
considered as a Euclidean diagram, it will be labeled $\widetilde{\Bbb A}_0$. The 
Euclidean quivers of type $\widetilde {\Bbb A}_n$ will be exhibited  
at the end of this preliminary section, the remaining ones 
(of type $\widetilde{\Bbb D}_n, \widetilde{\Bbb E}_6,\widetilde{\Bbb E}_7,
\widetilde{\Bbb E}_8$) in Part 3.

An {\it arm} $(Q',x)$ of length $t$ is a pair consisting of a quiver $Q'$ 
of type $\Bbb A_t$ and
a vertex $x$ of $Q'$ which has at most one neighbor in $Q'$ 
(such an arm is said to be {\it proper} provided its length is at least 2).
We say that a quiver $Q$ {\it has an arm} $(Q',x)$ provided $Q'$ is a full subquiver of $Q$
and there are no arrows between the vertices of $Q$ outside of $Q'$ and the vertices in $Q'$
different from $x$. Thus $Q$ has an arm of length $t$ provided $\overline Q$ has the form
$$
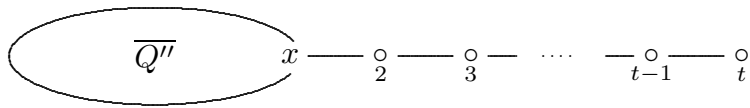
{\beginpicture
\setcoordinatesystem units <1cm,1cm>
\put{$x$} at 0 0
\put{$\ssize 2$} at 1 -.2
\put{$\ssize 3$} at 2 -.2
\put{$\ssize t-1$} at 4  -.2
\put{$\ssize t$} at 5 -.2
\multiput{$\circ$} at 1 0  2 0  4 0  5 0 /
\put{$\overline{Q''}$} at -1.5 0
\plot 0.2 0  0.8 0 /
\plot 1.2 0  1.8 0 /
\plot 2.2 0  2.5 0 /
\plot 3.5 0  3.8 0 /
\plot 4.2 0  4.8 0 /
\setdots <1mm>
\plot 2.8 0  3.2 0 /
\setsolid
\ellipticalarc axes ratio 3:1 -320 degrees from 0 -0.2 center at -1.5 0
\multiput{} at 0 -.7  0 .7 /
\endpicture}
$$
for some quiver $Q''$ (we say that we have attached 
at $x$ in $Q''$ an arm of length $t$).
    \medskip 
A quiver $Q$ is a {\it tree quiver} provided the underlying graph $\overline Q$ is a tree.
Tree quivers can be constructed inductively: The quiver $\Bbb A_1$ (one vertex, no arrow) is
a tree quiver. A quiver with at least two vertices is a tree quiver, 
provided it is obtained from a tree quiver by attaching an arm of length 2. 
Note that all the Dynkin quivers are tree quivers.

If $X$ is a set of vertices of the quiver $Q$, we write $\langle X\rangle$ for the full
subquiver of $Q$ with vertices in $X$ (and $\{X\}$ for the subquiver with 
vertix set $X$ and no arrows).
If the vertices of $Q$ are labeled by the natural
numbers $1,2,\dots,n$, and $1\le i \le j \le n$, we write $[i,j]$ for $\langle i,i\!+\!1,\dots,j\rangle.$
	\medskip

{\bf 0.2. Representations.}
If $Q$ is a quiver, a representation $M = (M_x,M_\alpha)$ is given by vector spaces $M_x$ for
$x\in Q_0$ and linear maps $M_\alpha\:M_{s(\alpha)}\to M_{t(\alpha)}$ for $\alpha\in Q_1$;
instead of $M_\alpha$ we often write just $\alpha$. Note that the vector spaces considered here may be infinite-dimensional. We denote by $\bdim M$ the dimension vector of $M$, with
$(\bdim M)_x$ being the $k$-dimension of the vector space $M_x,$ for $x\in Q_0$.
The {\it support} of a representation $M$ is the set of vertices $x$ with $M_x\neq 0.$ 

For any subquiver $Q'$ of the quiver $Q$, we define the representation $M(Q')$ of $Q$ by 
$M(Q')_x = k$ for all vertices $x\in Q'$ and $M(Q')_\alpha = 1_k$ for all arrows $\alpha$ in $Q'$,
whereas $M(Q')_x = 0$, if $x$ is a vertex in $Q$, but not in $Q'$, and $M(Q')_\alpha = 0$ 
if $\alpha$ is an arrow in $Q$, but does not belong to $Q'$. If $x$ is a vertex of $Q$, 
let $S(x) = M(\{x\})$, this is called the simple representation corresponding to $x$.
In case the quiver $Q$ is finite and acyclic (that is, $Q$ has no proper cyclic path),
we denote by $P(x)$ the projective cover, by
$I(x)$ the injective envelope of $S(x)$.
       \medskip 
If $M, M'$ are representations of a quiver $Q$, a homomorphism $f = (f_x)\:M \to M'$ is given by
linear maps $f_x\:M_x \to M'_x$ for $x\in Q_0$
such that $M'_\alpha f_{s(\alpha)} = f_{t(\alpha)}M_\alpha$ 
holds for all $\alpha\in Q_1.$ The composition $gf$
of $f\:M \to M'$ and $g\:M'\to M''$ is defined by $(gf)_x = g_xf_x$ for all $x\in Q_0$. The category
of representations of $Q$ will be denoted by $\Rep Q$, the full subcategory of all finite-dimensional
representations by $\rep Q$; these are abelian categories. 
If $M(i),\ i\in I,$ is a family of representations of a quiver $Q$, the
direct sum $M = \bigoplus_{i\in I} M(i)$ is defined by $M_x = \bigoplus M(i)_x$ and $M_\alpha =
\bigoplus M(i)_\alpha.$ A representation $M$ is said to be indecomposable provided it is not zero and  
for any direct sum decomposition $M = M'\oplus M''$, either $M' = 0$ or $M'' = 0.$

    \medskip 
{\bf 0.3. Thin representations of a quiver.}
A representation $M$ is said to be {\it thin,} provided $\dim M_x \le 1$ for all vertices $x$.
Of course, for any subquiver $Q'$ of $Q$, 
the representation $M(Q')$ is thin.
   \medskip 
The following properties of the thin representations will be relevant.
    \medskip 
{\bf (T1)} {\it Let $Q$ be a quiver. Then $Q'  \mapsto M(Q')$ provides an
injective map from the set of connected subquivers of $Q$ into the set of isomorphism
classes of indecomposable thin representations of $Q$.}
	\medskip
{\bf (T2)} {\it If $Q$ is a tree quiver, then any thin indecomposable 
representation with support $Q$ is isomorphic to $M(Q)$.}
	       \medskip
As an immediate consequence of (T1) and (T2) we obtain: 
   \medskip 
{\bf (T3)} {\it If $Q$ is a tree quiver, then $Q' \mapsto M(Q')$ provides a
bijection between the set of connected subquivers of $Q$ and the set of isomorphism
classes of indecomposable thin representations of $Q$.}
	\medskip 
If $Q$ is not 
the disjoint union of tree quivers, then there are indecomposable thin representations
which are not isomorphic to representations of the form $M(Q')$.
Here is the essential assertion:
{\it Assume that the underlying graph of $Q$ is a quiver of type $\widetilde{\Bbb A}_n$, thus its underlying
graph is the cycle with $n+1$ vertices:
$$  
 
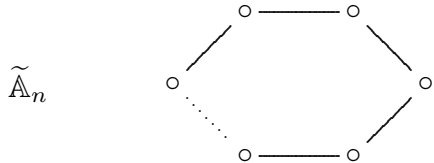
{\beginpicture
 \setcoordinatesystem units <.8cm,.8cm> 
\put{$\widetilde{\Bbb A}_n$} at -2 0 
\multiput{$\circ$} at 0 0  1 1  1 -1  2.5 1  2.5 -1 3.5 0 /
\plot 0.2 .2   0.8 0.8  /
\plot 1.2 -1  2.3 -1 /
\plot 1.2 1  2.3 1 /
\plot 2.7 0.8  3.3 0.2 /  
\plot 2.7 -.8  3.3 -.2 / 
\setdots <1mm> 
\plot 0.2 -.2   0.8 -.8  /
\endpicture}
$$
Let $\alpha$ be one of the arrows. Let $\lambda\in k$ be different from $0$ and $1$.
Define $M$ as follows: $M_x = k$ for all vertices $x$, let
$M_\alpha = \lambda$ and $M_\beta = 1$ for the remaining arrows $\beta$. Then
$M$ is thin, indecomposable, with support $Q$, but $M$ is not isomorphic to $M(Q)$.}
    \bigskip\bigskip 
\centerline{\bf Part 1. The quivers of type $\Bbb A$.}
		\bigskip 
Our aim is to provide a straightforward proof of the following well-known result:
    \medskip
{\bf 1.1. Theorem.} {\it Any representation of a quiver of type $\Bbb A$ is the direct
sum of thin indecomposable representations.}
    \medskip
We should recall that according to (T2), any thin indecomposable representation
is of the form $M(Q')$, where $Q'$ is a connected subquiver of $Q$.
   \medskip
Our proof will rely on a detailed study of the quivers of type $\Bbb A_3$.
For $n\ge 4$, we will first use induction and then we will invoke the knowledge about the
case $\Bbb A_3.$ Thus, let us look at the quivers of type $\Bbb A_3$. 
     \medskip
{\bf The case $\Bbb A_3$.}
There are three different orientations a quiver of type $\Bbb A_3$ may have. We discuss
these orientations one after the other. But the proof in all three cases is based on
the same result: given two subspaces $U, U'$ of a vector space $V$, there
is a basis $\bold B$ of $V$ which is compatible with 
each of the two subspaces (a basis $\bold B$ of
a vector space $V$ is said to be {\it compatible} with the subspace $U$ of $V$ provided 
$\bold B\cap U$ is a basis of $U$).
       \medskip 
{\bf (a) The $2$-subspace quiver.} By definition, this is the following quiver:
$$
{\beginpicture
\setcoordinatesystem units <1cm,1cm>
\multiput{$\circ$} at 0 0  1 0  2 0 /
\arr{0.2 0}{0.8 0}
\arr{1.8 0}{1.2 0}
\put{$1$} at 0 0.3
\put{$2$} at 1 0.3
\put{$3$} at 2 0.3
\put{$\alpha$} at 0.5 0.3
\put{$\beta$} at 1.5 0.3
\endpicture}
$$
Starting with a representation $M$, let $U$ be the image of $M_\alpha$ in $V = M_2$,
and $U'$ the image of $M_\beta$ in $V$.
Splitting off copies of $S(1)$, we can assume that $M_\alpha$ is the inclusion map of 
$U$ in $V$;
splitting off copies of $S(3)$, we can assume that $M_\beta$ is the inclusion map of 
$U'$ in $V$.
Thus we deal with a vector space $V$ with two subspaces $U,U'$:
$$
{\beginpicture
\setcoordinatesystem units <1cm,1cm>
\put{$V$} at 1 0
\put{$U$} at 0 0
\put{$U'$} at 2 0
\arr{0.2 0}{0.8 0}
\arr{1.8 0}{1.2 0}
\endpicture}
$$
Using a basis $\bold B$ of $V$ which is compatible with both $U$ and $U'$, we can write $M$ 
as follows:
$$
{\beginpicture
\setcoordinatesystem units <2.6cm,1cm>
\put{$\bigoplus_{b\in \bold B} kb$} at 1 -.05
\put{$\bigoplus_{b\in \bold B\cap U} kb$} at 0 -.05
\put{$\bigoplus_{b\in \bold B\cap U'} kb$} at 2 -.05
\arr{0.4 0}{0.6 0}
\arr{1.6 0}{1.4 0}
\endpicture}\ \ ,
$$
thus we obtain a decomposition of $M$
as a direct sum of copies of the following thin indecomposable representations:
$$
{\beginpicture
\setcoordinatesystem units <1cm,1cm>
\put{\beginpicture
\put{for any\strut} at -1.2 -.7
\put{$b\in U\cap U'$\strut} at 1 -.7
\multiput{$k$} at 0 0  1 0  2 0 /
\arr{0.2 0}{0.8 0}
\arr{1.8 0}{1.2 0}
\multiput{$1$} at 0.5 0.3  1.5 0.3 /
\endpicture} at -.9 0
\put{\beginpicture
\put{$b\in U\setminus U'$\strut} at 1 -.7
\multiput{$k$} at 0 0  1 0  /
\multiput{$0$} at 2 0 /
\arr{0.2 0}{0.8 0}
\arr{1.8 0}{1.2 0}
\multiput{$1$} at 0.5 0.3 /
\endpicture} at 3 0
\put{\beginpicture
\put{$b\in U'\setminus U$\strut} at 1 -.7
\multiput{$k$} at 2 0  1 0  /
\multiput{$0$} at 0 0 /
\arr{0.2 0}{0.8 0}
\arr{1.8 0}{1.2 0}
\multiput{$1$} at  1.5 0.3 /
\endpicture} at 6 0
\put{\beginpicture
\put{$b\notin U\cup U'$\strut} at 1 -.7
\multiput{$k$} at  1 0  /
\multiput{$0$} at 2 0 0 0 /
\arr{0.2 0}{0.8 0}
\arr{1.8 0}{1.2 0}
\multiput{$\phantom 1$} at  1.5 0.3 /
\endpicture} at 9 0
\endpicture}
$$
	\medskip
	{\bf (b) The quiver $Q$ of type $\Bbb A_3$ with linear orientation:}
	$$
	{\beginpicture
	\setcoordinatesystem units <1cm,1cm>
	\multiput{$\circ$} at 0 0  1 0  2 0 /
	\arr{0.8 0}{0.2 0}
	\arr{1.8 0}{1.2 0}
	\put{$1$} at 0 0.3
	\put{$2$} at 1 0.3
	\put{$3$} at 2 0.3
	\put{$\alpha$} at 0.5 0.3
	\put{$\beta$} at 1.5 0.3
	\endpicture}
	$$
Let $V = M_2.$
Splitting off copies of $S(1)$, we can assume 
that $M_\alpha$ is the canonical projection $V \to V/U$, where $U$ is the kernel of $M_\alpha$. 
Splitting off copies of $S(3)$, we can assume that $M_\beta$ is the 
inclusion of a subspace $U'$ of $V$.
Thus we deal with a vector space $V$ with two
subspaces $U,U'$ and consider the corresponding representation of $Q$:
$$
{\beginpicture
\setcoordinatesystem units <1cm,1cm>
\put{$V$} at 1 0
\put{$V/U$} at -0.3 0
\put{$U'$} at 2.1 0
\arr{0.8 0}{0.2 0}
\arr{1.8 0}{1.2 0}
\endpicture}
$$
As in case (1), we take a basis $\bold B$ of $V$ which is compatible with $U, U'$, in order to write $M$ 
as the direct sum of copies of the following (thin indecomposable) representations
$$
{\beginpicture
\setcoordinatesystem units <1cm,1cm>
\put{\beginpicture
\put{for any\strut} at -1.2 -.7
\put{$b\in U\cap U'$\strut} at 1 -.7
\multiput{$k$} at 2 0  1 0  /
\multiput{$0$} at 0 0 /
\arr{0.8 0}{0.2 0}
\arr{1.8 0}{1.2 0}
\multiput{$1$} at  1.5 0.3 /
\endpicture} at -.9 0
\put{\beginpicture
\put{$b\in U\setminus U'$\strut} at 1 -.7
\multiput{$k$} at  1 0  /
\multiput{$0$} at 2 0 0 0 /
\arr{0.8 0}{0.2 0}
\arr{1.8 0}{1.2 0}
\multiput{$\phantom 1$} at  1.5 0.3 /
\endpicture} at 3 0
\put{\beginpicture
\put{$b\in U'\setminus U$\strut} at 1 -.7
\multiput{$k$} at 0 0  1 0  2 0 /
\arr{0.8 0}{0.2 0}
\arr{1.8 0}{1.2 0}
\multiput{$1$} at 0.5 0.3  1.5 0.3 /
\endpicture} at 6 0
\put{\beginpicture
\put{$b\notin U\cup U'$\strut} at 1 -.7
\multiput{$k$} at 0 0  1 0  /
\multiput{$0$} at 2 0 /
\arr{0.8 0}{0.2 0}
\arr{1.8 0}{1.2 0}
\multiput{$1$} at 0.5 0.3 /
\endpicture} at 9 0
\endpicture}
$$
   \medskip
   {\bf (c) The $2$-factor-space quiver.} This is the quiver
   $$
   {\beginpicture
   \setcoordinatesystem units <1cm,1cm>
   \multiput{$\circ$} at 0 0  1 0  2 0 /
   \arr{0.8 0}{0.2 0}
   \arr{1.2 0}{1.8 0}
   \put{$1$} at 0 0.3
   \put{$2$} at 1 0.3
   \put{$3$} at 2 0.3
   \put{$\alpha$} at 0.5 0.3
   \put{$\beta$} at 1.5 0.3
   \endpicture}
   $$
We may proceed as in the previous cases, by first splitting off copies of $S(1),S(3)$
so that we are left with a vector space $V$ and two factor spaces of $V$, say 
with $V/U$ and $V/U'$ where $U,U'$ are subspaces of $V$. A basis $\bold B$ of $V$
which is compatible both with $U$ and $U'$ provides the desired direct sum decomposition
with thin direct indecomposable summands. 

But there is a second possibility for dealing with the case (c): 
We may use $k$-duality in order to reduce this case to the case (a). 
      \bigskip
{\bf Proof of Theorem 1.1.} We deal with a quiver $Q$ with the following underlying graph
$$
{\beginpicture
\setcoordinatesystem units <1cm,1cm>
\multiput{$\circ$} at 0 0  1 0  4 0  5 0 /
\plot 0.2 0  0.8 0 /
\plot 1.2 0  1.8 0 /
\plot 3.2 0  3.8 0 /
\plot 4.2 0  4.8 0 /
\put{$1$} at 0 0.3
\put{$2$} at 1 0.3
\put{$n\!-\!1$} at 4 0.3
\put{$n$} at 5 0.3
\put{$\dots$} at 2 0
\endpicture}
$$
For $1\le s\le t \le n$, we have denoted by $[s,t]$ the full subquiver with vertices
$s,s+1,\dots,t.$ According to (T3), the representations
$M([s,t])$ (with $1\le s \le t \le n)$ furnish a complete list of the indecomposable thin
representations of $Q$, up to isomorphism. We have to show that 
any representation of $Q$ is isomorphic to the direct sum of representations of the form
$M([s,t])$.
	
As we know already, the assertion is true for $n\le 3.$ Thus, consider now some
$n\ge 4$. By induction, we may assume that any indecomposable representation
of a quiver of type $\Bbb A_{n-1}$ is thin. We use this assertion first for $[1,n-1]$
and then for $[2,n]$, these are the first two steps of the proof.
    
We need another definition. 
Let $c$ be a vertex of any acyclic $\Delta$. 
We say that a representation $M$ of $\Delta$ 
is {\it $c$-conical} provided $M_\alpha$ is injective for 
any arrow $\alpha$ of $\Delta$ which points to $c$, 
and $M_\beta$ is surjective, for the remaining arrows $\beta$ of $\Delta$
(an alternative terminology calls $c$ a {\it peak} for $M$ in case $M$ is
$c$-conical). In our quiver $Q$, we use the ascending labels $1,2,\dots,n$ 
for the vertices; the arrows pointing to $c$ are the
arrows of the form $\alpha\:i \to i+1$ with $i<c$ and the 
arrows $\alpha\:j\!+\!1 \to  j$ with $c < j$).
For example, the representation $M([s,t])$ of $Q$ is $c$-conical provided $s\le c \le t.$
If $1\le i \le c \le j \le n$, and $M$ is a representation of $Q$, we say that $M$ is {\it conical on
$([i,j],c)$} provided the restriction $M|[i,j]$ is $c$-conical.
	     \medskip 
{\bf Step 1.} {\it Any representation $M$ of $Q$ is a direct sum $M = M' \oplus M''$
such that $M'$ is conical on $([1,n-1],n-1)$, whereas  the support of $M''$ is contained in 
$[1,n-2].$}

Proof. Consider the restriction $N = M|[1,n-1]$. By induction, we know that $N$ can be written
as a direct sum of representations of the form $M([s,t])$, with $1\le s \le t \le n-1$.
We decompose $N = N'\oplus N''$, where $N'$ is a direct sum of representations of the form
$M([s,n-1])$, and $N''$ a direct sum of representations of the form
$M([s,t])$ with $t \le n-2$. Note that the representation $N'$ of $[1,n\!-\!1]$ is $(n\!-\!1)$-conical.

Let us stress that $N'_{n-1} = M_{n-1}$ and $N''_{n-1} = 0.$ We define a
subrepresentation $M'$ as follows: 
$$
 M'|[1,n\!-\!1] = N' \quad \text{and} \quad M'|[n\!-\!1,n] = M|[n\!-\!1,n]
$$
(since $N'_{n-1} = M_{n-1}$, we see that $M'$ is well-defined).
And we define a
subrepresentation $M''$ by  
$$
 M''|[1,n\!-\!1] = N'' \quad \text{and}\quad M''|[n\!-\!1,n] = 0
$$ 
(since $N''_{n-1} = 0$, we see that also $M''$ is well-defined).
Of course, we obtain in this way a direct decomposition $M = M'\oplus M''$.
We have: $M'|[1,n\!-\!1] = N'$ is $(n\!-\!1)$-conical, 
and the support of $M''$ is equal to the support
of $N''$, thus contained in $[1,n-2].$
This completes the proof.
     
{\bf Step 2.} {\it Any representation $M$ of $Q$ is a direct sum $M = M' \oplus M''$
such that $M'$ is conical on $([2,n],2)$, whereas  the support of $M''$ is contained in $[3,n].$}
     
The proof is similar to the proof of step 1, but this time, we consider the subquiver
$[2,n]$. 
Consider the restriction $N = M|[2,n]$. By induction, we know that $N$ can be written
as a direct sum of representation of the form $M([s,t])$, with $2\le s \le t \le n$.
We decompose $N = N'\oplus N''$, where $N'$ is a direct sum of representations of the form
$M([2,t])$, and $N''$ a direct sum of representations of the form
$M([s,t])$ with $3\le s$. Note that the representation $N'$ of $[2,n]$ is $2$-conical.
And we stress that $N'_2 = M_2$ and $N''_2 = 0.$
We define a
subrepresentation $M'$ as follows: 
$M'|[2,n] = N'$ and $M'|[1,2] = M|[1,2]$
(since $N'_{2} = M_{2}$, we see that $M'$ is well-defined).
And we define a
subrepresentation $M''$ by  $M''|[2,n] = N''$ and $M'|[1,2]$ the zero representation 
(since $N''_{2} = 0$, we see that also $M''$ is well-defined).
We obtain in this way a direct decomposition $M = M'\oplus M''$.
We have: $M'|[2,n] = N'$ is $2$-conical, and the support of $M''$ is equal to the support
of $N''$, thus contained in $[3,n].$
This completes the proof.
     
{\bf Combining step 1 and step 2:} {\it Any representation $M$ of $Q$ can be written as a direct sum 
$M = M'\oplus M''$, such that for any arrow $\gamma$ in $[2,n-1]$, the map $M'_\gamma$
is a vector space isomorphism, whereas $M''$ is a direct sum of representations 
of the form $[s,t]$ with $3 \le s$ or with $t \le n-2.$}
   
Proof. According to step 1, we can assume that $M$ is conical on $([1,n-1],n\!-\!1)$.
We apply step 2 and obtain a direct decomposition $M = M'\oplus M''$, where
$M'$ is conical on $([2,n],2)$ and $M''$ has support in $[3,n]$. Since $M$
is conical on $([1,n-1],n\!-\!1)$, the same is true for any direct summand of $M$, thus $M'$ is
conical on $([1,n-1],n-1)$, as well as on $([2,n],2)$. 
This means that for any arrow $\gamma$ in $[2,n-1]$, the map $M'_\gamma$
is both injective and surjective, thus a vector space isomorphism.
   
{\bf Step 3.} {\it Let $M$ be a representation of $Q$ such that 
for any arrow $\gamma$ in $[2,n-1]$, the map $M'_\gamma$
is a vector space isomorphism. Then $M$ is isomorphic to a representation $M'$
with $M'_x = M_2$ for all $2\le x \le n-1$ and such that 
for any arrow $\gamma$ in $[2,n-1]$, the map $M'_\gamma$ is the identity map.}
    
Proof. Let $M'|[1,2] = M|[1,2]$, let $M'_x = M_2$ for $2\le x \le n-1$ and $M'_n = M_n$.
for any arrow $\gamma$ in $[2,n-1]$, let $M'_\gamma$ be the identity map. We do not
yet define the map $M'_\beta$, where $\beta$ is the arrow between $n-1$ and $n$.
But we start already to define an isomorphism $f\:M \to M'$. We de this inductively,
that means: we define $f|[1,x]$ starting with $x = 2.$ For $x = 2,$ we take as
$f|[1,2]$ the identity isomorphism. Now assume, we have defined $f|[1,x]$ for some
$2\le x \le n-2$ and we want to define $f_{x+1}$. Let $\gamma$ be the arrow in-between
$x$ and $y = x+1$. There are two possible orientations of $\gamma$:
In case $\gamma\:x \to y$, we define $f_y = f_x\cdot (M_\gamma)^{-1}$.
In case $\gamma\:y \to x$, we define $f_y = f_x\cdot M_\gamma$
In this way, we obtain an isomorphism $M|[1,y] \to M'|[1,y].$

Now assume we have constructed $f_x$ for $1\le x \le n-1$ such these maps $f_x$
combine to an isomorphism $M|[1,n-1] \to M'|[1,n-1].$
We now have to consider the arrow $\beta$ in-between $n-1$ and $n$.
Our aim is to define $M'_\beta$ as a map between $M_2$ and $M_n$ so that $f|[1,n-1]$
can be extended by $f_n = 1$ to an isomorphism $f\:M \to M'.$ 
Again, we have to distinguish the two possible orientations:
In case $\beta\:n\!-\!1 \to n$, we define $M'_\beta =  M_\beta\cdot (f_{n-1})^{-1}$.
In case $\beta\:n \to n\!-\!1$, we define $M'_\beta = f_{n-1}\cdot M_\beta $.
In this way, we obtain an isomorphism $f\:M \to M'.$
   
{\bf Final step.} 
Denote the arrow in $[1,2]$ by $\alpha$, the arrow in $[n-1,n]$ by $\beta$.
We can assume that we deal with a representation $M$ 
with $M_x = M_2$ for all $2\le x \le n-1$ and such that 
for any arrow $\gamma$ in $[2,n-1]$, the map $M_\gamma$ is the identity map.

We attach to $M$ 
$$
{\beginpicture
\setcoordinatesystem units <1.5cm,1cm> 
\put{$M =$} at -.7 0
\put{$\biggl($} at -.3 0
\put{$M_1$} at 0 0
\put{$M_\alpha$} at 0.5 0.3
\multiput{$M_2$} at 1 0  2 0  4 0 /
\put{$M_\beta$} at 4.5 0.3
\put{$M_n$} at 5 0 
\put{$\biggl)$} at 5.3 0
\plot 0.3 0  0.7 0 /
\plot 1.3 0.04  1.7 0.04 /
\plot 1.3 -.04  1.7 -.04 /

\plot 2.3 0.04  2.6 0.04 /
\plot 2.3 -.04  2.6 -.04 /

\plot 3.4 0.04  3.7 0.04 /
\plot 3.4 -.04  3.7 -.04 /
\multiput{$1$} at 1.5 0.3  2.45 0.3  3.55 0.3 /
\put{$\cdots$} at 3 0
\plot 4.3 0  4.7 0 /
\endpicture}
$$
the following representation $N$
$$
{\beginpicture
\setcoordinatesystem units <1.5cm,1cm>
\put{$N =$} at -.7 0
\put{$\biggl($} at -.3 0
\put{$M_1$} at 0 0
\put{$M_\alpha$} at 0.5 0.3
\put{$M_2$} at 1 0
\put{$M_\beta$} at 1.5 0.3
\put{$M_n$} at 2 0 
\put{$\biggl)$} at 2.3 0
\plot 0.3 0  0.7 0 /
\plot 1.3 0  1.7 0 /
\endpicture}
$$
of a corresponding quiver of type $\Bbb A_3$, say with vertices labeled $1,2,n$,
and with arrows $\alpha,\beta$ as in the quiver $Q$.

As we know, we can write $N$ as a direct sum of thin representations. Of course,
such a direct decomposition leads to a corresponding direct decomposition of $M$.
This completes the proof of Theorem 1.1.
   \bigskip\bigskip 
{\bf 1.2. The converse of Theorem 1.1.}
{\it Let $Q$ be a finite connected quiver. If all indecomposable
representations of $Q$ are thin, then $Q$ is of type $\Bbb A$.}
		\medskip
		Proof: Let us assume that $Q$ is connected and not of type $\Bbb A$.
Then $Q$ has a subquiver which is a cycle or which is of type $\Bbb D_4.$
It is sufficient to show: {\it If $Q$ is a quiver which is a cycle or of type $\Bbb D_4$,
then there is an indecomposable representation which is not thin.}

First, consider the case where $Q$ is a cycle, say 
$$  
 
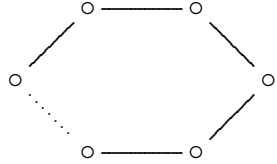
{\beginpicture
 \setcoordinatesystem units <.8cm,.8cm> 
\multiput{$\circ$} at 0 0  1 1  1 -1  2.5 1  2.5 -1 3.5 0 /
\plot 0.2 .2   0.8 0.8  /
\plot 1.2 -1  2.3 -1 /
\plot 1.2 1  2.3 1 /
\plot 2.7 0.8  3.3 0.2 /  
\plot 2.7 -.8  3.3 -.2 / 
\setdots <1mm> 
\plot 0.2 -.2   0.8 -.8  /
\endpicture}
$$
with vertices labeled $1,2,\dots, n$ and an arrow $\alpha(i)$ in-between $i$ and $i+1$, for
$1\le i < n$ and an arrow $\alpha(n)$ in-between $n$ and $1$. Without loss of generality, we can 
assume that $\alpha(n)\:n \to 1.$

Let $M_{i} = k^2,$ for all vertices $i$, take for all the arrows $\alpha(i)$
with $1\le i < n$ the identity map $k^2\to k^2$ and let $M_{\alpha(n)} =
\left[\smallmatrix 0 & 1 \cr  0 & 0 \endsmallmatrix\right]$. A straightforward calculation
shows that $M$ is indecomposable.
       
Second, we have to consider the quivers of type $\Bbb D_4$. It is easy to see (and
will be shown in Part 2) that such a quiver has a unique indecomposable representation which
is not thin. This completes the proof.
   \bigskip
Let us mention several applications of Theorem 1.1.

\bigskip
{\bf 1.3. Pairs of filtrations.}
The first application concerns vector spaces $V$ with 
two filtrations $(U_i)_{1\le i \le n}$ and $(U'_j)_{1\le j \le m}$:
$$
\gather
 0 = U_0 \subseteq  U_1 \subseteq U_2 \subseteq \cdots \subseteq U_n \subseteq  U_{n+1} = V, \cr
  0 = U'_0 \subseteq U'_1 \subseteq U'_2 \subseteq \cdots \subseteq U'_m \subseteq  U'_{m+1} = V.
  \endgather
  $$
  \medskip 
{\bf Corollary (first formulation).} {\it Let $V$ be a vector space with two filtrations
$(U_i)_{1\le i \le n}$ and $(U'_j)_{1\le j \le m}$.
Then there is a basis $\bold B$ of $V$ which is compatible with all the subspaces $U_i,U'_j.$}
   \medskip
   Proof: Consider the quiver $Q$ of type $\Bbb A_{n+m+1}$ with vertices labeled
   $1,\dots,n,\omega,  m',\dots,1'$ and the orientation
   $$
   {\beginpicture
   \setcoordinatesystem units <1.4cm,1cm>
   \multiput{$\circ$} at 0 0  1 0  3 0  4 0  5 0  7 0  8 0 /
   \put{$1$} at 0 0.3
   \put{$2$} at 1 0.3
   \put{$n$} at 3 0.3
   \put{$\omega$} at 4 0.3
   \put{$m'$} at 5 0.3
   \put{$2'$} at 7 0.3
   \put{$1'$} at 8 0.3
   \arr{0.2 0}{0.8 0}
   \arr{2.5 0}{2.8 0}
   \arr{3.2 0}{3.8 0}
   \arr{4.8 0}{4.2 0}
   \arr{5.5 0}{5.2 0}
   \arr{7.8 0}{7.2 0}
   \multiput{$\cdots$} at 2 0  6 0 /
   \plot 1.2 0  1.5 0 /
   \plot 6.8 0  6.5 0 /
   \endpicture}
   $$
   The two filtrations yield the following representation $M$ of $Q$ (all maps are
   inclusion maps):
   $$
   {\beginpicture
   \setcoordinatesystem units <1.4cm,1cm>
   \put{$U_1$} at 0 0
   \put{$U_2$} at 1 0
   \put{$U_n$} at 3 0
   \put{$V$} at 4 0
   \put{$U'_m$} at 5 0
   \put{$U'_2$} at 7 0
   \put{$U'_1$} at 8 0
   \arr{0.2 0}{0.8 0}
   \arr{2.5 0}{2.8 0}
   \arr{3.2 0}{3.8 0}
   \arr{4.8 0}{4.2 0}
   \arr{5.5 0}{5.2 0}
   \arr{7.8 0}{7.2 0}
   \multiput{$\cdots$} at 2 0  6 0 /
   \plot 1.2 0  1.5 0 /
   \plot 6.8 0  6.5 0 /
   \endpicture}
   $$
If we write this representation as a direct sum of thin representations $N$,
 and choose in any of these direct summands a non-zero element
 $b \in N_\omega$, we obtain the required basis $\bold B$. 
 \medskip
We may reformulate the Corollary as follows:
   \medskip 
{\bf Corollary (second formulation).} {\it Let $V$ be a vector space with two filtrations
$(U_i)_{1\le i \le n}$ and $(U'_j)_{1\le j \le m}$.
For  $1\le i\le n+1$ and $1\le j\le m+1$, let $C(i,j)$ be a subspace of $U_i\cap U'_j$ with
$$
  \left((U_i\cap U'_{j-1})+(U_{i-1}\cap U'_j)\right) \oplus C(i,j) = U_i\cap U'_j.
$$
Then}
$$ V = \bigoplus_{1\le i\le n+1 \atop
1\le j\le m+1} C(i,j).
$$

Here is the relationship between the assertions of the two formulations. 
     \medskip 
If $\bold B$ is a basis of $V$ compatible with all the subspaces $U_i,U'_j$, let
$\bold B(i,j)$ be the set of elements of $b\in\bold B$ with $b\in U_i\cap U'_j$ such that
$b$ does not belong to $U_{i-1}$ nor to $U'_{j-1}$. Then $\bold B$ is the disjoint union of the
subsets $\bold B(i,j)$ and $C(i,j)$ can be defined as the subspace generated by $\bold B(i,j)$. 

Conversely, assume we have written $V$ as the direct sum $\bigoplus C(i,j)$ with
subspaces $C(i,j)$ of $V$ such that 
$\left((U_i\cap U'_{j-1})+(U_{i-1}\cap U'_j)\right) \oplus C(i,j) = U_i\cap U'_j.$
Let $\bold B(i,j)$ be a basis of $C(i,j)$ and $\bold B$ 
the disjoint union of the sets $\bold B(i,j)$. 
Then $\bold B$ is a basis von $V$ which is compatible with all the subspaces $U_i,U'_j.$

The Corollaries have been obtained by looking at a quiver $Q$ of type $\Bbb A$ and writing 
a representation $M$ of $Q$ as a direct sum of thin representations, see the proof of
Corollary 1.3. Conversely, we may use
the assertions of the corollaries in order to recover the direct sum decomposition of $M$,
as follows. For any pair $i,j$, let $M(i,j)$
be the thin indecomposable representation with support the vertices $x$ between $i$ and $j'$. 
Then 
$$
 M = \bigoplus_{i,j} M(i,j)\otimes_k C(i,j). 
$$
	\medskip 
{\bf 1.4. Arms of a quiver.}
A second application concerns quivers with an arm. 
Assume that $Q$ has an arm $(Q',x)$, and that $M$ is a representation of $Q$.
We will say that $M$ is {\it conical on $(Q',x)$} provided 
$M_\alpha$ is injective for any arrow $\alpha$ of $Q'$ which points to $x$, 
and $M_\beta$ is surjective, for the remaining arrows $\beta$ of $Q'$ (thus if the restriction 
of $M$ to $Q'$ is $x$-conical, as defined in the proof of Theorem 1.1). 
   \medskip 
{\bf Corollary.} {\it Assume that $Q$ has an arm $(Q',x)$. 
Then any representation $M$ of $Q$ can be decomposed $M = M'\oplus M''$, where $M'$
is conical on $(Q',x)$ and the support of $M''$ is contained in $Q'\setminus\{x\}$.}
   \medskip 
In particular, this implies: {\it If  $Q$ has an arm $(Q',x)$ and
$M$ is an indecomposable representation of $Q$ 
such that $M_x \neq 0$, then $M$ is conical on $(Q',x)$.}
     \bigskip
Proof of Corollary 1.4.
Let $M$ be a representation of $Q$. According to Theorem 1.1, we can write the restriction 
$M|Q'$ of $M$ to $Q'$ as a direct sum of thin indecomposable representations $X(i)$ with
$i\in I$. Let $I'$ be the set of indices $i\in I$ such that $X(i)_x \neq 0$ and let 
$I''$ be the set of indices $i\in I$ such that $X(i)_x = 0$. 
For any vertex $y\in Q'$, let $M'_y = \bigoplus_{i\in I'} X(i)_y$ and 
$M''_y = \bigoplus_{i\in I''} X(i)_y$. If $y$ is a vertex of $Q\setminus Q'$, let 
$M'_y = M_y,$ and $M''_y = 0.$ It is easy to see that we obtain in this way
a direct decomposition $M = M'\oplus M''$. The representations $X(i)$ with $i\in I'$
are conical on $(Q',x)$, thus $M'$ is conical on $(Q',x)$. Since the representations
$X(i)$ with $i\in I''$ satisfy $X(i)_x = 0$, it follows that 
the support of $M''$ is contained in $Q'\setminus\{x\}$.
    \bigskip\bigskip  
Until now, we have tried to involve just knowledge from a Linear Algebra course. From now 
on, we will make references to some results in ring and module theory.
On the one hand, 
the considerations which follow may be seen 
in this way in a broader setting. On the other
hand, some of the further proofs do rely in an essential way 
on known results which have to be mentioned.
The methods to be used
will be the Auslander-Reiten theory including hammocks, one-point extensions 
as well as tilting theory and thick subcategories.
   \bigskip 
{\bf 1.5. Star quivers.}
We say that a quiver $Q$ is a {\it star quiver
with center $c$ and arms $(Q(1),c),\dots, (Q(s),c)$} provided 
$Q$ is the union of the quivers $Q(i)$, the pairs $(Q(i),c)$ with $1\le i \le s$ are arms of $Q$
and $Q(i)\cap Q(j) = \{c\}$ for all pairs $i\neq j.$ 
Thus, the underlying graph of $Q$ has the following shape:
$$
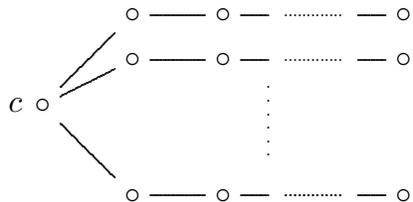
{\beginpicture
\setcoordinatesystem units <1cm,1cm>
\multiput{$\circ$} at 0 0  1 1  1 0.5  1 -1   2 1  2 .5  2 -1  4 1  4 .5  4 -1 /
\put{$c$} at -.3 0 
\plot 0.2 0.2  0.8  0.8 /
\plot 0.2 0.1  0.8  0.4 /
\plot 0.2 -.2  0.8  -.8 /
\plot 1.2 1  1.8  1 /
\plot 1.2 .5  1.8  .5 /
\plot 1.2 -1  1.8  -1 /

\plot 2.2 1  2.5  1 /
\plot 2.2 .5  2.5  .5 /
\plot 2.2 -1  2.5  -1 /

\plot 3.5 1  3.8  1 /
\plot 3.5 .5  3.8  .5 /
\plot 3.5 -1  3.8  -1 /
\setdots <.5mm>

\plot 2.7 1  3.3  1 /
\plot 2.7 .5  3.3  .5 /
\plot 2.7 -1  3.3  -1 / 
\plot 2.7 -1  3.3  -1 /

\setdots <1.5mm>
\plot 2.5 .2  2.5 -.7 /
\endpicture}
$$
If $Q$ is a star quiver with center $c$ and arms $(Q(i),c)$ of length $t(i)$, with $1\le i \le s$,
such that $t(1)\ge t(2)\ge \cdots \ge t(s)\ge 2$, then the sequence $(t(1),\cdots,t(s))$
is called the {\it type} of the star quiver $Q$ (with respect to $c$); 
the type of the star quiver $\Bbb A_1$ is 
by definition the empty sequence (or the sequence $(1)$).
A representation $M$ of a star quiver with center $c$ is said to be {\it conical} provided
it is conical on all the arms. (In case we deal with a star quiver with at least 3 proper arms,
the center of $c$ is uniquely determined; otherwise we need the reference to the center $c$.)
    \medskip 
{\bf Corollary.} {\it Assume that $Q$ is a star quiver
with center $c$.
Then any representation of $Q$ is the direct sum $M = M'\oplus M''$ such that
$M'$ is conical and $M''_c = 0.$}
     \medskip
Note that the support of an indecomposable  representation $M$ of $Q$ with $M_c = 0$
lies inside one of the arms. 
In particular, we see: {\it Assume that $Q$ is a star quiver
with center $c$. Then 
there are only finitely 
many indecomposable representations $M$ of $Q$ with $M_c = 0$ and 
any other indecomposable representation of $Q$ is conical.} 
We denote by $\Conic Q$ the full subcategory of $\Rep Q$ given by all
conical representations.
	\bigskip
If all the arrows of a star quiver $Q$ point
to the center $c$, and $M$ is a conical representation of $Q$, then
we can assume that $M$ is a subspace representation (meaning that all the maps $M_\alpha$ are
inclusion maps of subspaces).  
	  \bigskip

{\bf 1.6. Proposition.} {\it If $Q, Q'$ are star quivers of the same type, then the categories $\Conic Q$
and $\Conic Q'$ are equivalent.}
    \medskip
The easiest way to prove this assertion seems to be to use tilting theory.
Let $Q$ be a star quiver with center $c$ and $\Lambda$ its path algebra.
There is a unique multiplicity-free tilting $\Lambda$-module $T$ such that $P(c)$ is
a direct summand of $T$ and such that $\Hom(P(c),T(i)) \neq 0$ for any indecomposable
direct summand $T(i)$ of $T$ (namely the slice module for the preprojective component
of $\Rep Q$ for the slice with $P(c)$ as the only source). Let $\Lambda' = \End(T)^\op$,
this is the path algebra of the quiver $Q'$ with $\overline Q = \overline{Q'}$ and $c$
the unique sink. The tilting functor $F = \Hom(T,-)\:\Rep Q \to \Rep Q'$
provides an equivalence between $\Cal F(T) = \{M\in \Rep Q \mid \Ext^1(T,M) = 0\}$ and
$\Cal X(T) =\{X \in \Rep Q'\mid \Tor_1(T,X) = 0\}.$ It is easy to see that 
$\Conic \Lambda \subseteq \Cal F(T)$, and $\Conic \Lambda' \subseteq \Cal X(T)$, and 
that $F$ sends $\Conic \Lambda$ onto $\Conic \Lambda'.$

     \bigskip\bigskip

{\bf 1.7. Final remark for Part 1.} As we have mentioned at the beginning, the aim of Part 1 is to provide
a straightforward proof for the well-known result that all representations of a quiver of type $\Bbb A_n$
are direct sums of thin representations. Why should one care for such an approach? 
The representation theory of quivers
has to be considered as a sort of higher linear algebra, thus any introductory course devoted to
the representation theory of quivers should start with the discussion of quivers which describe
basic settings in linear algebra. Such quivers are the subspace quivers, and more generally the
star quivers, but also the quivers of type $\Bbb A$. Actually, for looking at star quivers, it is advisable
to have some knowledge concerning the behavior of representations on the arms, 
thus about quivers of type $\Bbb A$. Thus, we believe that there are good reasons to start 
with quivers of type $\Bbb A$. Some techniques to work with quivers can be trained nicely in this way,
and one has to be aware that the use of bases of vector spaces still is of interest in this setting,
in contrast to the usual procedure when dealing with representations of quivers.
Anyway, since the use of bases of vector spaces is omnipresent in any linear algebra course,
the quivers of type $\Bbb A$ may be seen in this way as an intermediate step.

The usual approach to quiver representations starts with positive roots. But the general
frame of Dynkin quivers hides the special features of the quivers of type $\Bbb A$, namely the
predominance of thin representations and the possibility to work with bases. 
It seems to be worthwhile not to neglect these features since they are helpful in the more
general setting of quivers of type $\widetilde {\Bbb A}$, or even of arbitrary
special biserial algebras with their string and band modules,
but also when looking at arms of quivers, and in particular at star quivers.
    \bigskip\bigskip 
\vfill\eject
\centerline {\bf Part 2. The quivers of type $\Bbb D$.}
	    	   \bigskip
We consider now the quivers with underlying graph of type $\Bbb D_n$. Whereas our interest 
in type $\Bbb A$ was devoted to thin indecomposable representations, here we are faced with categories
for which the thin representations play a minor role. The relevant indecomposable representations 
will be collected in a full subcategory $\check{\Cal B}$ which contains all the non-thin indecomposable
representations and only few thin representations. 

First, let us consider the following quiver $Q(n)$
$$
{\beginpicture
\setcoordinatesystem units <1cm,.5cm>
\put{$Q(n)$} at -2 0 
\put{$1$} at 0 1
\put{$2$} at 0 -1
\put{$3$} at 1 0
\put{$4$} at 2 0
\put{$\cdots$} at 3 0
\put{$n$} at 4 0
\arr{0.6 0.3}{0.3 0.7}
\arr{0.6 -.3}{0.3 -.7}
\arr{1.7 0}{1.3 0}
\arr{2.7 0}{2.3 0}
\arr{3.7 0}{3.3 0}
\put{$\pi_1$} at 0.6 0.85
\put{$\pi_2$} at 0.6 -.95
\endpicture}
$$
defined for $n\ge 3$. Of course, $Q(3)$ is a quiver of type $\Bbb A_3$ (the $2$-factor-space quiver),
thus our main interest will lie in the quivers $Q(n)$ with $n\ge 4.$
     \medskip
We denote by $\check{\Cal B}$ the full subcategory of all representations $M$ of $Q(n)$ such that the
map $\left[\smallmatrix \pi_1\cr\pi_2\endsmallmatrix\right]\:M_3 \to M_1\oplus M_2$ is bijective.
Let $\Cal B = \check{\Cal B}\cap \rep Q(n)$.
    \medskip
Again, our theme is to present a proof of a well-known result, namely the following theorem.
       \medskip 

{\bf 2.1. Theorem.} {\it Any representation of $Q(n)$
is the direct sum of a representation in $\check{\Cal B}$
and thin representations.}
    \medskip
Proof. There are several steps of splitting off direct sums of thin representations in order
to obtain a representation in $\check{\Cal B}.$ Two of the steps are based on 
Theorem 1.1.
	\medskip 
{\bf (0)} As a preliminary step, we consider the quiver $\sigma Q(n)$ 
$$
{\beginpicture
\setcoordinatesystem units <1cm,.5cm>
\put{$\sigma Q(n)$} at -2 0
\put{$1$} at 0 1
\put{$2$} at 0 -1
\put{$3$} at 1 0
\put{$4$} at 2 0
\put{$\cdots$} at 3 0
\put{$n$} at 4 0
\arr{0.3 0.7}{0.6 0.3}
\arr{0.3 -.7}{0.6 -.3}
\arr{1.7 0}{1.3 0}
\arr{2.7 0}{2.3 0}
\arr{3.7 0}{3.3 0}
\put{$\mu_1$} at 0.6 0.85
\put{$\mu_2$} at 0.6 -.95
\endpicture}
$$
and denote by $\Cal N$ the full subcategory of all representations of $\sigma Q(n)$ of the form
$$
{\beginpicture
\setcoordinatesystem units <1cm,.5cm>
\put{$N\ = $} at -1.7 0
\put{$\Biggl($} at -.6 0
\put{$\Biggr)$} at 4.5 0
\put{$U_1$} at 0 1
\put{$U_2$} at 0 -1
\put{$N_3$} at 1 0
\put{$N_4$} at 2 0
\put{$\cdots$} at 3 0
\put{$N_n$} at 4 0
\arr{0.3 0.7}{0.6 0.3}
\arr{0.3 -.7}{0.6 -.3}
\arr{1.7 0}{1.3 0}
\arr{2.7 0}{2.3 0}
\arr{3.7 0}{3.3 0}
\put{$\mu_1$} at 0.6 0.85
\put{$\mu_2$} at 0.6 -.95
\endpicture}
$$
where $U_1, U_2$ are subspaces of $N_3$ and $\mu_i\:U_i\to N_3$ the inclusion maps.

Let $N$ be a representation in $\Cal N$. We look at the representation
$$
  \overline N = \quad \left(\ U_1\oplus U_2 
   @>\left[\smallmatrix \mu_1 & \mu_2\endsmallmatrix\right]>> N_3 @<<< N_4 @<<< \cdots @<<< N_n \ \right)
$$ 
of the following quiver of type $\Bbb A_{n-1}$ 
$$
 b @>>> 3 @<<< 4 @<<< \cdots @<<< n
$$
and use the conification procedure of Corollary 1.3. In this way,
we write $\overline N$ as the direct sum of a representation of the form
$$
  0   @>>> N''_3 @<<< N''_4 @<<< \cdots @<<< N''_n.
$$
and a representation 
$$
 U_1\oplus U_2  @>>> N'_3 @<<< N'_4 @<<< \cdots @<<< N'_n
$$
which is conical with respect to $b$. In particular, the map $U_1\oplus U_2 \to N'_3$ 
has to be surjective. This means that
$N'_3 = \mu_1(U_1)+\mu_2(U_2)$, thus $N'_3$ is just the subspace $U_1+U_2$ of $N_3$.

It follows that $N$ is the direct sum $N = N'\oplus N''$ of
$$
{\beginpicture
\setcoordinatesystem units <1cm,.5cm>
\put{$N''\ = $} at -1.7 0
\put{$0$} at 0 1
\put{$0$} at 0 -1
\put{$N''_3$} at 1 0
\put{$N''_4$} at 2 0
\put{$\cdots$} at 3 0
\put{$N''_n$} at 4 0
\arr{0.3 0.7}{0.6 0.3}
\arr{0.3 -.7}{0.6 -.3}
\arr{1.7 0}{1.3 0}
\arr{2.7 0}{2.3 0}
\arr{3.7 0}{3.3 0}
\put{$\Biggl($} at -.5 0
\put{$\Biggr)$} at 4.5 0
\endpicture}.
$$
and 
$$
{\beginpicture
\setcoordinatesystem units <1cm,.5cm>
\put{$N'\ = $} at -1.7 0
\put{$U_1$} at 0 1
\put{$U_2$} at 0 -1
\put{$U_1\!+\!U_2$} at 1.2 0
\put{$N'_4$} at 2.6 0
\put{$\cdots$} at 3.6 0
\put{$N'_n$} at 4.6 0
\arr{0.3 0.7}{0.6 0.3}
\arr{0.3 -.7}{0.6 -.3}
\arr{2.3 0}{1.9 0}
\arr{3.3 0}{2.9 0}
\arr{4.3 0}{3.9 0}
\put{$\mu_1$} at 0.6 0.85
\put{$\mu_2$} at 0.6 -.95
\put{$\Biggl($} at -.6 0
\put{$\Biggr)$} at 5 0
\endpicture}
$$
with $\mu_i\:U_i \to U_1+U_2 = N'_3$ being the inclusion maps, 
and all the maps $N'_i \leftarrow N'_{i+1}$ for $3 \le i < n$ being injective.
    \medskip 
Given a representation $N$ in $\Cal N$, we denote by $\sigma^-N$ the representation
$$
{\beginpicture
\setcoordinatesystem units <1cm,.5cm>
\put{$\sigma^-N\ = $} at -2.2 0
\put{$\Biggl($} at -1 0
\put{$\Biggr)$} at 4.5 0
\put{$N_3/U_1$} at -.25 1
\put{$N_3/U_2$} at -.25 -1
\put{$N_3$} at 1 0
\put{$N_4$} at 2 0
\put{$\cdots$} at 3 0
\put{$N_n$} at 4 0
\arr{0.7 0.3}{0.3 0.7}
\arr{0.7 -.3}{0.3 -.7}
\arr{1.7 0}{1.3 0}
\arr{2.7 0}{2.3 0}
\arr{3.7 0}{3.3 0}
\put{$\pi_1$} at 0.65 0.85
\put{$\pi_2$} at 0.65 -.95
\endpicture}
$$
such that the maps $\pi_1,\pi_2$ are the canonical projections.
     \medskip 
Now, let us start with a representation $M$ of $Q(n).$
     \medskip 
{\bf (1)} We split off copies of $S(1)$ and $S(2)$ so that we can assume that $\pi_1$ and
$\pi_2$ are surjective. 
	
{\bf (2)} We look at the representation
$$
   M_1\oplus M_2 @<\left[\smallmatrix \pi_1\cr\pi_2\endsmallmatrix\right]<< M_3 @<<< M_4 @<<< \cdots @<<< M_n
$$ 
of the quiver $Q'$ of type $\Bbb A_{n-1}$ 
$$
 Q' \qquad \qquad a @<<< 3 @<<< 4 @<<< \cdots @<<< n.
$$
We use the conification procedure of Corollary 1.3 in order to write this representation
as the direct sum of a representation $(0 \leftarrow M''_3 \leftarrow \cdots \leftarrow M''_n)$
and a conical representation 
$(M_1\oplus M_2 \leftarrow M'_3 \leftarrow \cdots \leftarrow M'_n)$ of $(Q',a).$ 

It follows that $M$ is the direct sum $M = M'\oplus M''$, where
$$
{\beginpicture
\setcoordinatesystem units <1cm,.5cm>
\put{$M''\ = $} at -1.5 0
\put{$0$} at 0 1
\put{$0$} at 0 -1
\put{$M''_3$} at 1 0
\put{$M''_4$} at 2 0
\put{$\cdots$} at 3 0
\put{$M''_n$} at 4 0
\arr{0.6 0.3}{0.3 0.7}
\arr{0.6 -.3}{0.3 -.7}
\arr{1.7 0}{1.3 0}
\arr{2.7 0}{2.3 0}
\arr{3.7 0}{3.3 0}
\put{$\Biggl($} at -.45 0
\put{$\Biggr)$} at 4.5 0
\endpicture}.
$$
and 
$$
{\beginpicture
\setcoordinatesystem units <1cm,.5cm>
\put{$M' = $} at -1.4 0
\put{$M_1$} at 0 1
\put{$M_2$} at 0 -1
\put{$M'_3$} at 1 0
\put{$M'_4$} at 2 0
\put{$\cdots$} at 3 0
\put{$M'_n$} at 4 0
\arr{0.6 0.3}{0.3 0.7}
\arr{0.6 -.3}{0.3 -.7}
\arr{1.7 0}{1.3 0}
\arr{2.7 0}{2.3 0}
\arr{3.7 0}{3.3 0}
\put{$\pi_1$} at 0.6 0.85
\put{$\pi_2$} at 0.6 -.95
\put{$\Biggl($} at -.55 0
\put{$\Biggr)$} at 4.5 0
\endpicture}
$$
with the map $\left[\smallmatrix \pi_1\cr\pi_2\endsmallmatrix\right]\:M'_3 \to M_1\oplus M_2$ 
and all the maps $M'_{i+1} \to M'_i$ for $3 \le i < n$ being injective. Of course, the representation
$M''$ is a direct sum of thin representations, thus it remains to deal with $M'$. 
      \medskip 
{\bf (3)} We assume now that we deal with a representation $M$ 
$$
{\beginpicture
\setcoordinatesystem units <1cm,.5cm>
\put{$M_1$} at 0 1
\put{$M_2$} at 0 -1
\put{$M_3$} at 1 0
\put{$M_4$} at 2 0
\put{$\cdots$} at 3 0
\put{$M_n$} at 4 0
\arr{0.7 0.3}{0.3 0.7}
\arr{0.7 -.3}{0.3 -.7}
\arr{1.7 0}{1.3 0}
\arr{2.7 0}{2.3 0}
\arr{3.7 0}{3.3 0}
\put{$\pi_1$} at 0.6 0.85
\put{$\pi_2$} at 0.6 -.95
\endpicture}
$$
which is conical and that in addition the map
$\left[\smallmatrix \pi_1\cr\pi_2\endsmallmatrix\right]\:M_3 \to M_1\oplus M_2$ is injective. 
		    
Let $U_i$ be the kernel of $\pi_i$ with inclusion map $\mu_i\:U_i \to M_3.$ 
Since the map 
$\left[\smallmatrix \pi_1\cr\pi_2\endsmallmatrix\right]\:M_3 \to M_1\oplus M_2$  is injective, 
we have $U_1\cap U_2 = 0.$ 

Since $\pi_i$ is surjective, we may assume that
$M_i = M_3/U_i$ and that $\pi_i$ is the canonical projection map. The representation
$$
{\beginpicture
\setcoordinatesystem units <1cm,.5cm>
\put{$N\ = $} at -1.7 0
\put{$\Biggl($} at -.6 0
\put{$\Biggr)$} at 4.5 0
\put{$U_1$} at 0 1
\put{$U_2$} at 0 -1
\put{$M_3$} at 1 0
\put{$M_4$} at 2 0
\put{$\cdots$} at 3 0
\put{$M_n$} at 4 0
\arr{0.3 0.7}{0.6 0.3}
\arr{0.3 -.7}{0.6 -.3}
\arr{1.7 0}{1.3 0}
\arr{2.7 0}{2.3 0}
\arr{3.7 0}{3.3 0}
\put{$\mu_1$} at 0.6 0.85
\put{$\mu_2$} at 0.6 -.95
\endpicture}
$$
of $\sigma Q(n)$ belongs to $\Cal N$ and $\sigma^-N = M.$ Since $M$ is conical, also $N$
is conical. 

We use (0) in order to write $N = N'\oplus N''$ with
$$
{\beginpicture
\setcoordinatesystem units <1cm,.5cm>
\put{$N''\ = $} at -1.7 0
\put{$0$} at 0 1
\put{$0$} at 0 -1
\put{$N''_3$} at 1 0
\put{$N''_4$} at 2 0
\put{$\cdots$} at 3 0
\put{$N''_n$} at 4 0
\arr{0.3 0.7}{0.6 0.3}
\arr{0.3 -.7}{0.6 -.3}
\arr{1.7 0}{1.3 0}
\arr{2.7 0}{2.3 0}
\arr{3.7 0}{3.3 0}
\put{$\Biggl($} at -.5 0
\put{$\Biggr)$} at 4.5 0
\endpicture}.
$$
and 
$$
{\beginpicture
\setcoordinatesystem units <1cm,.5cm>
\put{$N'\ = $} at -1.7 0
\put{$U_1$} at 0 1
\put{$U_2$} at 0 -1
\put{$U_1\!+\!U_2$} at 1.2 0
\put{$N'_4$} at 2.6 0
\put{$\cdots$} at 3.6 0
\put{$N'_n$} at 4.6 0
\arr{0.3 0.7}{0.6 0.3}
\arr{0.3 -.7}{0.6 -.3}
\arr{2.3 0}{1.9 0}
\arr{3.3 0}{2.9 0}
\arr{4.3 0}{3.9 0}
\put{$\mu_1$} at 0.6 0.85
\put{$\mu_2$} at 0.6 -.95
\put{$\Biggl($} at -.6 0
\put{$\Biggr)$} at 5 0
\endpicture},
$$
where again the maps $\mu_i\:U_i \to U_1+U_2$ are the inclusion maps.
We have $M = \sigma^-N = \sigma^-N'\oplus \sigma^-N''.$ 
Since $M$ is conical, also the direct summands $\sigma^-N',\sigma^-N''$ of $M$ are conical. 
      \medskip 
      
Since $N''$ is conical and $N''_1 = N''_2 = 0$, we know that $N''$ is the direct sum
of indecomposable representations of $\sigma Q(n)$ with
dimension vectors of the form $\smallmatrix 0 \cr
                                 & 1 & \cdots & 1 & 0 & \cdots & 0 \cr
                               0 \endsmallmatrix$, thus 
$\sigma^-N''$ is the direct sum of indecomposable representations $W$ of $Q(n)$ with
dimension vectors of the form $\smallmatrix 1 \cr
                                 & 1 & \cdots & 1 & 0 & \cdots & 0 \cr
                               1 \endsmallmatrix$. In this way, we see that
$\sigma^-N''$ is a direct sum of thin representations.

It remains to look at $N'$ and $\sigma^-N'.$ Now $N'$ is a subrepresentation of $N$,
with $U_1 = N'_1 = N_1,\ U_2 = N'_2 = N_2$ and $U_1+U_2 = N'_3 \subseteq N_3.$
Since we  have $U_1\cap U_2 = 0,$ we see that $N'_3 = U_1\oplus U_2$ (and that the map
$\left[\matrix \mu_1 & \mu_2\endmatrix\right]\: N'_1 \oplus N'_2 \to N'_3$
is the identity map). It follows that
$$
{\beginpicture
\setcoordinatesystem units <1cm,.5cm>
\put{$\sigma^-N'\ = $} at -3.3 0
\put{$(U_1\!+\!U_2)/U_1$} at -1 1
\put{$(U_1\!+\!U_2)/U_2$} at -1 -1
\put{$U_1\!+\!U_2$} at 1.2 0
\put{$N'_4$} at 2.6 0
\put{$\cdots$} at 3.6 0
\put{$N'_n$} at 4.6 0
\arr{0.5 0.3}{0 0.7}
\arr{0.5 -.3}{0 -.7}
\arr{2.3 0}{1.9 0}
\arr{3.3 0}{2.9 0}
\arr{4.3 0}{3.9 0}
\put{$\pi_1$} at 0.4 0.85
\put{$\pi_2$} at 0.4 -.95
\put{$\Biggl($} at -2.2 0
\put{$\Biggr)$} at 5 0
\endpicture},
$$
belongs to $\check{\Cal B}.$ This completes the proof of Theorem 2.1.
	\bigskip 
{\bf Remark.} We have denoted by $\Cal N$ the category of representations 
$N$ of $\sigma Q(n)$ such that
the maps $\mu_1,\mu_2$ are inclusion maps. Similarly, we may denote 
by $\Cal M$ the category of representations $M$ of $Q(n)$ such that
the maps $\pi_1,\pi_2$ are surjective.
The functor $\sigma^-$ maps $\Cal N$ to $\Cal M$. Similarly, we may consider 
a functor $\sigma^+\:\Cal M \to \Cal N$ by using the kernels of the maps $\pi_1,\pi_2$.
Clearly, $\sigma^-$ is an equivalence of categories with inverse $\sigma^+$.
	 \bigskip 
Let us now investigate the subcategory $\check{\Cal B}$.
    \medskip 
{\bf 2.2. Theorem.} 
{\it The category $\check{\Cal B}$ is equivalent to $\Rep Q(n\!-\!1)$; it is the image of the
fully faithful functor 
$$
 \eta\:\Rep Q(n\!-\!1) \to \Rep Q(n)
$$ 
defined as follows:
$$
{\beginpicture
\setcoordinatesystem units <1cm,.5cm>
\put{\beginpicture
\put{$\Biggl($} at -.5 0
\put{$\Biggr)$} at 4.7 0
\put{$M_1$} at 0 1
\put{$M_2$} at 0 -1
\put{$M_3$} at 1 0
\put{$M_4$} at 2 0
\put{$\cdots$} at 3 0
\put{$M_{n-1}$} at 4.2 0
\arr{0.7 0.3}{0.3 0.7}
\arr{0.7 -.3}{0.3 -.7}
\arr{1.7 0}{1.3 0}
\arr{2.7 0}{2.3 0}
\arr{3.7 0}{3.3 0}
\put{$\pi_1$} at 0.6 0.85
\put{$\pi_2$} at 0.6 -.95
\endpicture} at 0 0 
\put{\beginpicture
\put{$\Biggl($} at -.95 0
\put{$\Biggr)$;} at 6.2 0
\put{$M_1$} at -.5 1
\put{$M_2$} at -.5 -1
\put{$M_1\!\oplus\! M_2$} at 1 0
\put{$M_3$} at 2.5 0
\put{$M_4$} at 3.5 0
\put{$\cdots$} at 4.5 0
\put{$M_{n-1}$} at 5.7 0
\arr{0.2 0.3}{-.2 0.7}
\arr{0.2 -.3}{-.2 -.7}
\arr{2.2 0}{1.8 0}
\arr{3.2 0}{2.8 0}
\arr{4.2 0}{3.8 0}
\arr{5.2 0}{4.8 0}
\put{$\left[\smallmatrix\pi_1\cr\pi_2\endsmallmatrix\right]$} at 2 0.7
\put{$\epsilon_1$} at 0.1 0.85
\put{$\epsilon_2$} at 0.1 -.95
\endpicture} at 7.3 0
\put{$\mapsto$} at 3.2 0 
\endpicture}
$$
here, for $i=1,2$, the map $\epsilon_i$ is the canonical projection $M_1\oplus M_2 \to M_i$.}
      \medskip
Proof. Clearly, the image of $\eta$ is just $\check{\Cal B}$. There is an inverse functor from 
$\Rep Q(n)$ to $\Rep Q(n\!-\!1)$ which deletes the vector space at the branching vertex:
$$
{\beginpicture
\setcoordinatesystem units <1cm,.5cm>
\put{\beginpicture
\put{$\Biggl($} at -.95 0
\put{$\Biggr)$} at 4.9 0
\put{$M_1$} at -.5 1
\put{$M_2$} at -.5 -1
\put{$M_3$} at .5 0
\put{$M_4$} at 1.8 0
\put{$M_{n}$} at 4.5 0
\arr{0.2 0.3}{-.2 0.7}
\arr{0.2 -.3}{-.2 -.7}
\arr{1.4 0}{.9 0}

\arr{2.7 0}{2.1 0}
\put{$\cdots$} at 3 0
\arr{3.9 0}{3.3 0}
\put{$\gamma$} at 1.15 0.3
\put{$\pi_1$} at 0.1 0.85
\put{$\pi_2$} at 0.1 -.95
\endpicture} at 0 0

\put{\beginpicture
\put{$\Biggl($} at -.5 0
\put{$\Biggr).$} at 4.5 0
\put{$M_1$} at 0 1
\put{$M_2$} at 0 -1
\put{$M_4$} at 1.4 0
\put{$\cdots$} at 2.7 0
\put{$M_{n}$} at 4 0
\arr{0.9 0.3}{0.3 0.7}
\arr{0.9 -.3}{0.3 -.7}
\arr{2.3 0}{1.8 0}
\arr{3.5 0}{3 0}
\put{$\pi_1\gamma$} at 0.8 0.85
\put{$\pi_2\gamma$} at 0.8 -.95
\endpicture} at 7 0 
\put{$\mapsto$} at 3.7 0 
\endpicture}
$$
This shows that $\eta$ is fully faithful.
     \bigskip 
{\bf 2.3. Corollary.} {\it Any representation of $Q(n)$ is the direct sum of finite-dimensional 
indecomposable representations.}
	       \medskip
Proof, by induction. For $n = 3$ we deal with a quiver of type $\Bbb A_3$, thus the assertion 
has been shown in Part 1. Thus, assume that $n\ge 4$. According to Theorem 2.1, any representation $M$ of $Q(n)$
is the direct sum of thin representations and a representation $M'$ in $\check{\Cal B}$. According to Theorem 2.2,
the functor $\eta$ provides an equivalence between $\check{\Cal B}$ and $\Rep Q(n\!-\!1)$,
and the representations in $\Cal B$ correspond under $\eta$ to the finite-dimensional
representations of $Q(n\!-\!1)$. By induction, $\eta^{-1}(M')$ is a direct sum of finite-dimensional
representations of $Q(n\!-\!1)$, thus $M'$ is the direct sum of finite-dimensional representations of $Q(n)$.

		\bigskip 
{\bf 2.4. Twin representations.}
We call an indecomposable representation $M$ of $Q(n)$
a {\it twin} representation provided $M_1\neq 0, M_2 \neq 0$.

Here are typical examples of twin representations: First of all, we will call the thin twin representations 
{\it $0$-twin} representations; they are the indecomposable representations with dimension vector of the form 
$\smallmatrix 1 \cr
                & 1 & \cdots & 1 & 0 & \cdots & 0 \cr
              1 \endsmallmatrix
$  
(and dimension at least 3). Second, the following representations
$$
{\beginpicture
\setcoordinatesystem units <1.3cm,1cm>
\multiput{$k$} at 0 1  0 -1  4 0  6 0   /
\multiput{$k^2$} at 1 0    3 0  /
\multiput{$0$} at 7 0  9 0   /
\multiput{$\cdots$} at 2 0  5 0  8 0 /
\put{$\delta$} at 3.5 0.25 
\arr{0.8 0.2}{0.2 0.8}
\arr{0.8 -.2}{0.2 -.8}
\arr{1.7 0}{1.3 0}
\arr{2.7 0}{2.3 0}
\arr{3.7 0}{3.3 0}
\arr{4.7 0}{4.3 0}
\arr{5.7 0}{5.3 0}
\arr{6.7 0}{6.3 0}
\arr{7.7 0}{7.3 0}
\arr{8.7 0}{8.3 0}
\put{$\pi_1$} at 0.6 0.65
\put{$\pi_2$} at 0.6 -.7
\multiput{$1$} at 1.5 0.25  2.5 0.25  4.5 0.25  5.5 0.25 /
 \endpicture}
$$
with $r\ge 1$ vector spaces of the form $k^2$, and at least three vector spaces of the
form $k$ will be called the {\it $r$-twin} representations;
here $\delta\:k\to k^2$ denotes the diagonal embedding, and $\pi_i\:k^2\to k$
with $i=1,2$ are the two canonical projections. In order to see that these $r$-twin representations
are twin representations, we need to know that they are indecomposable. But this is easy to check
(it is also a direct consequence of the following Lemma).
    \medskip
{\bf Lemma. (a)} {\it Let $n\ge 4$ and $1\le r \le n-3$. The functor 
$\eta$ provides a bijection between the $(r-1)$-twin representations of $Q(n\!-\!1)$
and the $r$-twin representations of $Q(n)$.} 

{\bf (b)} {\it A twin representation of $Q(n)$ is an $r$-twin representation for some $r\ge 0$.} 
     \medskip
Proof.  (a) is straightforward. 
(b) Let $M$ be a twin representation of $Q(n)$. If $M$ is thin, then it is a $0$-twin representation, by
definition. Thus, we many assume that $M$ is not thin. According to Theorem 2.1, $M$ belongs to
$\Cal B(n)$, thus it is of the form $\eta(N)$ for some indecomposable representation of $Q(n\!-\!1)$.
By the definition of $\eta$, we have $N_i = M_i \neq 0,$ for $i=1,2$, thus $N$
is a twin representation. By induction, $N$ is an $s$-twin  representation for some $s\ge 0.$ 
According to (a), $M = \eta(N)$ is an $(s+1)$-twin representation. 
	  \medskip 
{\bf Corollary.} {\it An indecomposable representation of $Q(n)$ is either thin or an
$r$-twin representation with $1\le r \le n-3$. Any representation of $Q(n)$ is the direct sum
of thin and twin representations.}
    \bigskip

From now on, we mainly will deal with finite-dimensional representations.
Thus, vector spaces and representations will be assumed to be finite-dimensional
unless otherwise stated.
       \bigskip
       {\bf 2.5. The position of some subcategories in the
       Auslander-Reiten quiver of $Q(n)$.}
Given a class $\Cal U$ of representations of a quiver, we denote by $\add \Cal U$ the class
of all direct summands of finite direct sums of representations in $\Cal U$. 
If $\Cal U', \Cal U''$ are subcategories, we write $\Cal U'\vee \Cal U''$ for the class of representations
of the form $C'\oplus C''$ with $C'\in \Cal U',\ C''\in \Cal U''$.

Here are the subcategories which we are interested in:
$$
\align
 \Cal V &= \add S(1)\vee \add S(2)\cr
 \Cal X &= \add\{M\mid \text{$M$ is a $0$-twin representation}\}\cr
 \Cal Y &= \{M\mid \bmatrix \pi_1\cr\pi_2\endbmatrix\:M_3 \to M_1\oplus M_2 \
   \text{bijective, and}\ 
   M_{i+1} \to M_i\ 
   \text{for $3\le i < n$ injective} \} \cr
 \Cal X' &= \add\{M\mid \text{$M$ indecomposable thin representation with $M_1 = 0 = M_2, M_3\neq 0$}\}\cr
 \Cal W &= \{\text{representations with support in $[4,n]$}\}
\endalign
$$

For the notion of an Auslander-Reiten quiver, we refer to [R3, ARS]. 
For $n = 6$, let us draw the Auslander-Reiten quiver of $Q(n)$.
We have encircled by solid lines the parts $\Cal X,\ \Cal Y$ and $\Cal X'$, and by
dotted lines the parts $\Cal V$ and $\Cal W$.
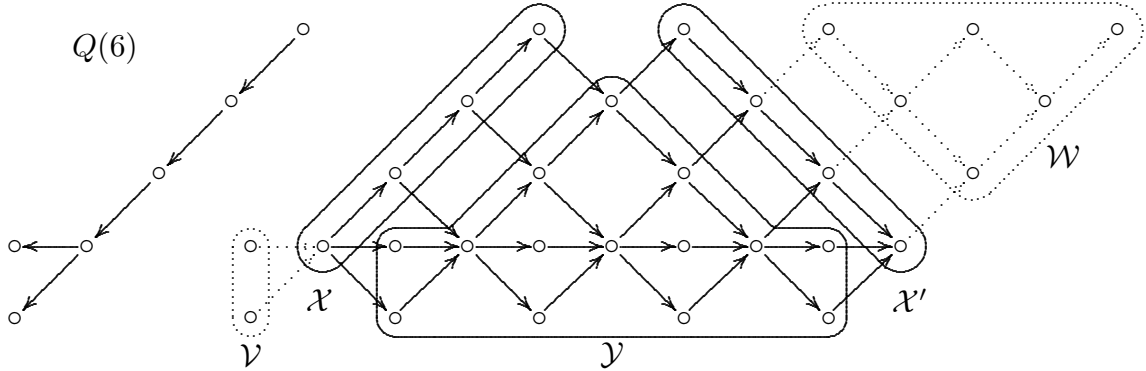
$$
{\beginpicture
\setcoordinatesystem units <.2cm,.2cm>
\put{\beginpicture
\multiput{} at 0 13  0 -7 /
\multiput{$\circ$} at 0 0  4 4  8 8  12 12  -4 0  -4 -4 /
\put{$Q(6)$} at 1 11
\arr{3.5 3.5}{0.5 0.5}
\arr{7.5 7.5}{4.5 4.5}
\arr{11.5 11.5}{8.5 8.5}
\arr{-0.5 0}{-3.5 0}
\arr{-0.5 -.5}{-3.5 -3.5}
\endpicture} at -25 0
\put{\beginpicture
\multiput{} at 0 13  0 -7 /
\multiput{$\circ$} at 0 0  4 0  8 0  12 0  16 0  20 0  24 0  28 0  32 0
       4 -4  12 -4  20 -4  28 -4
              4 4   12 4   20 4   28  4
	                8 8    16 8   24 8
			             12 12  20 12 /
				     \plot 1 -1  13 11 /
				     \plot -1 1  11 13 /
				     \circulararc -180 degrees from 1 -1 center at 0 0
				     \circulararc 180 degrees from 13 11 center at 12 12
				     \plot 31 -1  19 11 /
				     \plot 33 1  21 13 /
				     \circulararc -180 degrees from 19 11 center at 20 12
				     \circulararc 180 degrees from 31 -1 center at 32 0
				     \put{$\Cal V$} at -4 -6
				     \put{$\Cal X$} at -.2 -3
				     \put{$\Cal X'$} at 32.5 -3
				     \plot 17 9  25 1 /
				     \plot 15 9  7 1 /
				     \put{$\Cal Y$} at 16 -6
				     \circulararc 90 degrees from 17 9 center at 16 8
				     \plot 4 1  7 1 /
				     \circulararc 90 degrees from 4 1 center at 4 0
				     \plot 3 0  3 -4 /
				     \circulararc 90 degrees from 3 -4 center at 4 -4
				     \plot 4 -5  28 -5 /

				     \circulararc 90 degrees from 28 -5 center at 28 -4
				     \plot 29 -4  29 0 /
				     \circulararc 90 degrees from 29 0 center at 28 0
				     \plot 28 1  25 1 /
				     \arr{0.5 0}{3.5 0}
				     \arr{4.5 0}{7.5 0}
				     \arr{8.5 0}{11.5 0}
				     \arr{12.5 0}{15.5 0}
				     \arr{16.5 0}{19.5 0}
				     \arr{20.5 0}{23.5 0}
				     \arr{24.5 0}{27.5 0}
				     \arr{28.5 0}{31.5 0}

				     \arr{0.5 -.5}{3.5 -3.5}
				     \arr{4.3 -3.5}{7.5 -.5}
				     \arr{8.5 -.5}{11.5 -3.5}
				     \arr{12.3 -3.5}{15.5 -.5}
				     \arr{16.5 -.5}{19.5 -3.5}
				     \arr{20.3 -3.5}{23.5 -.5}
				     \arr{24.5 -.5}{27.5 -3.5}
				     \arr{28.3 -3.5}{31.5 -.5}

				     \arr{0.5 .5}{3.5 3.5}
				     \arr{4.3 3.5}{7.5 .5}
				     \arr{8.5 .5}{11.5 3.5}
				     \arr{12.3 3.5}{15.5 .5}
				     \arr{16.5 .5}{19.5 3.5}
				     \arr{20.3 3.5}{23.5 .5}
				     \arr{24.5 .5}{27.5 3.5}
				     \arr{28.3 3.5}{31.5 .5}

				     \arr{4.5 4.5}{7.5 7.5}
				     \arr{8.3 7.5}{11.5 4.5}
				     \arr{12.5 4.5}{15.5 7.5}
				     \arr{16.3 7.5}{19.5 4.5}
				     \arr{20.5 4.5}{23.5 7.5}
				     \arr{24.3 7.5}{27.5 4.5}

				     \arr{8.5 8.5}{11.5 11.5}
				     \arr{12.3 11.5}{15.5 8.5}
				     \arr{16.5 8.5}{19.5 11.5}
				     \arr{20.3 11.5}{23.5 8.5}

				     \multiput{$\circ$} at 28 12  32 8  36 4  36 12  40 8  44 12  -4 0  -4 -4 /
				     
\setdots <.7mm>
\plot 35 3  27 11 /
\plot 37 3  45 11 /
\plot 44 13.5  28 13.5 /
\circulararc 90 degrees from 35 3 center at 36 4
\circulararc 135 degrees from 45 11 center at 44 12
\circulararc -135 degrees from 27 11 center at 28 12
\put{$\Cal W$} at 41 5

\plot -5 0  -5 -4 /
\plot -3 0  -3 -4 /
\circulararc -180 degrees from -5 0 center at -4 0
\circulararc 180 degrees from -5 -4 center at -4 -4

\setdots <1mm>
\arr{-3.5 0}{-.5 0}
\arr{-3.5 -3.5}{-.5 -.5}

\arr{24.5 8.5}{27.5 11.5}
\arr{28.5 4.5}{31.5 7.5}
\arr{32.5 0.5}{35.5 3.5}
\arr{28.5 11.5}{31.5 8.5}
\arr{32.5 7.5}{35.5 4.5}

\arr{32.5 8.5}{35.5 11.5}
\arr{36.5 4.5}{39.5 7.5}

\arr{36.5 11.5}{39.5 8.5}

\arr{40.5 8.5}{43.5 11.5}

\arr{4.3 3.5}{7.5 .5}
\endpicture} at 4 0
\endpicture}
$$
These are the building blocks of $\rep Q(n)$. 
The dotted part on the left are the two simple representations $S(1)$ and $S(2)$, this is $\Cal V$.
The dotted part 
on the right shows the indecomposable representations with support in $[4,n]$, this is $\Cal W$.
We should stress that $\Cal X \vee \Cal Y \vee \Cal X'$ is just the category $\Conic Q(n)$. 
   \medskip 

Let us introduce the following
notation:
$$
  \Cal A = \Cal X \vee \Cal W, \quad
  \Cal B = \Cal Y \vee \Cal W, \quad
  \Cal A' = \Cal X' \vee \Cal W.
$$
(thus, $\Cal W$ belongs to $\Cal A, \Cal B, \Cal A'$).
       \medskip 
The category $\Cal A'$ are the  representations with support in $[3,n]$, thus it
is the category of representations of a quiver of type $\Bbb A_{n-2}$.

The category $\Cal A$ is equivalent to the category of  representations of a quiver 
of type $\Bbb A_{n-2}$, since the functor $\sigma^+$ (introduced in the proof of Theorem 2.1)
provides an equivalence between
$\Cal A$ and $\Cal A'.$

The shape of the category $\Cal B = \Cal Y \vee \Cal W$ illustrates nicely that 
$\Cal B = \Cal B(n)$ may be identified 
with  the category of representations of a quiver of
type $\Bbb D_{n-1}$, as we have shown in Theorem 2.2.
     \bigskip
{\bf Analysis of the proof of Theorem 2.1.}
Let us use the Auslander-Reiten quiver shown above  to
describe the steps in the proof of Theorem 2.1. In step (1), we have split off representations in
$\Cal V$. 
Step (2) was devoted to splitting off representations in $\Cal X'\vee\Cal W$.
In the final step (3), we have split off representations
in $\Cal X$. What remains turns out to be a direct sum of representations
in $\Cal Y$ (and $\Cal Y \subset \Cal B$).
   \bigskip
{\bf 2.6. Antichains in posets and additive categories.}
Given a poset $P$, a subset $A$ of $P$ is called an {\it antichain} provided $a\le a'$
for $a,a'$ in $A$ implies $a = a'$. 
An {\it antichain} in a $k$-category $\Cal C$ is a set of pairwise ($\Hom$-)orthogonal objects
with endomorphism rings being division rings.
(Starting with a poset $P$, one may consider its linearization $kP$, see for example [R5, N1.8],
this is an additive category and
the antichains in the poset $P$ are just the antichains in the additive category $kP$.)
An antichain set of cardinality 2 or 3 will be called
an {\it antichain pair} or an {\it antichain triple,} respectively. 
   \medskip 
{\bf Simplification.} Let us assume now that $\Cal C$ is an abelian category. 
Starting with an antichain $A$ in $\Cal C$, we can consider the full subcategory 
$\Cal E(A)$ 
of all objects with a filtration with factors in the antichain: this is a thick subcategory 
(an exact abelian subcategory which is closed under extensions).
Starting with an antichain
$A = \{A_1,\dots,A_t\}$ in an abelian $k$-category $\Cal C$ such that the endomorphism ring
of any object $A(i)$ is equal to $k$,  its 
{\it $\Ext$-quiver} $Q(A)$ has $t$ vertices, and the number of arrows $j\to i$ is the $k$-dimension
of $\Ext^1(A_j,A_i).$ In case the $\Ext$-quiver is finite and 
acyclic, the category $\Cal E(A)$ is
equivalent to the category $\mo \Lambda$ of all finite-dimensional $\Lambda$-modules, 
where $\Lambda$ is a finite-dimensional
$k$-algebra, see [R1]. 
Note that the simple objects in the category $\Cal E(A)$
are just the elements of $A$, thus to focus the attention to $\Cal E(A)$ is called the process of 
{\it simplification.} 
     \medskip
Let us return to the quiver $Q(n)$. It is easy to see that
the representations 
$$
  M(\langle 1,3\rangle), M(\langle 2,3\rangle), S(4),\dots, S(n)
$$ 
form an antichain and that $\Cal E(M(\langle 1,3\rangle), M(\langle 2,3\rangle), S(4),\dots, S(n))$
are the finite-dimensional representations in $\Cal B(n)$.

The $\Ext$-quiver of this antichain looks as follows:
$$
{\beginpicture
\setcoordinatesystem units <1.5cm,.6cm>
\put{$M(\langle 1,3\rangle)$} at -.2 1
\put{$M(\langle 2,3\rangle)$} at -.2 -1
\put{$S(4)$} at 1 0
\put{$\cdots$} at 2 0
\put{$S(n)$} at 3 0
\arr{0.7 0.2}{0.3 0.7}
\arr{0.7 -.2}{0.3 -.7}
\arr{1.7 0}{1.3 0}
\arr{2.7 0}{2.3 0}

\endpicture}
$$
thus, it is the quiver $Q(n\!-\!1)$. This shows (again) that the category of 
finite-dimensional representations in $\Cal B(n)$ is equivalent to $\rep Q(n\!-\!1)$.
		   \bigskip\bigskip 
{\bf 2.7. The conical representations.}
Let us recall that the category $\Conic Q(n)$ is the intrinsic object to look at,
we have $\Conic Q(n) = \Cal X\vee \Cal Y \vee \Cal X'.$
As we have mentioned in Proposition 1.6, the category of conical representations does not change,
if we change the orientation of the quiver. 
Here is the Auslander-Reiten quiver of a second quiver $Q'$ of type $\Bbb D_6$, it should facilitate
the reader to compare the embeddings of the category of 
conical representations into the representation
categories.
$$
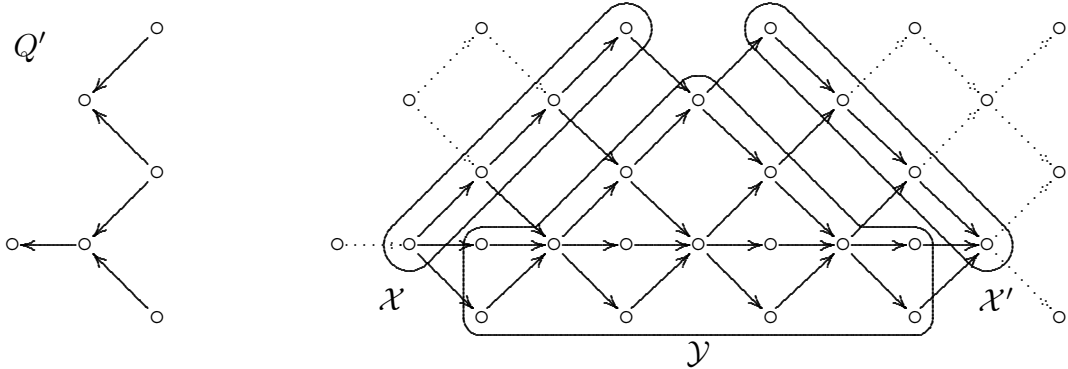
{\beginpicture
\setcoordinatesystem units <.2cm,.2cm>
\put{\beginpicture
\multiput{} at 0 13  0 -7 /
\multiput{$\circ$} at 0 0  4 4  0 8  4 12  -4 0  4 -4 /
\put{$Q'$} at -3 11
\arr{3.5 3.5}{0.5 0.5}
\arr{3.5 4.5}{.5 7.5}
\arr{3.5 11.5}{.5 8.5}
\arr{-.5 0}{-3.5 0}
\arr{3.5 -3.5}{0.5 -.5}
\endpicture} at -30 0
\put{\beginpicture
\multiput{} at 0 13  0 -7 /
\multiput{$\circ$} at 0 0  4 0  8 0  12 0  16 0  20 0  24 0  28 0  32 0
       4 -4  12 -4  20 -4  28 -4
              4 4   12 4   20 4   28  4
	                8 8    16 8   24 8
			             12 12  20 12 /
				     \plot 1 -1  13 11 /
				     \plot -1 1  11 13 /
				     \circulararc -180 degrees from 1 -1 center at 0 0
				     \circulararc 180 degrees from 13 11 center at 12 12
				     \plot 31 -1  19 11 /
				     \plot 33 1  21 13 /
				     \circulararc -180 degrees from 19 11 center at 20 12
				     \circulararc 180 degrees from 31 -1 center at 32 0
				     \put{$\Cal X$} at -1 -3
				     \put{$\Cal X'$} at 32.5 -3
				     \plot 17 9  25 1 /
				     \plot 15 9  7 1 /
				     \put{$\Cal Y$} at 16 -6
				     \circulararc 90 degrees from 17 9 center at 16 8
				     \plot 4 1  7 1 /
				     \circulararc 90 degrees from 4 1 center at 4 0
				     \plot 3 0  3 -4 /
				     \circulararc 90 degrees from 3 -4 center at 4 -4
				     \plot 4 -5  28 -5 /

				     \circulararc 90 degrees from 28 -5 center at 28 -4
				     \plot 29 -4  29 0 /
				     \circulararc 90 degrees from 29 0 center at 28 0
				     \plot 28 1  25 1 /
				     \arr{0.5 0}{3.5 0}
				     \arr{4.5 0}{7.5 0}
				     \arr{8.5 0}{11.5 0}
				     \arr{12.5 0}{15.5 0}
				     \arr{16.5 0}{19.5 0}
				     \arr{20.5 0}{23.5 0}
				     \arr{24.5 0}{27.5 0}
				     \arr{28.5 0}{31.5 0}

				     \arr{0.5 -.5}{3.5 -3.5}
				     \arr{4.3 -3.5}{7.5 -.5}
				     \arr{8.5 -.5}{11.5 -3.5}
				     \arr{12.3 -3.5}{15.5 -.5}
				     \arr{16.5 -.5}{19.5 -3.5}
				     \arr{20.3 -3.5}{23.5 -.5}
				     \arr{24.5 -.5}{27.5 -3.5}
				     \arr{28.3 -3.5}{31.5 -.5}

				     \arr{0.5 .5}{3.5 3.5}
				     \arr{4.3 3.5}{7.5 .5}
				     \arr{8.5 .5}{11.5 3.5}
				     \arr{12.3 3.5}{15.5 .5}
				     \arr{16.5 .5}{19.5 3.5}
				     \arr{20.3 3.5}{23.5 .5}
				     \arr{24.5 .5}{27.5 3.5}
				     \arr{28.3 3.5}{31.5 .5}

				     \arr{4.5 4.5}{7.5 7.5}
				     \arr{8.3 7.5}{11.5 4.5}
				     \arr{12.5 4.5}{15.5 7.5}
				     \arr{16.3 7.5}{19.5 4.5}
				     \arr{20.5 4.5}{23.5 7.5}
				     \arr{24.3 7.5}{27.5 4.5}

				     \arr{8.5 8.5}{11.5 11.5}
				     \arr{12.3 11.5}{15.5 8.5}
				     \arr{16.5 8.5}{19.5 11.5}
				     \arr{20.3 11.5}{23.5 8.5}

				     \multiput{$\circ$} at -4 0  0 8  4 12
				     28 12  32 8  36 4  36 12     36 -4 /

				     \setdots <1mm>
				     \arr{32.5 -.5}{35.5 -3.5}

				     \arr{24.5 8.5}{27.5 11.5}
				     \arr{28.5 4.5}{31.5 7.5}
				     \arr{32.5 0.5}{35.5 3.5}
				     \arr{28.5 11.5}{31.5 8.5}
				     \arr{32.5 7.5}{35.5 4.5}

				     \arr{32.5 8.5}{35.5 11.5}
				     \arr{4.3 3.5}{7.5 .5}

				     \arr{-3.5 0}{-.5 0}
				     \arr{0.5 7.5}{3.5 4.5}
				     \arr{4.5 11.5}{7.5 8.5}
				     \arr{0.5 8.5}{3.5 11.5}
				     \endpicture} at 4 0
				     \endpicture}
				     $$
					\bigskip
{\bf Change of orientation: Antichains.} For the quiver $Q(n)$, we have exhibited above
an antichain $A$ such that the objects of $\Cal E(A)$ are the finite-dimensional representations inside
$\Cal B(n)$. Similar antichains do exist for any orientation, but we have to distinguish two
cases, namely whether the two short arms have the same orientation or not.
In case they have the same orientation, the antichain to be considered consists (as in the case
$Q(n))$ of the representations:
$$
  M(\langle 1,3\rangle), M(\langle 2,3\rangle), S(4),\dots, S(n).
$$ 
Now assume that the short arms have different orientation. Let us look at the quiver $Q'(n)$ 
given as follows:
$$
{\beginpicture
\setcoordinatesystem units <1cm,.5cm>
\put{$Q'(n)$} at -2 0 
\put{$1$} at 0 1
\put{$2$} at 0 -1
\put{$3$} at 1 0
\put{$4$} at 2 0
\put{$\cdots$} at 3 0
\put{$n$} at 4 0
\put{$\mu$} at 0.6 0.8
\put{$\pi$} at 0.6 -0.8
\put{$\gamma$} at 1.55 0.3

\arr{0.3 0.7}{0.6 0.3}
\arr{0.6 -.3}{0.3 -.7}
\arr{1.7 0}{1.3 0}
\arr{2.7 0}{2.3 0}
\arr{3.7 0}{3.3 0}
\endpicture}
$$
Then the antichain to be considered is the following:
$$
   M(\langle 1,2,3\rangle), S(3), S(4),\dots, S(n).
$$
Of course, its $\Ext$-quiver 
$$
{\beginpicture
\setcoordinatesystem units <1.5cm,.6cm>
\put{$M(\langle 1,2,3\rangle)$} at -.4 1
\put{$S(3)$} at -.05 -1
\put{$S(4)$} at 1 0
\put{$\cdots$} at 2 0
\put{$S(n)$} at 3 0
\arr{0.7 0.2}{0.3 0.7}
\arr{0.7 -.2}{0.3 -.7}
\arr{1.7 0}{1.3 0}
\arr{2.7 0}{2.3 0}
\endpicture}
$$
is again of type $\Bbb D_{n-1}$ (here, the orientation of the two short arms is determined
by the orientation of $\gamma$). Note that $\Cal E(M(\langle 1,2,3\rangle), S(3), S(4),\dots, S(n))$ is the full
subcategory of all representations of $Q'(n)$ such that $\pi\mu$ is bijective. 
	    \medskip 
{\bf 2.8. Perpendicular categories.} Given a class $\Cal U$ of representation of a quiver $Q$, we denote
by $\Cal U^\perp$ the full subcategory of $\rep Q$ given by all representations $M$
with $\Hom(C,M) = 0 = \Ext^1(C,M)$ for all $C$ in $\Cal U$. 
Similarly,  we denote
by ${}^\perp \Cal U$ the full subcategory of $\rep Q$ given by all representations $M$
with $\Hom(M,C) = 0 = \Ext^1(M,C)$ for all $C$ in $\Cal U$.

The reduction from $\Bbb D_n$-quivers to $\Bbb D_{n-1}$-quivers can be described very
well using such perpendicular categories: 
Let $X$ be the indecomposable representation in $\Cal X$ which is a sink in $\Cal X$
(thus any non-zero homomorphism $X \to X'$ in $\Cal X$ is a split monomorphism).
Let $Z$ be the indecomposable representation in $\Cal X'$ which is a source in $\Cal X'$
(thus any non-zero homomorphism $Z' \to Z$ in $\Cal X'$ is a split epimorphism).
Then $\tau Z = X$. For any orientation, we may look at
$$
 \Cal B(n) = {}^\perp X = Z^\perp
           = \{M\mid \Hom(M,X) = 0 = \Hom(Z,M)\}
$$

Thus, in order to find $\Cal B(n)$, we have to delete on the one hand 
the indecomposable representations $M$ with 
$\Hom(M,X) \neq 0$, these are the indecomposable representations in $\Cal X$ as well as the simple projective
modules of the form $S(1),S(2)$, and on the other 
hand the indecomposable representations $M$ with 
$\Hom(Z,M) \neq 0$, these are the indecomposable representations in $\Cal X'$ as well as the simple injective
modules of the form $S(1),S(2)$. Thus $\Cal B(n)$ is obtained from $\rep Q(n)$ by deleting
$\Cal X, \Cal X'$ as well as $S(1)$ and $S(2).$ 

Here are the Auslander-Reiten quivers of $Q(6)$ and $Q'$ with $X,Z$ encircled
(the subcategories $\Cal X$ and $\Cal X'$ as well as the two simple representations in $\Cal V$
are indicated by dotted lines):
$$
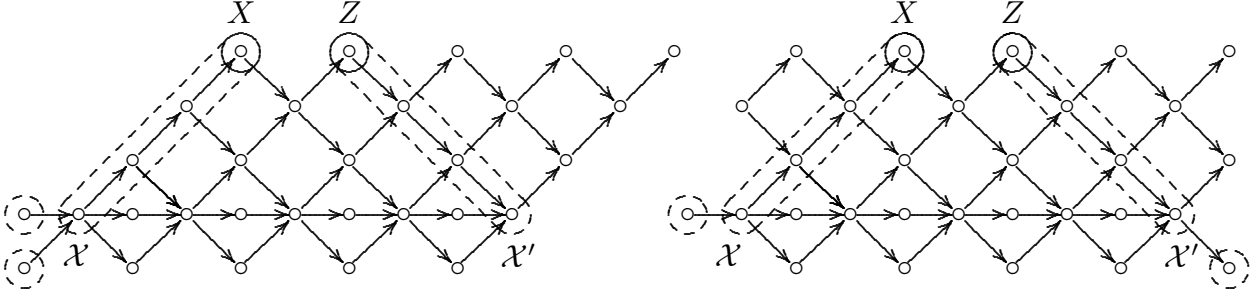
{\beginpicture
\setcoordinatesystem units <.15cm,.15cm>
\put{\beginpicture
\multiput{} at 0 13  0 -7 /
\multiput{$\circ$} at 0 0  4 0  8 0  12 0  16 0  20 0  24 0  28 0  32 0
       4 -4  12 -4  20 -4  28 -4
              4 4   12 4   20 4   28  4
	                8 8    16 8   24 8
			             12 12  20 12 /
\setdashes <1mm>
\plot 1 -1  13 11 /
\plot -1 1  11 13 /
\circulararc -180 degrees from 1 -1 center at 0 0
\plot 31 -1  19 11 /
\plot 33 1  21 13 /
\circulararc 180 degrees from 31 -1 center at 32 0
\put{$\Cal X$} at -.2 -3
\put{$\Cal X'$} at 32.5 -3

\circulararc 360 degrees from -3 1 center at -4 0
\circulararc 360 degrees from -3 -3 center at -4 -4

\setsolid
\circulararc 360 degrees from 13 11 center at 12 12
\circulararc 360 degrees from 19 11 center at 20 12
\put{$X$} at 12 15
\put{$Z$} at 20 15

\arr{0.5 0}{3.5 0}
\arr{4.5 0}{7.5 0}
\arr{8.5 0}{11.5 0}
\arr{12.5 0}{15.5 0}
\arr{16.5 0}{19.5 0}
\arr{20.5 0}{23.5 0}
\arr{24.5 0}{27.5 0}
\arr{28.5 0}{31.5 0}

\arr{0.5 -.5}{3.5 -3.5}
\arr{4.3 -3.5}{7.5 -.5}
\arr{8.5 -.5}{11.5 -3.5}
\arr{12.3 -3.5}{15.5 -.5}
\arr{16.5 -.5}{19.5 -3.5}
\arr{20.3 -3.5}{23.5 -.5}
\arr{24.5 -.5}{27.5 -3.5}
\arr{28.3 -3.5}{31.5 -.5}

\arr{0.5 .5}{3.5 3.5}
\arr{4.3 3.5}{7.5 .5}
\arr{8.5 .5}{11.5 3.5}
\arr{12.3 3.5}{15.5 .5}
\arr{16.5 .5}{19.5 3.5}
\arr{20.3 3.5}{23.5 .5}
\arr{24.5 .5}{27.5 3.5}
\arr{28.3 3.5}{31.5 .5}

\arr{4.5 4.5}{7.5 7.5}
\arr{8.3 7.5}{11.5 4.5}
\arr{12.5 4.5}{15.5 7.5}
\arr{16.3 7.5}{19.5 4.5}
\arr{20.5 4.5}{23.5 7.5}
\arr{24.3 7.5}{27.5 4.5}

\arr{8.5 8.5}{11.5 11.5}
\arr{12.3 11.5}{15.5 8.5}
\arr{16.5 8.5}{19.5 11.5}
\arr{20.3 11.5}{23.5 8.5}

\multiput{$\circ$} at 28 12  32 8  36 4  36 12  40 8  44 12  -4 0  -4 -4 /

\arr{-3.5 0}{-.5 0}
\arr{-3.5 -3.5}{-.5 -.5}

\arr{24.5 8.5}{27.5 11.5}
\arr{28.5 4.5}{31.5 7.5}
\arr{32.5 0.5}{35.5 3.5}
\arr{28.5 11.5}{31.5 8.5}
\arr{32.5 7.5}{35.5 4.5}

\arr{32.5 8.5}{35.5 11.5}
\arr{36.5 4.5}{39.5 7.5}

\arr{36.5 11.5}{39.5 8.5}

\arr{40.5 8.5}{43.5 11.5}

\arr{4.3 3.5}{7.5 .5}
\endpicture} at 0 0

\put{\beginpicture
\multiput{} at 0 13  0 -7 /
\multiput{$\circ$} at 0 0  4 0  8 0  12 0  16 0  20 0  24 0  28 0  32 0
       4 -4  12 -4  20 -4  28 -4
              4 4   12 4   20 4   28  4
	                8 8    16 8   24 8
			             12 12  20 12 /
				     
\circulararc 360 degrees from 13 11 center at 12 12
\circulararc 360 degrees from 19 11 center at 20 12 
\put{$X$} at 12 15
\put{$Z$} at 20 15

\setdashes <1mm>
\plot 1 -1  13 11 /
\plot -1 1  11 13 /
\circulararc -180 degrees from 1 -1 center at 0 0
\circulararc 180 degrees from 13 11 center at 12 12
\plot 31 -1  19 11 /
\plot 33 1  21 13 /
\circulararc -180 degrees from 19 11 center at 20 12
\circulararc 180 degrees from 31 -1 center at 32 0
\put{$\Cal X$} at -1 -3
\put{$\Cal X'$} at 32.5 -3

\circulararc 360 degrees from -3 1 center at -4 0
\circulararc 360 degrees from 37 -3 center at 36 -4

\setsolid

\arr{0.5 0}{3.5 0}
\arr{4.5 0}{7.5 0}
\arr{8.5 0}{11.5 0}
\arr{12.5 0}{15.5 0}
\arr{16.5 0}{19.5 0}
\arr{20.5 0}{23.5 0}
\arr{24.5 0}{27.5 0}
\arr{28.5 0}{31.5 0}

\arr{0.5 -.5}{3.5 -3.5}
\arr{4.3 -3.5}{7.5 -.5}
\arr{8.5 -.5}{11.5 -3.5}
\arr{12.3 -3.5}{15.5 -.5}
\arr{16.5 -.5}{19.5 -3.5}
\arr{20.3 -3.5}{23.5 -.5}
\arr{24.5 -.5}{27.5 -3.5}
\arr{28.3 -3.5}{31.5 -.5}

\arr{0.5 .5}{3.5 3.5}
\arr{4.3 3.5}{7.5 .5}
\arr{8.5 .5}{11.5 3.5}
\arr{12.3 3.5}{15.5 .5}
\arr{16.5 .5}{19.5 3.5}
\arr{20.3 3.5}{23.5 .5}
\arr{24.5 .5}{27.5 3.5}
\arr{28.3 3.5}{31.5 .5}

\arr{4.5 4.5}{7.5 7.5}
\arr{8.3 7.5}{11.5 4.5}
\arr{12.5 4.5}{15.5 7.5}
\arr{16.3 7.5}{19.5 4.5}
\arr{20.5 4.5}{23.5 7.5}
\arr{24.3 7.5}{27.5 4.5}

\arr{8.5 8.5}{11.5 11.5}
\arr{12.3 11.5}{15.5 8.5}
\arr{16.5 8.5}{19.5 11.5}
\arr{20.3 11.5}{23.5 8.5}

\multiput{$\circ$} at -4 0  0 8  4 12
28 12  32 8  36 4  36 12     36 -4 /

\arr{32.5 -.5}{35.5 -3.5}

\arr{24.5 8.5}{27.5 11.5}
\arr{28.5 4.5}{31.5 7.5}
\arr{32.5 0.5}{35.5 3.5}
\arr{28.5 11.5}{31.5 8.5}
\arr{32.5 7.5}{35.5 4.5}

\arr{32.5 8.5}{35.5 11.5}
\arr{4.3 3.5}{7.5 .5}

\arr{-3.5 0}{-.5 0}
\arr{0.5 7.5}{3.5 4.5}
\arr{4.5 11.5}{7.5 8.5}
\arr{0.5 8.5}{3.5 11.5}
\endpicture} at 45 0
\endpicture}
$$
	\bigskip
The subcategory $\Cal B(n)$ can be characterized as follows. 
    \medskip 

{\bf 2.9. Theorem.} {\it Let $Q$ be a quiver of type $\Bbb D_n$ with $n\ge 4.$ There exists 
a smallest thick subcategory $\Cal T$ of $\rep Q$ which contains all non-thin indecomposable
representations. If $n\ge 5$, we have $\Cal T = \Cal B(n)$, thus $\Cal T$  is 
equivalent to $\rep Q'$ for some quiver $Q'$ of type $\Bbb D_{n-1}$. 
If $n = 4$, we have $\Cal T = \add M$, where $M$ is the maximal indecomposable 
representation, and thus $\Cal T$ is of type $\Bbb A_1$.}
		\medskip
Proof. For $n = 4,$ there is just one non-thin indecomposable representation $M$. Thus $\add M$ is the
smallest thick subcategory which contains all non-thin indecomposable representations. 

Thus, let us assume that $n\ge 5.$ We know that $\Cal B = \Cal B(n)$ is a proper thick subcategory which
contains all non-thin indecomposable representations.
Conversely, let $\Cal T$ be a proper thick subcategory of $\rep Q$ which contains all non-thin 
indecomposable representations.  We have to show that $\Cal T = \Cal B.$ 
Since $\Cal T$ is a proper subcategory, $\Cal T^\perp$ is not zero (see for example [R5]). 
Let $X'$ be an indecomposable representation in $\Cal T^\perp,$ and $Z' = \tau^{-}X'.$
We want to show that $X' = X$ (and $Z' = Z).$ 
Since $X'$ belongs to $\Cal T^\perp$, we know the following:
If $Y$ is an indecomposable representation which is not thin, then 
$\Hom(Y,X') = 0$ and $\Hom(Z',Y) = D\Ext(Y,X') = 0$ (since $Y$ belongs to
$\Cal T$).

Let us look at the successors of $\tau^{-}P(3)$ and at the predecessors of
$\tau I(3)$.
If $N$ is indecomposable and a successor of $\tau^-P(3)$, then one easily sees
that there is a non-thin indecomposable representation $Y$ with $\Hom(Y,N) \neq 0.$
Dually, 
if $N$ is indecomposable and a predecessor of $\tau I(3)$, then 
there is a non-thin indecomposable representation $Y$ with $\Hom(N,Y) \neq 0.$
 
It follows that $X'$ is not a successor of $\tau^-P(3)$, and $Z'$ is
not a predecessor of $\tau I(3)$. 
But the only such indecomposable representation is $X' = X$ (and then $Z' = Z$). 

Let us indicate the position of the representations $\tau^-P(3)$ and $\tau I(3)$
as well as $X$ and $Z$ 
inside the category $\Conic Q(n) = \Cal X\vee \Cal Y \vee \Cal X'$:
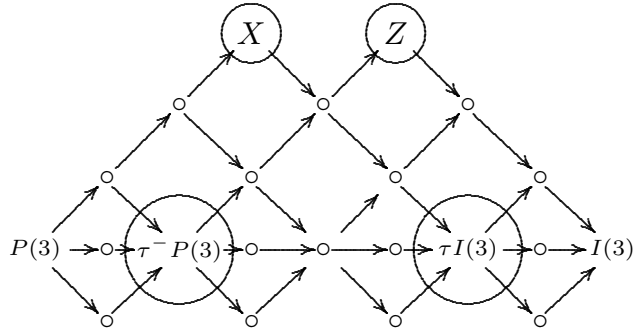
$$
{\beginpicture
\setcoordinatesystem units <.2cm,.2cm>
\multiput{} at 0 13  0 -7 /
\multiput{$\circ$} at  4 0    12 0  16 0  20 0    28 0  
       4 -4  12 -4  20 -4  28 -4
              4 4   12 4   20 4   28  4
	                8 8    16 8   24 8
			              /
\put{$\ssize P(3)$} at 0 0 
\put{$\ssize\tau^-P(3)$} at 8 0 
\put{$\ssize I(3)$} at 32 0 
\put{$\ssize \tau I(3)$} at 24 0 
\setsolid
\circulararc 360 degrees from 13.3 11 center at 12 12
\circulararc 360 degrees from 21.3 11 center at 20 12
\put{$X$} at 12 12
\put{$Z$} at 20 12

\circulararc 360 degrees from 11 0 center at 8 0
\circulararc 360 degrees from 21 0 center at 24 0

\arr{2 0}{3.5 0}
\arr{4.5 0}{5.5 0}
\arr{10.5 0}{11.5 0}
\arr{12.5 0}{15.5 0}
\arr{16.5 0}{19.5 0}
\arr{20.5 0}{22 0}
\arr{26 0}{27.5 0}
\arr{28.5 0}{30.5 0}

\arr{1 -1}{3.5 -3.5}
\arr{4.3 -3.5}{7 -1}
\arr{9 -1}{11.5 -3.5}
\arr{12.3 -3.5}{15 -1}
\arr{17 -1}{19.5 -3.5}
\arr{20.3 -3.5}{23 -1}
\arr{25 -1}{27.5 -3.5}
\arr{28.3 -3.5}{31 -1}

\arr{1 1}{3.5 3.5}
\arr{4.3 3.5}{7 1}
\arr{9 1}{11.5 3.5}
\arr{12.3 3.5}{15 1}
\arr{17 1}{19 3}
\arr{20.3 3.5}{23 1}
\arr{25 1}{27.5 3.5}
\arr{28.3 3.5}{31 1}

\arr{4.5 4.5}{7.5 7.5}
\arr{8.3 7.5}{11.5 4.5}
\arr{12.5 4.5}{15.5 7.5}
\arr{16.3 7.5}{19.5 4.5}
\arr{20.5 4.5}{23.5 7.5}
\arr{24.3 7.5}{27.5 4.5}

\arr{8.5 8.5}{11 11}
\arr{13 11}{15.5 8.5}
\arr{16.5 8.5}{19 11}
\arr{21 11}{23.5 8.5}

\endpicture} 
$$

Altogether, we have shown that $\Cal T^\perp = \add X$ and therefore $\Cal T = {}^\perp X = \Cal B(n).$
This completes the proof. 
     \bigskip 
If $Q$ is a quiver of type $\Bbb D_n$ with $n\ge 5$, then $\Cal B$ itself and 
$\rep Q$ are the only thick subcategories which contain $\Cal B$.

If $Q$ is a quiver of type $\Bbb D_4$, then the thick subcategories of $\rep Q$ which contain the 
non-thin indecomposable representation $M$ form a lattice $L$ of the form
$$
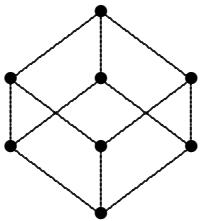
{\beginpicture
\setcoordinatesystem units <1cm,.75cm> 
\multiput{$\bullet$} at 0 0  0 1  0 2  0 3  -1 1  -1 2  1 1  1 2 /
\plot 0 0  -1 1  -1 2  0 3  1 2  1 1  0 0 /
\plot 0 0  0 1 /
\plot 0 2  0 3 /
\plot -1 1  0 2  1 1 /
\plot -1 2  0 1  1 2 /
\endpicture}
$$
Namely, let $Q'$ be obtained from $Q$ by deleting the center, thus $Q'$ is a quiver of
type $\Bbb A_1\sqcup  \Bbb A_1\sqcup  \Bbb A_1$. The lattice of the thick subcategories of $\rep Q$
which contain $M$ is isomorphic to the lattice of thick subcategories of $\rep Q'$ (see for example
[R5]). 

For example, let $Q$ be the $3$-subspace quiver with simple injective representations $S(1),S(2),S(3)$.
As we have mentioned, $P(i)$ denotes the projective cover of $S(i)$, and we let
$V(i) = \tau S(i) = \tau^{-1}P(i)$. Then $L$ looks as follows:
$$
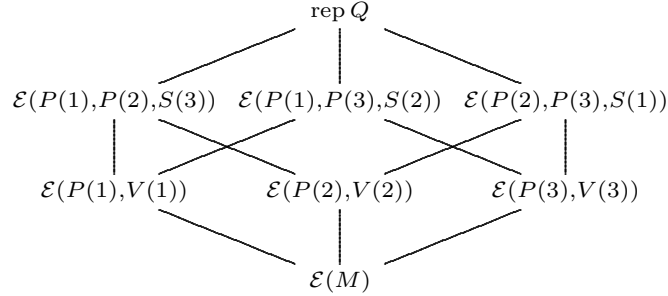
{\beginpicture
\setcoordinatesystem units <2.5cm,1cm> 
\put{$\ssize \ssize\Cal E(M)$} at 0 0  
\put{$\ssize \Cal E(P(1),V(1))$} at -1 1  
\put{$\ssize \Cal E(P(2),V(2))$} at 0 1  
\put{$\ssize \Cal E(P(3),V(3))$} at 1 1  
\put{$\ssize \Cal E(P(1),P(2),S(3))$} at -1 2  
\put{$\ssize \Cal E(P(1),P(3),S(2))$} at 0 2
\put{$\ssize \Cal E(P(2),P(3),S(1))$} at 1 2
\put{$\ssize \rep Q$} at 0 3  
\plot -.2 .2  -.8 .8 /
\plot -1 1.2  -1 1.8  /
\plot -.8 2.2  -.2 2.8 /
\plot  0.2 2.8   0.8 2.2 / 
\plot 1 1.8   1 1.2  /
\plot 0.8 0.8  0.2 0.2 /
\plot 0 0.2  0 .8 /
\plot 0 2.2  0 2.8 /
\plot -0.8 1.2  -.2 1.8  /
\plot 0.2 1.8   .8 1.2 /
\plot -.8 1.8  -.2 1.2  /
\plot 0.2 1.2  0.8 1.8 /
\endpicture}
$$
	\bigskip
We may proceed inductively: Let $Q$ be a quiver of type $\Bbb D_n$. For $0 \le t \le n-3$, 
let $\Cal B^{(t)}$ be the smallest thick subcategory which contains all $s$-twin representations
with $s\ge t$. Then
$$
 \rep Q = \Cal B^{(0)} \supset \Cal B = \Cal B^{(1)} \supset \cdots \supset \Cal B^{(n-3)},
$$
and $\Cal B^{(t)}$ is of type $\Bbb D_{n-t}$ for $0\le t \le n-4$ and of type $\Bbb A_1$ for $t=n-3.$ 
    \bigskip\bigskip 
{\bf 2.10. Final remark for Part 2.} The aim of Part 2 was to collect the relevant indecomposable
representations as a subcategory $\Cal B$ which contains all the non-thin indecomposable
representations and only few thin representations. 

For a long time there was the strong belief 
that one of the first themes of representation theory should be to provide
lists of all the indecomposable representations, whenever this is possible, and, 
only then, as a second step,
to determine homomorphisms and extensions between the indecomposable objects. Of course, such a listing 
will be possible in case we deal with a representation-finite artin algebras, but it turns out that
usually the list of indecomposables is quite uninteresting: it is the internal categorical structure and the
interplay between indecomposable representations which should be described. In order to do so, one may
look at sets of indecomposables which are related either by small changes of parameters or by the existence of
irreducible maps. 

This procedure is well accepted
in case one deals with tame artin algebras, where the one-parameter families always are considered as units.
There are three classes of artin algebras with good descriptions of all the indecomposable modules:
the tame concealed algebras, the tubular algebras and the special biserial algebras. 
For the special biserial algebras, the string modules 
may be considered individually, but the band modules for a primitive cyclic word always are considered as a unit. 
For the tame concealed and the tubular algebras, one looks at the preprojective and the preinjective 
components, as well as the tubular families. 

Our discussion of the quivers of type $\Bbb D$ uses a similar principle. The full subcategory of 
the conical representations is a unit which is independent of the orientation and which draws the attention to
the representations which are not thin. As we have seen, we can decompose this unit into smaller blocks
$\Cal X, \Cal Y, \Cal X'$ whose internal structure is important, as well. Clearly, 
these blocks are the data of interest, whereas the individual indecomposables gain their relevance by
their position inside these blocks. 
      \bigskip\bigskip 
\centerline{\bf Part 3. The quivers of type $\Bbb E_6,\ \Bbb E_7,\ \Bbb E_8$ (and again $\Bbb D_n$).}
		       \bigskip
Recall that any Dynkin quiver has exceptional vertices, those of type $\Bbb D_n$ and 
$\Bbb E_m$ have
just one exceptional vertex. 
If $\Delta$ is a Dynkin quiver of type $\Bbb D_n$ or 
$\Bbb E_m$, and $y$ is its exceptional vertex, we denote by $\Delta'$ the quiver
obtained from $\Delta$ by deleting $y$, and by $\Delta''$ the quiver obtained from
$\Delta$ be deleting $y$ and the neighbors of $y$.

A main result of Part 3 is the following theorem.  
     \medskip 
{\bf 3.1. Theorem.} {\it 
Let $\Delta$ be a Dynkin quiver of type $\Bbb D_n$ or $\Bbb E_m.$  
There is a unique antichain triple $A(1), A(2), A(3)$ in $\rep \Delta$ such that
the support of any $A(i)$ is in $\Delta'$, but not in $\Delta''$.

In addition, we have: $\Ext^1(A(i),A(j)) = 0$ for all $i,j$. For any $i$, 
there is only one neighbor
$x$ of the exceptional vertex
with $A(i)_x \neq 0$, and $\dim A(i)_x = 1$ for this vertex $x$.

Finally, if $M$ is the maximal indecomposable representation of $\Delta$, then 
$$
 M|\Delta' = A(1)\oplus A(2)\oplus A(3)
$$ 
(and $\dim M_y = 2$, for the exceptional vertex $y$).}
    \medskip
The triple $A(1),A(2),A(3)$ will be called the {\it special antichain triple} in $\rep \Delta$. 
The proof will be based on a careful study of the quiver $\Delta'$. It will
be completed in section 3.6.
The further considerations in part 3 
will outline some consequences. In addition, we will draw the 
attention also to $M|\Delta''.$

	\bigskip
We add two remarks, using the following definition.
Given a dimension vector $\bold d$, let $\n_a(\bold d)$ be the sum
of the numbers $\bold d_b$ with $b$ a neighbor of $a$.
	\medskip
{\bf Remark 1.} If $\Delta$ is a Dynkin quiver of type $\Bbb D_n$ or $\Bbb E_m$ and
$M$ is its maximal indecomposable representation, then 
it is easy to see that $M|\Delta'$ is the direct sum of 
an antichain triple, since $\n_y(\bdim M) = 3.$

Namely, there is the following general fact: {\it Let $\Delta$ be a Dynkin quiver and 
$a\in\Delta_0$. If $N$ is an indecomposable
representation of $\Delta$, with $N_a \neq 0$, and 
$\n_a(\bdim N) \le 3$, then $N|\Delta\setminus\{a\}$ is the direct sum of
an antichain of cardinality $\n_a(\bdim N)$.}

Proof. If $X$ is an indecomposable direct summand of $N|\Delta\setminus\{a\}$,
then $\n_a(\bdim X) = 1$, since $\n_a(\bdim X) \ge 2$ would imply that  
$\Cal E(X,S(a))$ is representation-infinite.
This shows that $N|\Delta\setminus\{a\}$ is the direct sum of at most 3
indecomposable representations. In order to see that these representations form
an antichain, one may, for example, refer to Kleiner's list of posets of finite type
[R3].
	\medskip
{\bf Remark 2.} {\it If $M$ is the maximal indecomposable representation of the
Dynkin quiver $\Delta$ and $a\in\Delta_0$, then $\n_a(\bdim M) = 3$
if and only if $\Delta$ is of type $\Bbb D_n$ with $n\ge 4$ or of type $\Bbb E_m$
and $a$ is the exceptional vertex.}

For the proof, we just have to calculate the numbers $\n_a(\bdim M)$. The following table
shows on the left the dimension vector $\bdim M$; on the right, we display in the
same way the numbers $\n_a(\bdim M)$.
\vfill

$$
{\beginpicture
\setcoordinatesystem units <0.6cm,0.6cm>
\put{$\Delta$} at -4 1.5 
\put{$\bdim M$} at 0 1.5
\put{$\n(\bdim M)$} at 8 1.5 

\put{\beginpicture
\multiput{$1$} at  1 0  2 0  4 0 /
\multiput{$1$} at 0 0  5 0 /
\plot 0.2 0  0.8 0 /
\plot 1.2 0  1.8 0 /
\plot 2.2 0  2.5 0 /
\plot 3.5 0  3.8 0 /
\plot 4.2 0  4.8 0 /
\put{$\cdots$} at 3 0 
\put{$\Bbb A_n$} at -2 0
\put{} at 6 0
\endpicture} at 0 -.1

\put{\beginpicture
\multiput{$2$} at  1 0  2 0  4 0 /
\multiput{$1$} at 0 0  5 0 /
\plot 0.2 0  0.8 0 /
\plot 1.2 0  1.8 0 /
\plot 2.2 0  2.5 0 /
\plot 3.5 0  3.8 0 /
\plot 4.2 0  4.8 0 /
\put{$\cdots$} at 3 0 
\put{} at 6 0
\endpicture} at 9 -.1

\put{\beginpicture
\setcoordinatesystem units <0.6cm,0.4cm>
\multiput{$1$} at 0 1  0 -1   5 0 /
\multiput{$2$} at 1 0  3 0  4 0 /
\plot 0.2 0.8  0.8 0.2 /
\plot 0.2 -.8  0.8 -.2 /
\plot 1.2 0  1.5 0 /
\plot 2.5 0  2.8 0 /
\plot 3.2 0  3.8 0 /
\plot 4.2 0  4.8 0 /
\put{} at 6 0
\put{$\cdots$} at 2 0 
\put{$\Bbb D_n$} at -2 0

\endpicture} at 0 -2
\put{\beginpicture
\setcoordinatesystem units <0.6cm,0.4cm>
\multiput{$2$} at 0 1  0 -1   5 0 /
\multiput{$4$} at 1 0  3 0  /
\put{$3$} at   4 0 
\plot 0.2 0.8  0.8 0.2 /
\plot 0.2 -.8  0.8 -.2 /
\plot 1.2 0  1.5 0 /
\plot 2.5 0  2.8 0 /
\plot 3.2 0  3.8 0 /
\plot 4.2 0  4.8 0 /
\put{} at 6 0
\put{$\cdots$} at 2 0

\endpicture} at 9 -2
\put{\beginpicture
\setcoordinatesystem units <0.6cm,0.5cm>
\put{$1$} at 0 0
\put{$2$} at 1 0
\put{$3$} at 2 0
\put{$2$} at 3 0
\put{$1$} at 4 0
\put{$2$} at 2 1 
\plot 0.2 0  0.8 0 /
\plot 1.2 0  1.8 0 /
\plot 2.2 0  2.8 0 /
\plot 3.2 0  3.8 0 /
\plot 2 0.3  2 0.7 /
\put{$\Bbb E_6$} at -2 0
\put{} at 6 0
\endpicture} at 0 -4
\put{\beginpicture
\setcoordinatesystem units <0.6cm,0.5cm>
\put{$2$} at 0 0
\put{$4$} at 1 0
\put{$6$} at 2 0
\put{$4$} at 3 0
\put{$2$} at 4 0
\put{$3$} at 2 1 
\plot 0.2 0  0.8 0 /
\plot 1.2 0  1.8 0 /
\plot 2.2 0  2.8 0 /
\plot 3.2 0  3.8 0 /
\plot 2 0.3  2 0.7 /
\put{} at 6 0
\endpicture} at 9 -4
\put{\beginpicture
\setcoordinatesystem units <0.6cm,0.5cm>
\put{$2$} at 0 0
\put{$3$} at 1 0
\put{$4$} at 2 0
\put{$3$} at 3 0
\put{$2$} at 4 0
\put{$1$} at 5 0 
\put{$2$} at 2 1 
\plot 0.2 0  0.8 0 /
\plot 1.2 0  1.8 0 /
\plot 2.2 0  2.8 0 /
\plot 3.2 0  3.8 0 /
\plot 4.2 0  4.8 0 /
\plot 2 0.3  2 0.7 /
\put{$\Bbb E_7$} at -2 0
\put{} at 6 0
\endpicture} at 0 -6
\put{\beginpicture
\setcoordinatesystem units <0.6cm,0.5cm>
\put{$3$} at 0 0
\put{$6$} at 1 0
\put{$8$} at 2 0
\put{$6$} at 3 0
\put{$4$} at 4 0
\put{$2$} at 5 0 
\put{$4$} at 2 1 
\plot 0.2 0  0.8 0 /
\plot 1.2 0  1.8 0 /
\plot 2.2 0  2.8 0 /
\plot 3.2 0  3.8 0 /
\plot 4.2 0  4.8 0 /
\plot 2 0.3  2 0.7 /
\put{} at 6 0
\endpicture} at 9 -6
\put{\beginpicture
\setcoordinatesystem units <0.6cm,0.5cm>
\put{$2$} at 0 0
\put{$4$} at 1 0
\put{$6$} at 2 0
\put{$5$} at 3 0
\put{$4$} at 4 0
\put{$3$} at 5 0 
\put{$2$} at 6 0
\put{$3$} at 2 1 
\plot 0.2 0  0.8 0 /
\plot 1.2 0  1.8 0 /
\plot 2.2 0  2.8 0 /
\plot 3.2 0  3.8 0 /
\plot 4.2 0  4.8 0 /
\plot 5.2 0  5.8 0 /
\plot 2 0.3  2 0.7 /
\put{$\Bbb E_8$} at -2 0
\endpicture} at 0 -8
\put{\beginpicture
\setcoordinatesystem units <0.6cm,0.5cm>
\put{$4$} at 0 0
\put{$8$} at 1 0
\put{$12$} at 2 0
\put{$10$} at 3 0
\put{$8$} at 4 0
\put{$6$} at 5 0 
\put{$3$} at 6 0
\put{$6$} at 2 1 
\plot 0.2 0  0.8 0 /
\plot 1.2 0  1.7 0 /
\plot 2.3 0  2.7 0 /
\plot 3.3 0  3.8 0 /
\plot 4.2 0  4.8 0 /
\plot 5.2 0  5.8 0 /
\plot 2 0.3  2 0.7 /
\endpicture} at 9 -8
\put{} at 0 -8.9
\put{} at 0 -8.9
\endpicture}
$$

We see that the value 3 occurs in the right column only in 
the cases $\Bbb D_n$ and $\Bbb E_m$, and only at the exceptional vertex. 
Actually, for $\Bbb D_n$ and $\Bbb E_m$ the value
3 is the only odd value which occurs (and it is the only value $\n_a(\bdim M)$ which
is different from $2\cdot \bdim M_a$).

In view of Remark 1, we can say: {\it  If $\Delta$ is of type $\Bbb D_n$ or $\Bbb E_m$,
there is a unique vertex $a$ of $\Delta$
such that $\n(M|\Delta\setminus\{a\}) = 3$, namely the exceptional vertex of $\Delta.$}
    	      \bigskip
{\bf 3.2. Hammocks.}
Hammocks have been introduced by Brenner [B] (for a general theory see [HV], the special case
of dealing with quiver representations was already considered by Gabriel [G2]). If $Q$ is a 
Dynkin quiver, and $x$ a vertex of $Q$, the hammock $H(Q,x)$ is a subset of $\Bbb ZQ$ 
which is obtained using the knitting algorithm in the same way as one constructs a preprojective
component (however, not
using representations or dimension vectors, but just single numbers; actually, each number 
can be interpreted as the dimension of the vector space $M_x$, where
$M$ is an indecomposable  preprojective representation). Here is the definition:
    \medskip
{\bf Knitting algorithm for the hammock function $h_{(Q,x)}$,} where $Q$ is a finite connected acyclic
quiver and $x\in Q_0$.  The hammock function $h = h_{(Q,x)}$ is a function $(\Bbb ZQ)_0 \to \Bbb N_0$
defined as follows: Start: Let $h(a) = 0$ for any vertex $a$ in $\Bbb ZQ = Q_0\times \Bbb Z$ 
with a proper path from $a$ to $(x,0);$ and let   
$h(x,0) = 1$. Inductive procedure: Assume that $c$ is a vertex of $(\Bbb ZQ)_0$ such that
$h( \tau c)$ as well as $h(b)$ for any arrow $b\to c$ in $\Bbb ZQ$ are defined. Then either
$-h(\tau(c)) + \sum_{\alpha\:t(\alpha) \to c} h(s(\alpha))$ is non-negative, then this should be
the value $h(c)$. Otherwise put $h(c) = 0.$ By definition, the {\it hammock} $H(Q,x)$ is
the support of the hammock function $h_{(Q,x)}$, this is a translation subquiver of $\Bbb Z Q$.
Note that $H(Q,x)$ depends only on the underlying graph $\overline Q$ of the quiver $Q$ (and not on
the orientation of the edges).
    \medskip
Let us exhibit two hammock functions for a quiver of type $\Bbb D_6$:
$$  
 
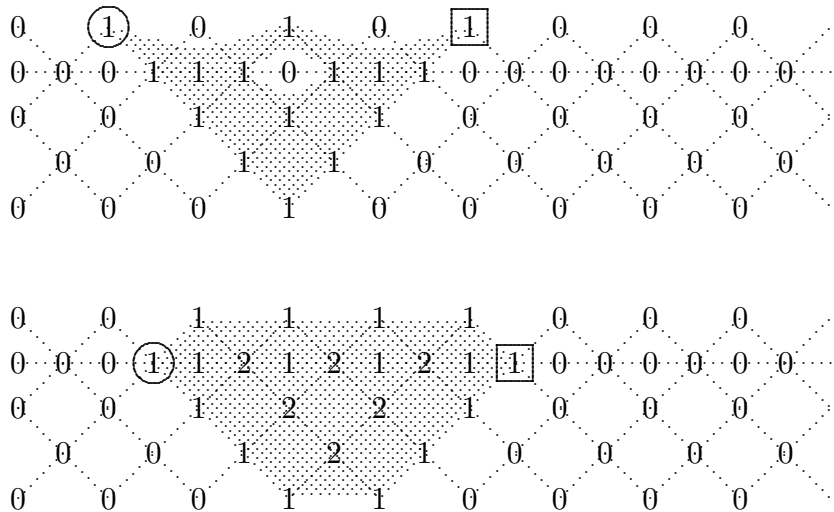
{\beginpicture
 \setcoordinatesystem units <.5cm,.5cm>
\circulararc 360 degrees from 2.45 4 center at 2 4 
\plot 10.4 3.6  10.4  4.4  9.6 4.4  9.6 3.6  10.4 3.6 /
\multiput{$1$} at 2 4  3 3  4 2  5 1  6 0   4 3  5 3  6 4  6 2  7 3  7 1
    8 2  8 3  9 3   10 4  /
\multiput{$0$} at  0 0  0 2  0 4  1 1  1 3 
   2 0  2 2   3 1   
    4 0   4 4  
  8 0   8 4  9 1  
 10 0  10 2   11 1  11 3 
 12 0  12 2  12 4  13 1  13 3 
 14 0  14 2  14 4  15 1  15 3 
 16 0  16 2  16 4  17 1  17 3 /
\multiput{$0$} at 0 3  2 3   6 3   10 3  12 3  14 3  16 3 /
\setdots <1mm>
\plot 0 0   4 4  8 0  12 4  16 0  18 2 /
\plot 0 4  4 0  8 4  12 0  16 4  18 2 /
\plot 0 2  2 4   6 0  10 4  14 0  18 4 /
\plot 0 2  2 0  6 4  10 0  14 4  18 0 /
\plot 0 3  18 3 /
\setshadegrid span <.5mm>
\vshade  2 4 4  <,z,,> 4 2 3.5 <z,z,,>  5 1 3.75  /
\vshade  5 3 3.75 <z,z,,> 6 3.5 4.1 <z,z,,>   7 3 3.75 /
\vshade 5 1 3 <z,z,,>   6 0 2.5 <z,z,,>   7 1 3 /
\vshade  7 1 3.75 <z,z,,> 8 2 3.5  <z,,,>  10 4 4 /
\endpicture}
$$
	\medskip 
$$  
 {\beginpicture
 \setcoordinatesystem units <.5cm,.5cm>
\circulararc 360 degrees from 3.45 3 center at 3 3 
\multiput{$2$} at 5 3  6 2  7 1  7 3  8 2  9 3 /
\multiput{$1$} at 3 3  4 2  4 4   5 1  6 0   4 3    6 3   6 4  
    8 3   8 4   10 3  10 4   8 0  9 1  10 2  11 3 /
\multiput{$0$} at  0 0  0 2  0 4  1 1  1 3 
   2 0  2 2   2 4  3 1   
    4 0     
 10 0   11 1 
 12 0  12 2   12 4  13 1  13 3 
 14 0  14 2  14 4  15 1  15 3 
 16 0  16 2  16 4  17 1  17 3 /
\plot 11.4 2.6  11.4  3.4  10.6 3.4  10.6 2.6  11.4 2.6 /
\multiput{$0$} at 0 3  2 3  12 3  14 3  16 3 /
\setdots <1mm>
\plot 0 0   4 4  8 0  12 4  16 0  18 2 /
\plot 0 4  4 0  8 4  12 0  16 4  18 2 /
\plot 0 2  2 4   6 0  10 4  14 0  18 4 /
\plot 0 2  2 0  6 4  10 0  14 4  18 0 /
\plot 0 3  18 3 /
\setshadegrid span <.5mm>
\vshade  3 3 3 <,z,,> 4 2 4 <z,z,,> 6 0 4 <z,z,,> 8 0 4 <z,z,,> 10 2 4 <z,,,> 11 3 3 /
\put{} at 0 -.5
\endpicture}
$$
In both examples, the corresponding hammock is shaded; the vertex $(x,0)$ 
(this is the source of the hammock) is encircled, the last vertex $c$ obtained with non-zero 
value $h(c)$ is marked by a square (this is the sink of the hammock).
      \medskip
In case $h_{(Q,x)}(y) \le 1$ for all vertices $y$, we define on the set $H(Q,x)$ a relation $\le$ as
follows: $y\le z$ provided there is a path $y = y_0 \to y_1 \to \cdots y_t = z$ in $\Bbb ZQ$ such that all
vertices $y_i$ belong to $H(Q,x)$. It is well-known (and easy to see) that in this way $H(Q,x)$
becomes a poset. By abuse of language, we will say that $H(Q,x)$ {\it is a poset} provided 
$h_{(Q,x)}(y) \le 1$ for all vertices $y$. 
	      \medskip 
If we look at the two hammock functions for a quiver of type $\Bbb D_6$ exhibited above,
we see that the first hammock satisfies the poset condition $h_{(Q,x)}(a) \le 1$ for all
vertices of $\Bbb ZQ$; the shape of this poset will be shown below. The second hammock 
function does not satisfy the poset condition.

	 \medskip 
{\bf Hammock sets.} Let $Q(1),\dots, Q(t)$ be Dynkin quivers and let $Q$ be the disjoint union
of the quivers $Q(i)$. Let $x(i)$ be a vertex of $Q(i)$, for $1\le i \le t.$ The disjoint union 
$$
 H(Q,x(1),\dots,x(t)) = \bigsqcup\nolimits_i H(Q(i),x(i))
$$
of the hammocks $H(Q(i),x(i))$ will be called a {\it hammock set.} In case all the hammocks
$H(Q(i),x(i))$ are posets, we consider also $H(Q,x(1),\dots,x(t))$ as a poset and we will say
that $H(Q,x(1),\dots,x(t))$ {\it is a poset.}
     \bigskip
{\bf Hammock categories.} 
Given a representation $M$ of a quiver $Q$, we denote by $\Cal H(M)$ the following category:
its objects are the indecomposable representations $X$ of $Q$ and $\Hom_{\Cal H(M)}(X,X')$
is the set of equivalence classes of maps $f\:X \to X'$ where maps $f,f'\:X \to X'$ are equivalent
provided $\Hom(M,f-f') = 0.$ Note that the isomorphism classes of indecomposable objects in
$\Cal H(M)$ are just the isomorphism classes of indecomposable representations $X$ of $Q$
such that $\Hom(M,X) \neq 0.$ If $M = P(x)$, then $\Cal H(M)$ is the image of the forgetful functor
from $\rep Q$ to the category of $k$-spaces which send a representation $X$ to $X_x$ and a map
$f$ to $f_x$. 

Let us consider the case $M = P(x)$, where $x$ is a vertex of a Dynkin quiver $Q$. 
Since $\dim \Hom(P(x),M) = \dim M_x$, for any representation $M$ of $Q$,
the hammock function $h = h_{(Q,x)}$ describes $\Cal H(P(x)):$ 
{\it The  hammock $H(Q,x)$ is the quiver of the hammock category $\Cal H(P(x))$;}
the vertices $a$ of the hammock $H(Q,x)$ may be considered as the
isomorphism classes $[M]$ 
of the indecomposable representations $N$ of $Q$ with 
$N_x \neq 0$ {\it and $h_{(Q,x))}(a)$ is just the dimension of the vector space $N_x$.}

     \medskip 
Proof. The knitting algorithm for the hammock function is nothing else than the knitting
algorithm for one of the coordinates of the dimension vector $\bdim N$ of the 
indecomposable representations $N$, namely for $(\bdim N)_x = \dim N_x.$
	       \medskip 
Always, the hammock $H(Q,x)$ has a unique source $a$, namely the isomorphism class of the
representation $P(x)$. 
If $Q$ is a Dynkin quiver, then $H(Q,x)$ also has a (unique) sink $z$, namely the isomorphism class
of the indecomposable injective representation $M(z) = I(x)$.
In the two $\Bbb D_6$-examples
above, the encircled vertex $a$ is the position of $P(x)$, the vertex $z$ marked by a square
is the position of $I(x)$.
   \medskip
{\it If $Q$ is a star quiver with center $c$ and representation-finite, then the hammock category
$\Cal H(Q,c)$ is just the full subcategory of all conical representations.} The second example
above is of this kind.   
      \medskip 
More generally, assume that we deal with a hammock set
$$
 H(Q,x(1),\dots,x(t)) = \bigsqcup\nolimits_i H(Q(i),x(i)),
$$
thus $Q(1),\dots, Q(t)$ are the connected components of a quiver $Q$ and 
$x(i)$ is a vertex of the Dynkin quiver $Q(i)$, for $1\le i \le t.$ Then we have: 
{\it The  hammock set 
$H(Q,x(1),\dots,x(t))$ is the quiver of the hammock category $\Cal H(\bigoplus_i P(x(i)))$.} 
		       \bigskip 
{\bf 3.3. Quivers with a special vertex set.}
Let $Q(1),\dots, Q(t)$ be Dynkin quiver and let $Q$ be the disjoint union
of the quivers $Q(i)$. Let $x(i)$ be a vertex of $Q(i)$, for $1\le i \le t.$ We say that
$x(1),\dots,x(t)$ is a {\it special vertex set of $Q$} provided 
$H(Q,x(1),\dots,x(t))$ is a poset with precisely one antichain triple.
		       \medskip 

{\bf Proposition.} {\it 
Let $Q(1),\dots, Q(t)$ be Dynkin quivers and let $Q$ be the disjoint union
of the quivers $Q(i)$. Let $x(i)$ be a vertex of $Q(i)$, for $1\le i \le t.$ Then
$x(1),\dots,x(t)$ is s special vertex set if and only if one of the following 
cases holds:

The first possibility $\Bbb A_1\sqcup \Bbb A_1\sqcup \Bbb A_1$ is 
$$
  Q = Q(1)\sqcup Q(2)\sqcup Q(3) \quad\text{with}\quad Q(i) = \{x(i)\} \quad\text{for\ \ } 1\le i \le 3,
$$ 
thus all the connected components $Q(i)$ of $Q$ are of type $\Bbb A_1.$ 

The second possibility $\Bbb A_1\sqcup \Bbb D_m$ is 
$$
{\beginpicture
\setcoordinatesystem units <2cm,.6cm> 
\put{$Q = Q(1)\sqcup Q(2)$\strut} at 0 0 
\put{with\strut} at 0.9 0
\put{$Q(1) = \{x(1)\}$\strut} at 1.7 0 
\put{and \ $Q(2)=$\strut} at 2.8 0
\put{\beginpicture
\setcoordinatesystem units <.6cm,.6cm> 
\multiput{$\circ$} at 0 1  0 -1  1 0  2 0  4 0 /
\put{$\cdots$} at 3 0
\put{$x(2)$} at 5.3 0
\plot 0.3 0.7 0.7 0.3 /
\plot 0.3 -.7 0.7 -.3 /
\plot 1.3 0  1.7 0 /
\plot 2.3 0  2.5 0 /
\plot 3.5 0  3.7 0 /
\plot 4.3 0  4.7 0 /
\endpicture} at 4.3 0
\endpicture}
$$
here, the component $Q(1)$ of $Q$ is of type $\Bbb A_1$, the component $Q(2)$ is of type
$\Bbb D_m$ with $m\ge 3$} (recall that by definition $\Bbb D_3 = \Bbb A_3$
and $x(2)$ is the central vertex).

{\it And finally,  there are the possibilities that $Q$ is  of type $\Bbb A_5, \Bbb D_6, \Bbb E_7$
and $x = x(1)$ is the following encircled vertex:}
$$
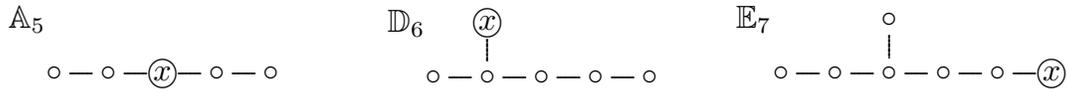
{\beginpicture
\setcoordinatesystem units <.6cm,.6cm> 
\put{\beginpicture
\multiput{$\circ$} at 0 0  1 0  3 0  4 0 /
\put{$x$} at 2 0
\put{$\bigcirc$} at 2 0
\plot 0.3 0  0.7 0 /
\plot 1.3 0  1.7 0 / 
\plot 2.3 0  2.7 0 /
\plot 3.3 0  3.7 0 /
\put{} at 0 1.05
\put{$\Bbb A_5$} at -.5 1
\endpicture} at 0 0 
\put{\beginpicture
\multiput{$\circ$} at 0 0  1 0  2 0  3 0  4 0 /
\put{$x$} at 1 1
\put{$\bigcirc$} at 1 1
\plot 0.3 0  0.7 0 /
\plot 1.3 0  1.7 0 /
\plot 2.3 0  2.7 0 /
\plot 3.3 0  3.7 0 /
\plot 1 0.3  1 0.7 /
\put{$\Bbb D_6$} at -.5 1
\endpicture} at 7 0 
\put{\beginpicture
\multiput{$\circ$} at -1 0  0 0  1 0  1 1  2 0  3 0  /
\put{$x$} at 4 0
\put{$\bigcirc$} at 4 0
\plot -.3 0  -.7 0 /
\plot 0.3 0  0.7 0 /
\plot 1.3 0  1.7 0 /
\plot 2.3 0  2.7 0 /
\plot 3.3 0  3.7 0 /
\plot 1 0.3  1 0.7 /
\put{$\Bbb E_7$} at -1.5 1
\endpicture} at 14 0 
\endpicture}
$$
Note that {\it the vertices $x(1),\dots,x(t)$ are determined by $Q$ up to symmetry} 
(and uniquely determined except in case $\Bbb D_6$). 
     \medskip 
Sketch of the proof. We assume that 
$H = H(Q,x(1),\dots,x(t))$ is a poset with precisely one antichain triple.
Since $H(Q,x(1),\dots,x(t))$  is the disjoint union of the hammocks $H(Q(i),x(i)),$ 
and $H$ has width 3, we see that $t\le 3.$ 

First assume that $t=3$. We want to show that any $Q(i)$ is of type $\Bbb A_1.$ If say
$Q(t)$ is not of type $\Bbb A_1$, then $H(Q(t),x(t))$
has at least two elements. This shows that there are
at least 2 antichain triples  (using one element from $H(Q(1),x(1))$ and one from $H(Q(2),x(2))$).

Second, let $t=2$, then we may assume that $H(Q(1),x(1))$ is a chain, whereas $H(Q(2),x(2))$ has
an antichain pair. As in the case $t=3$, we see that $H(Q(1),x(1))$ has to consist of a single
vertex, thus $Q(1)$ has to be of type $\Bbb A_1.$ In addition, $H(Q(2),x(2))$ can have only
one antichain pair. Since a hammock has a unique source and a unique sink, we see that
$H(Q(2),x(2))$ must have the form 
$$
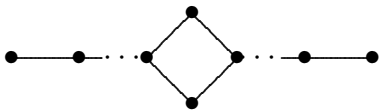
{\beginpicture
\setcoordinatesystem units <.5cm,.5cm> 
\multiput{$\bullet$} at 2 -2  2 0  /
\multiput{$\bullet$} at 1 -1   -.5 -1 -2 -1   3 -1  4.5 -1  6 -1  /

\multiput{$\cdots$} at 0.5 -1  3.5 -1 /
\plot -2 -1  0 -1 /
\plot 1 -1  2 0  3 -1 /
\plot 4 -1  6 -1 /

\plot 1 -1  2 -2  3 -1 /
\endpicture}
$$
But then there are only the following possibilities: $Q(2)$ is of type $\Bbb A_3$ and $x$ is 
the central vertex, 
or  $Q(2)$ is of type $\Bbb D_4$ and $x$ is a leaf vertex, or
$Q(2)$ is of type $\Bbb D_n$ with $n\ge 5$ and $x$ is the leaf vertex of the long arm. In the
formulation of Theorem 3.3, these cases have been collected under the label
$\Bbb D_m$ with $m\ge 3.$ 

It remains to consider the case that $Q$ is connected and $x = x(1)$ is a vertex such that the
hammock $H(Q,x)$ is a poset with precisely one antichain triple.
There are only few cases such that a hammock $H(Q,x)$ is 
a poset, namely the following cases.

First of all, $Q$ can be of type $\Bbb A$ and $x$ arbitrary.
If $x$ is a leaf vertex, then $H(Q,x)$ has width 1, if $x$ is the neighbor of a
leaf vertex (and not a leaf vertex), then $H(Q,x)$ has width 2. Otherwise,
the width of $H(Q,x)$ is greater or equal to $3$, as we want, however usually
there are several antichain triples, the only exception is the
$\Bbb A_5$ case.
      
The remaining possibilities for $x$ are the leaf vertices of $\Bbb D_n$, 
two leaf vertices for $\Bbb E_6$ and one for $\Bbb E_7$.
$$
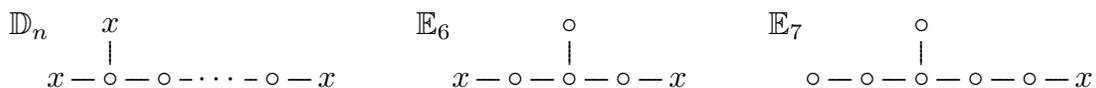
{\beginpicture
\setcoordinatesystem units <.6cm,.6cm> 
\put{\beginpicture
\multiput{$\circ$} at   1 0  2 0   4 0  /
\multiput{$x$} at 0 0  1 1  5 0 /
\put{$\cdots$} at 3 0 
\plot 0.3 0  0.7 0 /
\plot 1.3 0  1.7 0 /
\plot 2.3 0  2.5 0 /
\plot 3.5 0  3.7 0 /
\plot 4.3 0  4.7 0 /
\plot 1 0.3  1 0.7 /
\put{$\Bbb D_n$} at -.5 1
\endpicture} at 0 0 
\put{\beginpicture
\multiput{$\circ$} at   0 0  1 0  1 1  2 0  /
\multiput{$x$} at  -1 0  3 0 /
\plot -.3 0  -.7 0 /
\plot 0.3 0  0.7 0 /
\plot 1.3 0  1.7 0 /
\plot 2.3 0  2.7 0 /
\plot 1 0.3  1 0.7 /
\put{$\Bbb E_6$} at -1.5 1
\endpicture} at 7 0 
\put{\beginpicture
\multiput{$\circ$} at -1 0  0 0  1 0  1 1  2 0  3 0  /
\put{$x$} at 4 0
\plot -.3 0  -.7 0 /
\plot 0.3 0  0.7 0 /
\plot 1.3 0  1.7 0 /
\plot 2.3 0  2.7 0 /
\plot 3.3 0  3.7 0 /
\plot 1 0.3  1 0.7 /
\put{$\Bbb E_7$} at -1.5 1
\endpicture} at 14 0 
\endpicture}
$$
For $\Bbb D_4$ and $x$ any leaf vertex, 
for $\Bbb D_5$ and $x$ the leaf vertex of a short arm, 
as well as for $\Bbb D_n$ with $n\ge 5$ and $x$ the leaf vertex of the long arm, 
the hammock $H(Q,x)$
has width 2. For $n > 6$ and $x$ the leaf vertex on a short arm, 
the hammock $H(Q,x)$ contains
several antichain triples. Finally, for $\Bbb E_6$, and $x$ 
the leaf vertex of a long arm, the hammock
$H(Q,x)$ has width 2. This shows that only the cases mentioned in the theorem
remain. 
	\bigskip 
Here are the posets $H(Q,x_1,\dots,x_t)$ for the various special vertex sets $x(1),\dots,x(t)$ of 
a quiver $Q$. We draw the 
Hasse diagram of each poset, with increase from left to right; we also may interpret the pictures as
describing the quiver of $\Cal H(\bigoplus P(x(i)))$, but deleting the
arrow heads (all arrows point to the right, thus are
northeast, east or southeast arrows).
$$
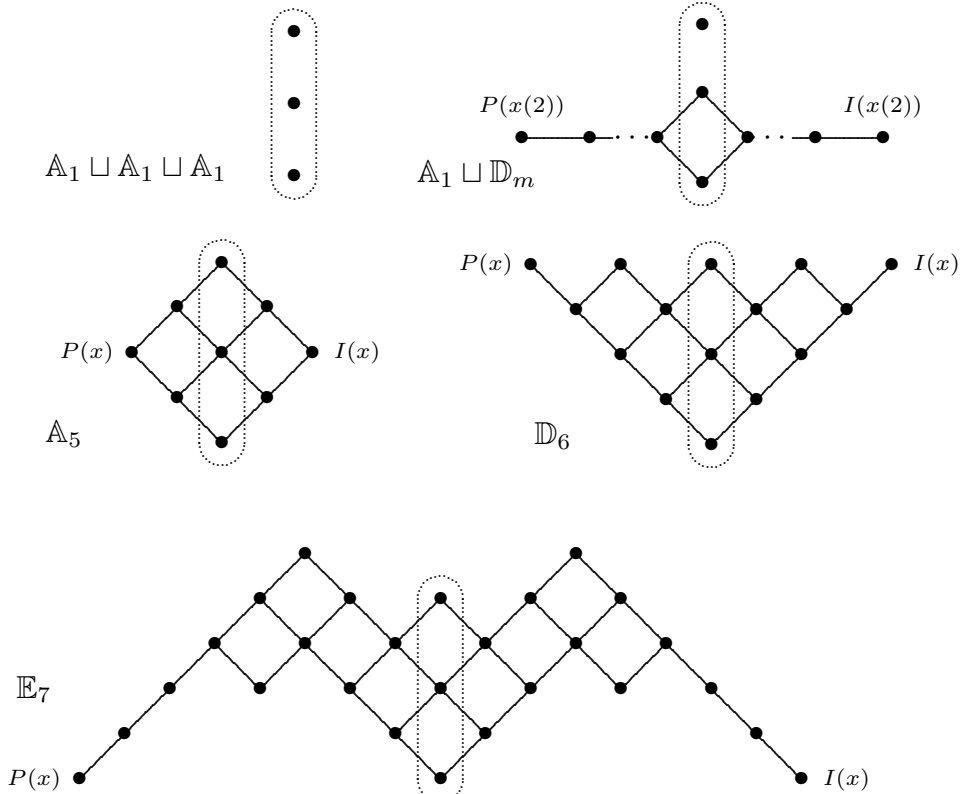
{\beginpicture
\setcoordinatesystem units <.5cm,.5cm> 
\put{\beginpicture
\setcoordinatesystem units <.5cm,.5cm> 
\multiput{$\bullet$} at 2 -2 2 -.4  2 1.2   /
\setdots <.5mm>
\plot 1.5 -2  1.5 1.2 /
\plot 2.5 -2  2.5 1.2 /

\circulararc 180 degrees from 1.5 -2  center at 2 -2 
\circulararc 180 degrees from 2.5 1.2  center at 2 1.2 
\put{$\Bbb A_1\sqcup\Bbb A_1\sqcup\Bbb A_1$} at -1.5 -1.8 
\endpicture} at 0 0 
\put{\beginpicture
\setcoordinatesystem units <.5cm,.5cm> 
\multiput{$\bullet$} at 2 -2  2 0  2 1.5   /
\multiput{$\bullet$} at 1 -1   -.5 -1 -2 -1   3 -1  4.5 -1  6 -1  /

\multiput{$\cdots$} at 0.5 -1  3.5 -1 /
\plot -2 -1  0 -1 /
\plot 1 -1  2 0  3 -1 /
\plot 4 -1  6 -1 /

\plot 1 -1  2 -2  3 -1 /
\setdots <.5mm>
\plot 1.5 -2  1.5 1.5 /
\plot 2.5 -2  2.5 1.5 /
\circulararc 180 degrees from 1.5 -2  center at 2 -2 
\circulararc 180 degrees from 2.5 1.5  center at 2 1.5 
\put{$\Bbb A_1\sqcup \Bbb D_m$} at -3 -1.8 
\put{$\ssize P(x(2))$} at -2 -.3
\put{$\ssize I(x(2))$} at 6 -.3
\endpicture} at 11 0 
\endpicture}
$$
\bigskip
$$
{\beginpicture
\setcoordinatesystem units <.5cm,.5cm> 
\put{\beginpicture
\put{$\ssize P(x)$} at -1 0 
\put{$\ssize I(x)$} at 5 0 
\multiput{$\bullet$} at 0 0  1 1  2 2  1 -1  2 0  3 1  2 -2  3 -1  4 0 /
\plot 0 0  2 2  4 0  2 -2  0 0 /
\plot 1 1  3 -1 /
\plot 1 -1  3 1 /
\setdots <.5mm>
\plot 1.5 -2  1.5 2 /
\plot 2.5 -2  2.5 2 /
\circulararc 180 degrees from 1.5 -2  center at 2 -2 
\circulararc 180 degrees from 2.5 2  center at 2 2 
\put{$\Bbb A_5$} at -1.5 -1.8 
\endpicture} at 0 0 
\put{\beginpicture
\put{$\ssize P(x)$} at -1 4 
\put{$\ssize I(x)$} at 9 4
\multiput{$\bullet$} at 0 4  1 3  2 2  3 1  4 0  5 1  6 2  6 4  7 3  8 4  2 4  3 3
   4 2  4 4  5 3  /
\plot 0 4  4 0  8 4 /
\plot 1 3  2 4  3 3  5 1 /
\plot 3 1  5 3  6 4  7 3 /
\plot 2 2  4 4  6 2 /
\setdots <.5mm>
\plot 3.5 0  3.5 4 /
\plot 4.5 0  4.5 4 /
\circulararc 180 degrees from 3.5 0  center at 4 0 
\circulararc 180 degrees from 4.5 4  center at 4 4 
\put{$\Bbb D_6$} at .5 0.2 
\endpicture} at 11 0 
\put{\beginpicture
\multiput{$\bullet$} at 0 0  1 1  2 2  3 3  4 4  5 5  4 2  5 3  6 2  6 4  7 1  7 3  8 0  8 2 
    8 4  9 1  9 3  10 2  10 4  11 3  11 5  12 2  12 4  13 3  14 2  15 1  16 0  /
\put{$\ssize P(x)$} at -1 0 
\put{$\ssize I(x)$} at 17 0 
\plot  0 0  5 5  9 1 /
\plot  3 3  4 2  5 3  6 4 /
\plot  4 4  8 0  12 4  /
\plot 6 2  7 3  8 4  9 3  10 2 /
\plot 7 1  11 5  16 0  /
\plot 10 4  11 3  12 2  13 3 /
\put{$\Bbb E_7$} at -1 2 
\setdots <.5mm>
\plot 7.5 0  7.5 4 /
\plot 8.5 0  8.5 4 /
\circulararc 180 degrees from 7.5 0  center at 8 0 
\circulararc 180 degrees from 8.5 4  center at 8 4 
\endpicture} at 5 -7 
\endpicture}
$$
the dotted lines show the unique antichain triple. 
In the case $\Bbb A_1\sqcup \Bbb D_m$, the number of 
elements of $H(Q,x(1),x(2))$ is $2m-2$, and the poset is self-dual.
We also should remark that we have poset embeddings 
$$
 H(\Bbb A_5,x) \subset H(\Bbb D_6,x) \subset H(\Bbb E_7,x).
$$
(Note that the hammock $H(\Bbb E_7,x)$ was used
already in [R4], Example 7, with a similar aim, namely to describe the largest 
indecomposable representation of
a quiver of type $\Bbb E_8.$) 
	\bigskip 
Actually, if we want to stress that these posets arise as subquivers of a translation quiver
of the form $\Bbb Z Q$, we should draw them in the cases 
$\Bbb A_1\sqcup\Bbb D_m,\ \Bbb D_6,\ \Bbb E_7$ differently, namely as follows:
$$
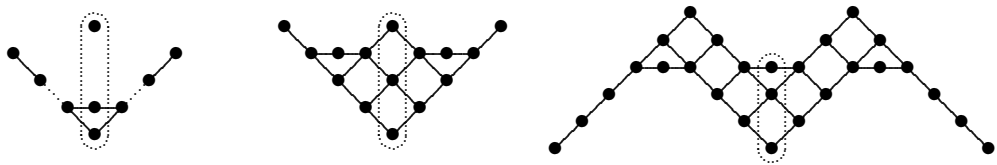
{\beginpicture
\setcoordinatesystem units <.3cm,.3cm> 
\put{\beginpicture
\multiput{$\bullet$} at 0 0  -1 1  2 2  1 -1  2 -1  5 1  2 -2  3 -1  4 0 /
\plot -1 1  0 0 /
\plot  3 -1 1 -1 2 -2  3 -1 /
\plot  4 0  5 1 /
\setdots <.7mm>
\plot 0 0  1 -1 /
\plot 4 0  3 -1 /
\setdots <.5mm>
\plot 1.5 -2  1.5 2 /
\plot 2.5 -2  2.5 2 /
\circulararc 180 degrees from 1.5 -2  center at 2 -2 
\circulararc 180 degrees from 2.5 2  center at 2 2 
\endpicture} at 0 0 
\put{\beginpicture
\multiput{$\bullet$} at 0 4    1 3  2 2  3 1  4 0  5 1  6 2  6 3  7 3  8 4  2 3  3 3
   4 2  4 4  5 3  /
\plot 0 4  4 0  8 4 /
\plot 1 3  3 3  5 1 /
\plot 3 1  5 3  7 3 /
\plot 2 2  4 4  6 2 /
\setdots <.5mm>
\plot 3.5 0  3.5 4 /
\plot 4.5 0  4.5 4 /
\circulararc 180 degrees from 3.5 0  center at 4 0 
\circulararc 180 degrees from 4.5 4  center at 4 4 
\endpicture} at 11 0 
\put{\beginpicture
\multiput{$\bullet$} at 0 0  1 1  2 2  3 3  4 4  5 5  4 3  5 3  6 2  6 4  7 1  7 3  8 0  8 2 
    8 3  9 1  9 3  10 2  10 4  11 3  11 5  12 3  12 4  13 3  14 2  15 1  16 0  /
\plot  0 0  5 5  9 1 /
\plot  3 3  5 3  6 4 /
\plot  4 4  8 0  12 4  /
\plot 6 2  7 3  9 3  10 2 /
\plot 7 1  11 5  16 0  /
\plot 10 4  11 3  13 3 /
\setdots <.5mm>
\plot 7.5 0  7.5 3 /
\plot 8.5 0  8.5 3 /
\circulararc 180 degrees from 7.5 0  center at 8 0 
\circulararc 180 degrees from 8.5 3  center at 8 3 
\endpicture} at 25 0 
\endpicture}
$$
	\bigskip 
If $x(1),\dots,x(t)$ is a special vertex set of the quiver $Q$, the 
three representations $A(1),A(2),A(3)$ in $\rep Q$ which form an antichain triple in 
$H(Q,x(1),\dots,x(t))$ will be called the {\it special triple} of $(Q,x(1),\dots,x(t))$.
The hammocks $H(Q(i),x(i))$ show: 

{\it If $A(1),A(2),A(3)$ is the special triple of $(Q,x(1),\dots,x(t))$, then for $i\neq j$,
there is no path in $\rep Q$ from $A(i)$ to $A(j)$.} It follows that $A(1),A(2),A(3)$ 
is not only an antichain triple in $H(Q,x(1),\dots,x(t))$, but we even have
$\Hom(A(i),A(j)) = 0$ for $i\neq j$ (it is an antichain triple in $\rep Q$), as well as
$\Ext^1(A(i),A(j)) = 0$ (by the previous argument for $i\neq j$, but, of course, also for $i=j$).
		   \bigskip 
Let us exhibit the dimension vectors of such an antichain triple for some fixed orientation
of the quiver $Q$. In case $Q = \Bbb A_1\sqcup\Bbb A_1\sqcup\Bbb A_1$, 
the special triple 
is given by the simple representations. In case $Q = Q(1)\sqcup Q(2)$, where $Q(1) = \Bbb A_1$
and $Q(2)$ is the quiver of type $\Bbb D_m$ with subspace orientation, 
the first representation
$A(1)$ is the simple representation of the component of type $\Bbb A_1$, the remaining two
representations $A(2),A(3)$ are as follows:
$$
{\beginpicture
\setcoordinatesystem units <.5cm,.5cm> 
\put{\beginpicture
\put{$\Bbb D_m$} at -2 0
\put{$\cdots$} at 3.5 0
\multiput{$\circ$} at 0 1  0 -1  1 0  2 0  5 0 /
\arr{0.3 0.7}{0.7 0.3}
\arr{0.3 -.7}{0.7 -.3}
\arr{1.7 0}{1.3 0}
\arr{2.7 0}{2.3 0}
\arr{4.7 0}{4.3 0}
\put{$A(2) \qquad \smallmatrix 1 \cr
                     &1 & 1 & 1 & \cdots & 1 \cr
                   0 \endsmallmatrix$} at 11 1
\put{$A(3) \qquad  \smallmatrix 0 \cr
                     &1 & 1 & 1 & \cdots & 1 \cr
                   1 \endsmallmatrix$} at 11 -1

\endpicture} at 0 0
\endpicture}
$$
Finally we look at quivers of type $\Bbb A_5, \Bbb D_6, \Bbb E_7$, again dealing with
a subspace orientation:
$$
{\beginpicture
\setcoordinatesystem units <1cm,.75cm> 
\put{\beginpicture
\put{$\Bbb A_5$} at 0 2
\put{\beginpicture
  \setcoordinatesystem units <.5cm,.5cm> 
  \multiput{$\circ$} at 0 0  1 0   3 0  4 0 /
  \put{$x$} at 2 0
  \put{$\bigcirc$} at 2 0
  \arr{0.2 0}{0.8 0} 
  \arr{1.2 0}{1.8 0} 
  \arr{2.8 0}{2.2 0} 
  \arr{3.8 0}{3.2 0} 
  \endpicture} at 0 1 
\put{$A(1)$} at -2 -.3
\put{$A(2)$} at -2 -1.3
\put{$A(3)$} at -2 -2.3
\put{$\smallmatrix 1 & 1 & 1 & 0 & 0\endsmallmatrix$} at 0 -.3
\put{$\smallmatrix 0 & 1 & 1 & 1 & 0\endsmallmatrix$} at 0 -1.3
\put{$\smallmatrix 0 & 0 & 1 & 1 & 1\endsmallmatrix$} at 0 -2.3
\endpicture} at 0 0
\put{\beginpicture
\put{$\Bbb D_6$} at 0 2
\put{\beginpicture
  \setcoordinatesystem units <.5cm,.5cm> 
  \multiput{$\circ$} at 1 1  1 0  2 0  3 0  4 0 /
  \put{$x$} at 0 0 
  \put{$\bigcirc$} at 0 0 
  \arr{0.2 0}{0.8 0} 
  \arr{1.8 0}{1.2 0} 
  \arr{2.8 0}{2.2 0} 
  \arr{3.8 0}{3.2 0} 
  \arr{1 0.8}{1 0.2}
  \endpicture} at 0 1
\put{$\smallmatrix   & 0 \cr
                   1 & 1 & 1 & 1 & 0 \endsmallmatrix$} at 0 -.2 
\put{$\smallmatrix   & 1 \cr
                   1 & 2 & 2 & 1 & 1\endsmallmatrix$} at 0 -1.2
\put{$\smallmatrix   & 1 \cr
                   1 & 1 & 0 & 0 & 0\endsmallmatrix$} at 0 -2.2
\endpicture} at 4 0
\put{\beginpicture
\put{$\Bbb E_7$} at 0 2
\put{\beginpicture
  \setcoordinatesystem units <.5cm,.5cm> 
  \multiput{$\circ$} at -1 0  0 0  1 1  1 0  2 0  3 0   /
  \put{$x$} at 4 0 
  \put{$\bigcirc$} at 4 0 
  \arr{-.8 0}{-.2 0} 
  \arr{0.2 0}{0.8 0} 
  \arr{1.8 0}{1.2 0} 
  \arr{2.8 0}{2.2 0} 
  \arr{3.8 0}{3.2 0} 
  \arr{1 0.8}{1 0.2}
  \endpicture} at 0 1
\put{$\smallmatrix &  & 1 \cr
                   1 & 1 & 2 & 2 & 1 & 1 \endsmallmatrix$} at 0 -.2 
\put{$\smallmatrix &  & 2 \cr
                   1 & 2 & 3 & 2 & 2 & 1\endsmallmatrix$} at 0 -1.2 
\put{$\smallmatrix &  & 0 \cr
                   0 & 1 & 1 & 1 & 1 & 1\endsmallmatrix$} at 0 -2.2
\endpicture} at 8 0
\endpicture}
$$
	\bigskip 
{\bf 3.4. The one-point extension for a special vertex set.}
     \medskip 
{\bf The one-point extension for a special vertex set.} {\it Let 
$x(1),\dots,x(t)$ be a special vertex set of the acyclic
quiver $Q$. Then $\Delta = Q[x(1),\dots,x(t)]$ is 
a Dynkin quiver of type $\Bbb D_n$ or $\Bbb E_m$ and the extension vertex of 
$\Delta$ is its exceptional vertex and is a source. 

Conversely, let $\Delta$ be a Dynkin quiver of type $\Bbb D_n$ or $\Bbb E_m$ and $y$ its
exceptional vertex. Recall that $\Delta'$ is obtained from $\Delta$ by deleting $y$ and let $x(1),\dots,x(t)$
be the neighbors of $y$ in $\Delta$. Then 
$x(1),\dots,x(t)$ is a special vertex set of the quiver $\Delta'$.}
		  \medskip
Proof: Direct verification. 
       \medskip 
{\bf Corollary.} {\it The Dynkin quivers of type $\Bbb D_n$ or $\Bbb E_m$ are
representation-finite.}
	\medskip
Proof. Let $\Delta$ be of a quiver of type $\Bbb D_n, \Bbb E_6, \Bbb E_7$ or $\Bbb E_8$. Changing if
necessary the orientation, we may assume that the vertex $y$ is a source. 
Thus, $\Delta$ is the one-point extension $\Delta'[x(1),\dots,x(t)]$
of $\Delta'$, where $x(1),\dots,x(t)$ is a special vertex set. Its representation type is determined by the
poset $H(\Delta',x(1),\dots,x(t))$, see for example [R2] or [R3], sections 2.5 and 2.6.
Since $H(\Delta',x(1),\dots,x(t))$ has a unique antichain
triple, $\Delta'[x(1),\dots,x(t)]$ has to be representation-finite.
	\medskip
{\bf 3.5.}
Actually, the poset $H(\Delta',x(1),\dots,x(t))$ provides an easy way to construct explicitly all the
indecomposable representations of $\Delta$. Let us
distinguish 6 different kinds:
\item{(1)} The indecomposable representations $Y$ of $\Delta$ with $N_y = 0$ and 
  $N_{x(i)} = 0$
   for all $1\le i \le t$, these are just the 
 indecomposable representations of the quiver $\Delta''$.
\item{(2)} The indecomposable representations $Y$ of $\Delta$ with 
  $N_y = 0$ and  $N_{x(i)} \neq 0$
  for some $1\le i \le t;$ these are the representations which belong to the
  hammocks $H(\Delta',x(i))$ with $1\le i \le t.$
\item{(3)} For any representation $N$ in $H(\Delta',x(1),\dots,x(t)),$ there is an indecomposable representation
  $Y$ of $\Delta$ with $Y|\Delta' = N$ and $\dim \overline Y_y = 1.$
\item{(4)} For any antichain pair $B(1),B(2)$ in $H(\Delta',x(1),\dots,x(t))$, 
  there is an indecomposable representation $Y$ of $\Delta$
  with $Y|\Delta' = B(1)\oplus B(2)$ and $\dim Y_y = 1.$
\item{(5)} Two representations $Y$ of $\Delta$ with $Y|\Delta'  
  = A(1)\oplus A(2)\oplus A(3)$
  such that $\dim Y_y$ is equal to $1$ or $2$.
\item{(6)} The simple representation $S(y)$.
 	   \medskip 
Here is the table which lists the number of indecomposable representations $Y$ of
$\Delta$ of the various kinds.
Given a representation $N$ of any quiver, 
we denote by $\m(N)$ the number of direct summands
when $N$ is written as a direct sum of indecomposable representations.  The upper three
rows name the types of $\Delta$, of $\Delta' = \Delta\setminus\{y\}$ and of
$\Delta'' = \Delta\setminus\{x,y\}$.
$$
{\beginpicture
\setcoordinatesystem units <1cm,.7cm>
\put{} at -4.5 3
\put{$\Delta$} at -3 3
\put{$\Bbb D_n\ (n\ge 4)$} at -.7 3 
\put{$\Bbb E_6$} at 1 3
\put{$\Bbb E_7$} at 2 3  
\put{$\Bbb E_8$} at 3 3

\put{$\Delta'$} at -3 2
\put{$\Bbb D_{n-2}$} at -.7 2 
\put{$\Bbb A_5$} at 1 2
\put{$\Bbb D_6$} at 2 2  
\put{$\Bbb E_7$} at 3 2

\put{$\Delta''$} at -3 1
\put{$\Bbb D_{n-3}$} at -.7 1 
\put{$\Bbb A_2\!\sqcup\!\Bbb A_2$} at 1 1
\put{$\Bbb A_5$} at 2 1  
\put{$\Bbb E_6$} at 3 1

\put{$Y \in \rep \Delta''$} at -4 0 
\put{$Y \in \rep \Delta'\setminus \rep \Delta''$} at -4 -1 
\put{$\m(Y|\Delta') = 1,\ Y_y\neq 0$} at -4 -2 
\put{$\m(Y|\Delta') = 2$} at -4 -3 
\put{$\m(Y|\Delta') = 3$} at -4 -4 
\put{$Y|\Delta' = 0$} at -4 -5

\put{(1)} at -6.5 0 
\put{(2)} at -6.5 -1 
\put{(3)} at -6.5 -2 
\put{(4)} at -6.5 -3 
\put{(5)} at -6.5 -4 
\put{(6)} at -6.5 -5 
\put{sum} at -6.5 -6 

\put{$(n-4)(n-3)$} at -.7 0
\put{$2n-5$} at -.7 -1 
\put{$2n-5$} at -.7 -2 
\put{$2n-5$} at -.7 -3 
\put{2} at -.7 -4 
\put{1} at -.7 -5 
\put{$(n-1)n$} at -.7 -6

\put{6} at 1 0 
\put{9} at 1 -1 
\put{9} at 1 -2 
\put{9} at 1 -3 
\put{2} at 1 -4 
\put{1} at 1 -5 
\put{36} at 1 -6 
  
\put{15} at 2 0 
\put{15} at 2 -1 
\put{15} at 2 -2 
\put{15} at 2 -3 
\put{2} at 2 -4 
\put{1} at 2 -5 
\put{63} at 2 -6 

\put{36} at 3 0 
\put{27} at 3 -1 
\put{27} at 3 -2 
\put{27} at 3 -3 
\put{2} at 3 -4 
\put{1} at 3 -5 
\put{120} at 3 -6
 
\plot -7 .5  3.5 .5 /
\plot -7 -5.5  3.5 -5.5 /
\plot -2 3.4  -2 -6.4 /
\endpicture}
$$
The Dynkin graphs $\Bbb D_n$ are usually considered only for $n\ge 4$. But our
induction procedure requires to define also $\Bbb D_n$ with $n = 1,2,3$: we set
$\Bbb D_1 = \emptyset,\ \Bbb D_2 = \Bbb A_1\sqcup \Bbb A_1$ and $\Bbb D_3 = \Bbb A_3.$
	       \medskip 
It is clear that the numbers in the rows (2) and (3) coincide: this number is
just the cardinality of the vertex set of the poset  $H(\Delta',x(1),\dots,x(t))$. However,
it seems to be a curious coincidence that in all cases the number of
vertices of the poset $H(\Delta',x(1),\dots,x(t))$  is equal to the number of antichain pairs
in $H(\Delta',x(1),\dots,x(t))$ (thus, that the number of the representations in row (4) 
is equal to the number in row (2) and (3)).
	 \medskip 
The indecomposable representation $Y$ in (5) 
with $Y|\Delta' = A(1)\oplus A(2)\oplus A(3)$
and $\dim Y_y = 3$ is the maximal indecomposable representation of $\Delta$, we usually
label it $M$.
	\bigskip 
{\bf 3.6. The special antichain triple of a Dynkin quiver of type $\Bbb D_n$ and $\Bbb E_m$.}
 Let $\Delta$ be a Dynkin quivers of type $\Bbb D_n$ or $\Bbb E_m$. Let $y$ be its  exceptional vertex. By definition,
$\Delta'$ is obtained from $\Delta$ by deleting $y$ and we denote by $x(1),\dots,x(t)$
the neighbors of $y$ in $\Delta$. 
Then $x(1),\dots,x(t)$ is a special vertex set of the quiver $\Delta'$, 
thus there is a unique antichain triple $A(1),A(2),A(3)$ in 
$H(\Delta',x(1),\dots,x(t)).$ 
We call $A(1),A(2),A(3)$ the {\it special antichain triple} of $\Delta$. 
The previous considerations show that the special antichain triple is 
the unique antichain triple $A(1),A(2),A(3)$ in 
$\rep \Delta'$ such that for $1\le i \le 3$, there is at least one 
vertex $x$ in the set $\{x(1),\dots,x(t)\}$ with
$A(i)_{x} \neq 0.$ 
To say that $A(i)$ is a representation of $\Delta'$ is the same as to say that it is
a representation of $\Delta$ with $A(i)_y = 0$. 
In this way, we have shown the main assertion of Theorem 3.1. 

Of course, we know that for any $i$, there
is just one vertex $x$ in the set  $\{x(1),\dots,x(t)\}$
with $A(i)_{x)} \neq 0$ and for this vertex $x$, we have $\dim A(i)_{x} = 1.$ 
Also, we know that for $i\neq j$ there is no path in $\rep \Delta'$ from $A(j)$ to $A(i)$.
In particular, we have $\Ext^1(A(i),A(j)) = 0.$  Since also $\Ext^1(X,X) = 0$ for 
any indecomposable representation of $\Delta'$, 
we see that $\Ext^1(A(i),A(j)) = 0$ for all $i,j$. 

In order to complete the proof of Theorem 3.1, it remains to be seen that $M|\Delta' =
A(1)\oplus A(2)\oplus A(3)$. This will be shown now.
     \medskip 
Let us consider the special antichain triple in more detail. 
    \medskip
{\bf Theorem.} {\it Let $\Delta$ be a Dynkin quiver of type $\Bbb D_n$ or $\Bbb E_m$ with exceptional
vertex $y$ and special antichain triple $A(1),A(2),A(3).$ The
$\Ext$-quiver $C(\Delta)$ of 
$$
 \Cal C(\Delta) = \Cal E(A(1),A(2),A(3),S(y))
$$ has the following form
$$
{\beginpicture
\setcoordinatesystem units <2cm,1cm> 
\put{$A(1)$} at 0 1 
\put{$A(2)$} at 0 0 
\put{$A(3)$} at 0 -1 
\put{$S(y)$} at 1 0 
\plot 0.7 0.3 0.3 0.8 /
\plot 0.7 0 0.3 0 /
\plot 0.7 -.3 0.3 -.8 /
\endpicture}
$$
It is of type $\Bbb D_4$ and $\Cal C(\Delta)$ contains the maximal indecomposable representation $M$
of $\Delta$.}
   \medskip 
As an element of $\Cal C$, the representation $M$ has dimension vector 
$$
 \bdim_{\Cal C} M = \smallmatrix 1 \cr
                                 1 & 2 \cr
                                 1 \endsmallmatrix\ ,
$$
thus $M|\Delta' = A(1)\oplus A(2)\oplus A(3)$ and $\dim M_y = 2.$
     \medskip
Proof. In order to determine the quiver $C(\Delta)$, we first note that the indecomposable
representations of a Dynkin quiver have no self-extensions, thus there is no loop in
$C(\Delta)$. 
Second, since the representations $A(1), A(2),A(3)$ are $\Ext$-orthogonal, there is no
arrow between these representations. Finally, for any $i$, there is a unique neighbor $x$ of $y$
with $A(i)_x \neq 0$, and one has $\dim A(i)_x = 1$. In case the arrow between $x$ and $y$ points to $y$,
we have $\dim\Ext^1(A(i),S(y)) = 1$ and $\Ext^1(S(y),A(i)) = 0.$ Otherwise we
have $\dim\Ext^1(S(y),A(i)) = 1$ and $\Ext^1(A(i),S(y)) = 0.$

It remains to look at the maximal indecomposable representation $M$ of $\Delta$. It is the only
indecomposable representation of $\Delta$ with $\dim M_y = 2$ 
(since $y$ is the exceptional vertex of $\Delta$). 
Since $\Cal C(\Delta)$ is of type $\Bbb D_4$, there is an indecomposable object $X$ in
$\Cal C(\Delta)$ with dimension vector 
$ \bdim_{\Cal C} X = \smallmatrix 1 \cr
                                 1 & 2 \cr
                                 1 \endsmallmatrix\ .
$
In particular, $X$ is an indecomposable representation with $\dim X_y = 2$, thus $X$ is isomorphic to $M$.
   \medskip 
As we have seen in the proof, there is the following {\bf recipe} 
for obtaining the  orientation of the quiver $C(\Delta)$: 
For $1\le i \le 3$, the orientation of the edge between $[A(i)]$ and $[S(y)]$ is the same
as the orientation between $x$ and $y$ in $\Delta$, where $A(i)_x \neq 0.$
   \bigskip 
We call $\Cal C(\Delta)$ the {\it core} of $\rep\Delta$.
   \medskip 
{\bf Corollary.} {\it If $\Delta$ is a Dynkin quiver of type $\Bbb E_m$, then
the core of $\rep\Delta$ is the smallest thick subcategory of
$\rep\Delta$ which contains $M$ and $S(y)$.}
	\medskip
Proof. Of course, the core $\Cal C(\Delta)$ contains both $M$ and $S(y)$.
Conversely, assume that $\Cal T$ is a thick subcategory which contains $M$ and $S(y)$.
If $y$ is a source, then $A(1)\oplus A(2)\oplus A(3)$ is the kernel of a map
$M \to S(y)^2$, thus $A(1),A(2),A(3)$ belong to $\Cal T.$ A similar agument works if 
$y$ is a sink.
	     \bigskip 
{\bf 3.7. The corresponding Euclidean quivers.}
     \medskip
We consider now pairs $\Delta \subset \widetilde\Delta$,
where $\Delta$ is a Dynkin quiver of type $\Bbb D_n$ or $\Bbb E_m$ with exceptional vertex $y$,
and $\widetilde\Delta$ is a corresponding Euclidean quiver, thus there is a unique vertex $z$ 
of $\widetilde\Delta$ outside of $\Delta$ and $z$ is a neighbor of $y$
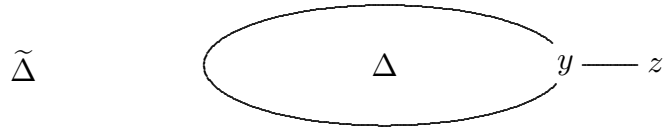
$$  
 {\beginpicture
 \setcoordinatesystem units <1cm,.7cm> 
\multiput{} at 1 1.2  1 -1.2 /
\put{} at -1.5 0 
\put{$y$} at 3 0
\put{$z$} at 4 0
\plot 3.8 0  3.2 0 /
\put{$\Delta$} at 1 0
\put{$\widetilde\Delta$} at -3 0
\ellipticalarc axes ratio 3:1 -320 degrees from 2.9 -.3 center at 1 0  
\endpicture} 
$$
(the orientation of the arrow between $y$ and $z$ is not fixed). 
     \medskip 
{\bf Theorem.} {\it Let $A(1),A(2),A(3)$ be the special antichain triple for $\Delta$ and let
$$
 \Cal C_z(\widetilde\Delta) = \Cal E(A(1),A(2),A(3),S(y),S(z)).
$$
Then, the $\Ext$-quiver of $\Cal C_z(\widetilde\Delta)$
has the following form 
$$
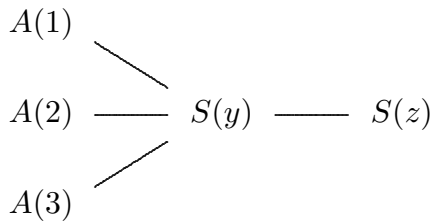
{\beginpicture
\setcoordinatesystem units <2cm,1cm> 
\put{$A(1)$} at 0 1 
\put{$A(2)$} at 0 0 
\put{$A(3)$} at 0 -1 
\put{$S(y)$} at 1 0 
\put{$S(z)$} at 2 0 
\plot 0.7 0.3 0.3 0.8 /
\plot 0.7 0 0.3 0 /
\plot 0.7 -.3 0.3 -.8 /
\plot 1.7 0 1.3 0 /
\endpicture}
$$
It is of type $\widetilde{\Bbb D}_4$ and the 
subcategory $\Cal C_z(\widetilde\Delta)$ contains an element 
from any $\tau$-orbit of simple regular representations of $\widetilde \Delta$.} 
     \medskip
We should remark that actually $\Cal C = \Cal C_z(\widetilde\Delta)$ 
only depends on $y$, not on $z$ itself:
For quivers of type $\widetilde{\Bbb D}_4$, there is just one subcategory $\Cal C$,
namely $\rep\Delta$ itself (but
four different choices of $z$). For quivers of type $\widetilde{\Bbb D}_n$, with $n\ge 5$,
there are two subcategories of the form $\Cal C$ (but again 
four different choices of $z$). For the quivers of type $\widetilde{\Bbb E}_m$, the possible vertices $z$
correspond bijectively to the vertices $y$ (there are three choices for $\widetilde{\Bbb E}_6$, 
two choices for $\widetilde{\Bbb E}_7$, and only one choice for $\widetilde{\Bbb E}_8$).

    \medskip 
Proof of the Theorem. We know from Theorem 3.6 the shape of the $\Ext$-quiver $C(\Delta)$.
Clearly, the $\Ext$-quiver for $\Cal C_z(\widetilde\Delta)$ is obtained from $C(\Delta)$
by adding the vertex $[S(z)]$ and an arrow between $[S(y)]$ and $[S(z)]$.
The orientation of the edge between $[S(y)]$ and $[S(z)]$ is the same as 
the orientation of the edge between $y$ and $z$ in
the quiver $\Delta$. 

It remains to show that any $\tau$-orbit of a simple regular representation $N$ of 
$\widetilde\Delta$
contains a representation in $\Cal C_z(\widetilde\Delta)$. 
If $N$ is homogeneous (that means $\tau N = N$), then $N|\Delta$ is a direct sum of copies of 
the maximal indecomposable representation $M$.
Since $M$ belongs to $\Cal C(\Delta),$ we see that $N$ belongs to $\Cal C_z(\widetilde \Delta).$
Thus, we now may assume that $N$ is exceptional. We distinguish according to the
number $t$ of connected components of $\Delta'$.
       \medskip 
{\bf The types $\widetilde{\Bbb D}_4$} (the case $t=3$). 
In case $\widetilde\Delta$ is of type $\widetilde{\Bbb D}_4$, we have
$\Cal C_z(\widetilde\Delta) = \rep\widetilde\Delta$ (since $\rep\widetilde\Delta$ is its
only thick subcategory of rank 5, 
see for example [R5]), thus nothing has to be shown in this case.
    \medskip
{\bf The types $\widetilde{\Bbb E}_m$} (the cases $t=1$). 
Here, $\widetilde\Delta$ has the following shape:
$$  
 
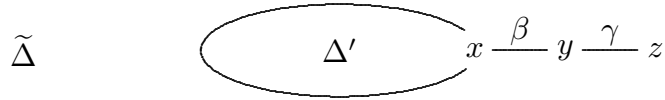
{\beginpicture
 \setcoordinatesystem units <1cm,.7cm> 
\multiput{} at 1 1.2  1 -1.2 /
\put{} at -1.5 0 
\put{$x$} at 2 0
\put{$y$} at 3 0
\put{$z$} at 4 0
\plot 2.8 0  2.2 0 /
\plot 3.8 0  3.2 0 /
\put{$\beta\strut$} at 2.5 0.3
\put{$\gamma\strut$} at 3.5 0.3
\put{$\Delta'$} at 0.5 0
\put{$\widetilde\Delta$} at -3 0
\ellipticalarc axes ratio 3:1 -310 degrees from 1.9 -.3 center at .5 0  
\endpicture} 
$$
with $\{x\}$ a special vertex set of $\Delta'$. Using BGP reflection functors at vertices $a\in\widetilde\Delta_0
\setminus\{x,y\}$, we may change the orientation of all the arrows different from $\beta$. 
Under such a change of orientation, the image of the special antichain triple of $\Delta$ will again 
be a special antichain triple. In addition, using duality, we may replace the direction of $\beta$ by
the opposite direction. Altogether, we see that it is sufficient to look at a single orientation, say
at the subspace orientation as shown below (with $\Delta'$ the full subquiver with vertex set 
$\widetilde\Delta_0\setminus\{y,z\}$). In particular, the edges $\beta$ and $\gamma$ are oriented as
follows:
$$  
 
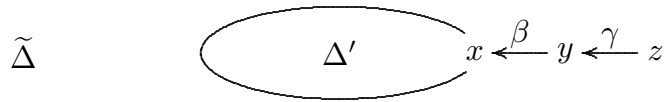
{\beginpicture
 \setcoordinatesystem units <1cm,.7cm> 
\multiput{} at 1 1.2  1 -1.2 /
\put{} at -1.5 0 
\put{$x$} at 2 0
\put{$y$} at 3 0
\put{$z$} at 4 0
\arr{2.8 0}{2.2 0}
\arr{3.8 0}{3.2 0}
\put{$\beta\strut$} at 2.5 0.3
\put{$\gamma\strut$} at 3.5 0.3
\put{$\Delta'$} at 0.5 0
\put{$\widetilde\Delta$} at -3 0
\ellipticalarc axes ratio 3:1 -310 degrees from 1.9 -.3 center at .5 0  
\endpicture} 
$$

For such an orientation, the $\Ext$-quiver $R$ of $\Cal C = \Cal C_z(\widetilde\Delta)$ is as follows:
$$
{\beginpicture
\setcoordinatesystem units <1.5cm,.7cm> 
\put{$R$} at -1.5 0
\put{$A(1)$} at 0 1 
\put{$A(2)$} at 0 0 
\put{$A(3)$} at 0 -1 
\put{$S(y)$} at 1 0
\put{$S(z)$} at 2 0 
\arr{0.7 0.3}{0.3 0.8}
\arr{0.7 0}{0.3 0}
\arr{0.7 -.3}{0.3 -.8}
\arr{1.7 0}{1.3 0}
\endpicture}
$$

Given a representation $X$ of $\Delta'$, we denote by $\overline X$ the representation of 
$\widetilde \Delta$
with $\overline X|\Delta' = X,$ $\overline X_\beta\:\overline X_y \to \overline X_x = X_x$ the identity map, and $\overline X_z = 0.$ There is the following interesting observation:
{\it Let $X$ be an indecomposable representation of $\Delta'$.
Then $\overline X$ is a regular representation of $\widetilde\Delta$
if and only if $X$ is one of the representations 
$A(1), A(2),A(3)$.} 
       \medskip
As elements of the subcategory $\Cal C$, the representations $\overline{A(i)}$
have the dimension vectors
$$
  \bdim_{\Cal C} \overline {A(1)} = 
  \smallmatrix 1 \cr
              0 & 1 & 0 \cr
              0 \endsmallmatrix, \quad  
  \bdim_{\Cal C} \overline {A(2)} = 
  \smallmatrix 0 \cr
              1 & 1 & 0 \cr
              0 \endsmallmatrix, \quad  
  \bdim_{\Cal C} \overline {A(3)} = 
  \smallmatrix 0 \cr
              0 & 1 & 0 \cr
              1 \endsmallmatrix;
$$
these are simple regular objects of $\Cal C$. 

At the end of section 3.3, the antichain triple $A(1),A(2),A(3)$ has been exhibited for our quivers $\Delta'$.
Thus, the representations $\overline{A(i)}$ are as follows:
$$
{\beginpicture
\setcoordinatesystem units <1cm,.75cm> 
\put{\beginpicture
\put{$\widetilde{\Bbb E}_6$} at -.7 2
\put{\beginpicture
  \setcoordinatesystem units <.6cm,.6cm> 
  \multiput{$\circ$} at 0 0  1 0   3 0  4 0 /
  \put{$y$} at 2 1 
  \put{$z$} at 2 2
  \put{$\ssize\beta$} at 2.3 0.5 
\arr{2 0.7}{2 0.3} 
\arr{2 1.7}{2 1.3} 
  \put{$x$} at 2 0
  \arr{0.2 0}{0.8 0} 
  \arr{1.2 0}{1.8 0} 
  \arr{2.8 0}{2.2 0} 
  \arr{3.8 0}{3.2 0} 
  \endpicture} at 0 1.17 
\put{$\overline{A(1)}$} at -2.5 -.7
\put{$\overline{A(2)}$} at -2.5 -2.1
\put{$\overline{A(3)}$} at -2.5 -3.5
\put{$\smallmatrix & & 0 \cr
                   & & 1 \cr 
         1 & 1 & 1 & 0 & 0\endsmallmatrix$} at 0 -.7
\put{$\smallmatrix  & & 0 \cr
                    & & 1 \cr 
         0 & 1 & 1 & 1 & 0\endsmallmatrix$} at 0 -2.1
\put{$\smallmatrix  & & 0 \cr
                    & & 1 \cr 
        0 & 0 & 1 & 1 & 1\endsmallmatrix$} at 0 -3.5
\endpicture} at -.3 0
\put{\beginpicture
\put{$\widetilde{\Bbb E}_7$} at -1 3
\put{\beginpicture
  \setcoordinatesystem units <.6cm,.6cm> 
  \multiput{$\circ$} at 1 1  1 0  2 0  3 0  4 0 /
  \put{$x\strut$} at 0 0 
  \put{$y\strut$} at -1 0 
  \put{$z\strut$} at -2 0 
  \put{$\ssize\beta$} at -.5 0.3 
  \arr{0.2 0}{0.8 0} 
  \arr{1.8 0}{1.2 0} 
  \arr{2.8 0}{2.2 0} 
  \arr{3.8 0}{3.2 0} 
  \arr{1 0.8}{1 0.2}
  \arr{-.8 0}{-.2 0} 
  \arr{-1.8 0}{-1.2 0} 

  \endpicture} at 0 1.7
\put{$\smallmatrix   &&& 0 \cr
                  0 & 1 & 1 & 1 & 1 & 1 & 0 \endsmallmatrix$} at 0 .3
\put{$\smallmatrix   &&& 1 \cr
                  0 & 1 &  1 & 2 & 2 & 1 & 1\endsmallmatrix$} at 0 -1.2
\put{$\smallmatrix   &&& 1 \cr
                  0 & 1 &  1 & 1 & 0 & 0 & 0\endsmallmatrix$} at 0 -2.7
\endpicture} at 4.5 0
\put{\beginpicture
\put{$\widetilde{\Bbb E}_8$} at -1.6 3
\put{\beginpicture
  \setcoordinatesystem units <.6cm,.6cm> 
  \multiput{$\circ$} at -1 0  0 0  1 1  1 0  2 0  3 0   /
  \put{$\strut x$} at 4 0 
  \arr{-.8 0}{-.2 0} 
  \arr{0.2 0}{0.8 0} 
  \arr{1.8 0}{1.2 0} 
  \arr{2.8 0}{2.2 0} 
  \arr{3.8 0}{3.2 0} 
  \arr{1 0.8}{1 0.2}
\put{$\strut y$} at 5 0 
\put{$\strut z$} at 6 0 
  \put{$\ssize\beta$} at 4.6 0.3 
  \arr{4.8 0}{4.2 0} 
  \arr{5.8 0}{5.2 0} 
  \endpicture} at 0 1.7
\put{$\smallmatrix &  & 1 \cr
                   1 & 1 & 2 & 2 & 1 & 1 & 1 & 0 \endsmallmatrix$} at -.4 .3 
\put{$\smallmatrix &  & 2 \cr
                   1 & 2 & 3 & 2 & 2 & 1 & 1 & 0 \endsmallmatrix$} at -.4 -1.2 
\put{$\smallmatrix &  & 0 \cr
                   0 & 1 & 1 & 1 & 1 & 1 & 1 & 0 \endsmallmatrix$} at -.4 -2.7
\endpicture} at 9 0
\endpicture}
$$

The essential claim is the following:
{\it The representations $\overline{A(i)}$ of $\widetilde\Delta$
are simple regular and belong to pairwise different $\tau$-orbits. 
For $i=1$, the $\tau$-period of $\overline {A(i)}$ in $\rep \widetilde\Delta$ is
$3$, for $i=2$, it is $2$, and for $i=3$ it is $3,4,5$
in the cases $\widetilde E_6, \widetilde E_7, \widetilde E_8$, respectively.}

For a proof, one should write down the corresponding $\tau$-orbits.
Actually, in case we deal with a star quiver (as we do), an indecomposable regular representation
which is thin has to be simple regular (and all but three of the listed representations are thin).
      \medskip 
{\bf The types $\widetilde{\Bbb D}_n$ with $n\ge 5$} (the cases $t=2$).  
We can use a similar
procedure as in the cases $\widetilde {\Bbb E}_m$. We consider the quiver $\widetilde\Delta$ with 
the following orientation:
$$  
 
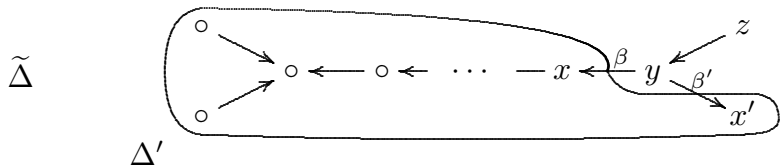
{\beginpicture
 \setcoordinatesystem units <1cm,.5cm> 
\put{$\widetilde\Delta$} at -2 0 
\multiput{$\circ$} at 0 1  0 -1  1 0  2 0   /
\put{$x\strut$} at 4 0  
\put{$y\strut$} at 5 0  
\put{$z\strut$} at 6 1 
\put{$x'\strut$} at 6 -1 
\put{$\cdots$} at 3 0

\arr{0.2 0.8}{0.8 0.2}
\arr{0.2 -.8}{0.8 -.2}
\arr{1.8  0}{1.2 0}
\arr{2.5  0}{2.2 0}
\plot 3.5 0 3.8  0 /
\arr{4.8  0}{4.2 0}
\put{$\ssize\beta$} at 4.65 0.3 
\put{$\ssize\beta'$} at 5.55 -.2 
\arr{5.8  0.8}{5.2 0.2}
\arr{5.2 -.2}{5.8  -.8}
\setquadratic
\plot 0 -1.4  -.3 -1  -.4 0  -.3 1  0 1.4  3.5 1  4.5 0  4.7 -.4  5 -.5  5.5 -.5  
    6 -.5  6.3 -.6  6.4 -1.05  6.3  -1.3  6 -1.4  3 -1.5  0 -1.4 /
\put{$\Delta'$} at -.6 -1.8 

\endpicture}
$$
(since it is sufficient to work up to duality, we may deal with a fixed orientation of the arrow $\beta$
between $x$ and $y$ and we may use BGP reflection functors at vertices $a\in\widetilde\Delta_0
\setminus\{x,y\}$, in order to change the orientation of all the arrows different from $\beta$). 

The special antichain triple of $\Delta'$ is given by the following representations of $\widetilde\Delta$:
$$
 \bdim A(1) = \smallmatrix 1 &&&&& 0 \cr
                             & 1 & \cdots & 1 &0 \cr
                           0 &   &        &  & \ & 0 \endsmallmatrix\ , \quad
 \bdim A(2) = \smallmatrix 0 &&&&& 0 \cr
                             & 1 & \cdots & 1 &0 \cr
                           1 &   &        &  & \ & 0 \endsmallmatrix\ , \quad 
 \bdim A(3) = \smallmatrix 0 &&&&& 0 \cr
                             & 0 & \cdots & 0 &0 \cr
                           0 &   &        &  & \ & 1 \endsmallmatrix\ .
$$
Thus, it turns out that the $\Ext$-quiver of $\Cal C = \Cal C_z(\Delta)$ is again $R$ and we
consider again the representations $\overline{A(i)}$ of $\widetilde\Delta$ which coincide on $\Delta'$ with
$A(i)$, with $\overline{A(1)}_\beta,\ \overline{A(2)}_\beta,\ \overline{A(3)}_{\beta'}$ 
the identity map and with $\overline{A(i)}_z = 0$:
$$
 \bdim \overline{A(1)} = \smallmatrix 1 &&&&& 0 \cr
                             & 1 & \cdots & 1 &1 \cr
                           0 &   &        &  & \ & 0 \endsmallmatrix\ , \quad
 \bdim \overline{A(2)} = \smallmatrix 0 &&&&& 0 \cr
                             & 1 & \cdots & 1 &1 \cr
                           1 &   &        &  & \ & 0 \endsmallmatrix\ , \quad 
 \bdim \overline{A(3)} = \smallmatrix 0 &&&&& 0 \cr
                             & 0 & \cdots & 0 &1 \cr
                           0 &   &        &  & \ & 1 \endsmallmatrix
$$
Again, we see: {\it The representations $\overline{A(i)}$ of $\widetilde\Delta$
are simple regular and belong to pairwise different $\tau$-orbits. 
The $\tau$-period of $\overline {A(1)}$ and $\overline {A(2)}$ in $\rep \widetilde\Delta$ is $2$, 
the $\tau$-period of $\overline {A(3)}$ in $\rep \widetilde\Delta$ is $2$ is $n-2.$}
    \bigskip
We call $\Cal C_y(\widetilde \Delta)$ the {\it $y$-core} of $\rep\widetilde\Delta$.
   \medskip 
{\bf Corollary.} {\it If $\Delta$ is a Dynkin quiver of type $\Bbb E_m$, then 
the $y$-core of $\rep\widetilde\Delta$ is the smallest 
thick subcategory of
$\rep\widetilde\Delta$ which contains $M, S(y),S(z)$.}
		       \bigskip
{\bf 3.8. The quivers $\Delta''$ for $\Delta$ a Dynkin quiver of type $\Bbb D_n$ and
$\Bbb E_m$.}
	\medskip 
{\bf Theorem.} {\it Let $\Delta$ be a Dynkin quiver and $M$ its maximal indecomposable
representation. If $\Delta$ is of type $\Bbb D_4,$ then $M|\Delta'' = 0$.
If $\Delta$ is of type $\Bbb D_n$ with $n\ge 5$, then $M|\Delta''$ is the direct sum
of an antichain pair. If $\Delta$ is of type $\Bbb E_m$,
then $M|\Delta''$ is the direct sum of four indecomposable representations $U, V, U', V'$
with $\Hom(U,V) \neq 0,\ \Hom(U',V') \neq 0$ and
$\dim\End(M|\Delta'') = 6.$}
	\medskip
The proof is similar to the proof of Theorem 3.1. We indicate the main
steps. Details are left to the reader.
	\medskip
{\bf Proposition 1.} {\it 
Let $Q(1),\dots, Q(t)$ be Dynkin quivers and let $Q$ be the disjoint union
of the quivers $Q(i)$. Let $x(i)$ be a vertex of $Q(i)$, for $1\le i \le s.$ Then
$H(Q,w(1),\dots,w(s)$ has precisely one antichain pair
if and only if one of the following 
cases holds: $Q = \Bbb A_1\sqcup \Bbb A_1$ and $w(1),w(2)$ are the two vertices, or
$Q = \Bbb A_3$ and $w(1),w(2)$ the two leaves, or
$Q = \Bbb D_n$ with $n\ge 5$ and $w = w(1)$ is the leaf of the long arm.}
	\medskip 
In all these cases, we have $Q = \Bbb D_n$ with $n\ge 2$, where $\Bbb D_2 =  
\Bbb A_1\sqcup \Bbb A_1$ and $\Bbb D_3 = \Bbb A_3$. The hammock set for $\Bbb D_n$ 
mentioned in Proposition 1 is the following self-dual poset with $2n-2$ vertices:
$$
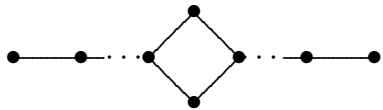
{\beginpicture
\setcoordinatesystem units <.5cm,.5cm> 
\multiput{$\bullet$} at 2 -2  2 0  /
\multiput{$\bullet$} at 1 -1   -.5 -1 -2 -1   3 -1  4.5 -1  6 -1  /

\multiput{$\cdots$} at 0.5 -1  3.5 -1 /
\plot -2 -1  0 -1 /
\plot 1 -1  2 0  3 -1 /
\plot 4 -1  6 -1 /

\plot 1 -1  2 -2  3 -1 /
\endpicture}
$$
	\medskip 
Given a poset $P$, a {\it $(2,2)$-set} in $P$ is by definition a full subposet 
which is the disjoint union of two chains of cardinality $2.$ (More generally, 
an $(n_1,\dots,n_t)$-set in $P$ is a full subposet of $P$ which is the 
disjoint union of chains of cardinalities $n_1,\dots, n_t.$ In particular, an
antichain in $P$ is a $(1,1,\dots,1)$-set in $P$.)
	\medskip 
{\bf Proposition 2.} {\it 
Let $Q(1),\dots, Q(t)$ be Dynkin quivers and let $Q$ be the disjoint union
of the quivers $Q(i)$. Let $w(i)$ be a vertex of $Q(i)$, for $1\le i \le s.$ Then
$H(Q,w(1),\dots,w(s))$ has width $2$ and contains a unique $(2,2)$-set 
if and only $Q$ is a quiver of the following form, and $w(1),\dots, w(s)$ are the
encircled vertices:}
$$
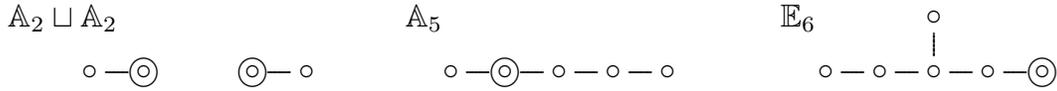
{\beginpicture
\setcoordinatesystem units <.6cm,.6cm> 
\put{\beginpicture
\multiput{$\circ$} at 0 0  1 0  3 0  4 0 /
\multiput{$\bigcirc$} at 1 0  3 0 /
\plot 0.3 0  0.7 0 /
\plot 3.3 0  3.7 0 /
\put{} at 0 1.05
\put{$\Bbb A_2\sqcup\Bbb A_2$} at -.5 1
\endpicture} at 0 0 
\put{\beginpicture
\multiput{$\circ$} at 0 0  1 0  2 0  3 0  4 0 /
\put{$\bigcirc$} at 1 0
\plot 0.3 0  0.7 0 /
\plot 1.3 0  1.7 0 /
\plot 2.3 0  2.7 0 /
\plot 3.3 0  3.7 0 /
\put{$\Bbb A_5$} at -.5 1
\endpicture} at 7 0 
\put{\beginpicture
\multiput{$\circ$} at -1 0  0 0  1 0  1 1  2 0  3 0  /
\put{$\bigcirc$} at 3 0
\plot -.3 0  -.7 0 /
\plot 0.3 0  0.7 0 /
\plot 1.3 0  1.7 0 /
\plot 2.3 0  2.7 0 /
\plot 1 0.3  1 0.7 /
\put{$\Bbb E_6$} at -1.5 1
\endpicture} at 14 0 
\endpicture}
$$
	\medskip 
Here are the hammock sets $H(Q,w(1),\dots,w(s))$ of width 2 with precisely one
$(2,2)$-set. Of course, $s \le 2$ and we
write $w = w(1)$ and $w' = w(2)$. 
$$
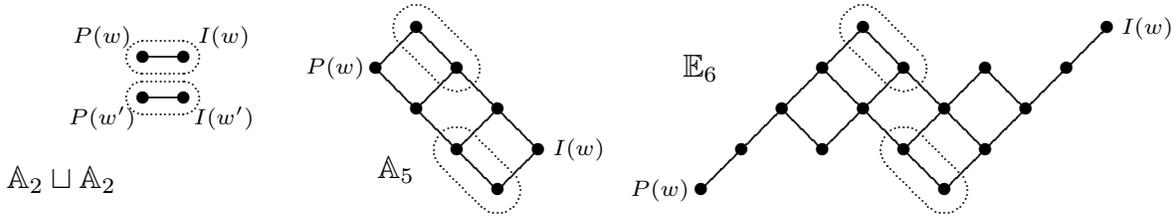
{\beginpicture
\setcoordinatesystem units <.45cm,.45cm> 
\put{\beginpicture
\multiput{$\bullet$} at 0 0  1 0   0 1  1 1   /
\plot 0 0  1 0 /
\plot 0 1  1 1 /
\setdots <.5mm>
\plot 0 0.4  1 0.4 /
\plot 0 -.4  1 -.4 /
\circulararc -180 degrees from 0 -.4  center at 0 0 
\circulararc 180 degrees from 1 -.4  center at 1 0
\put{$\ssize P(w)$} at -1 1.5
\put{$\ssize P(w')$} at -1 -.5
\put{$\ssize I(w)$} at 2 1.5
\put{$\ssize I(w')$} at 2 -.5
\plot 0 1.4  1 1.4 /
\plot 0 .6  1 .6 /

\circulararc 180 degrees from 0 1.4  center at 0 1 
\circulararc -180 degrees from 1 1.4  center at 1 1
\put{$\Bbb A_2\sqcup\Bbb A_2$} at -2 -2 
\endpicture} at -2 0
\put{\beginpicture
\multiput{$\bullet$} at 0 3  1 2  1 4  2 1 2 3  3 0  3 2  4 1  /
\put{$\ssize P(w)$} at -1 3
\put{$\ssize I(w)$} at 5 1
\plot 0 3  1 4  4 1  3 0  0 3 /
\plot 1 2  2 3 /
\plot 2 1  3 2 /
\setdots <.5mm>
\plot 0.6 3.6  1.6 2.6 /
\plot 1.4 4.4  2.4 3.4 /
\circulararc -180 degrees from .6 3.6  center at 1 4 
\circulararc 180 degrees from 1.6 2.6  center at 2 3  

\plot 1.6 .6  2.6 -.4 /
\plot 2.4 1.4  3.4 .4 /
\circulararc -180 degrees from 1.6 0.6  center at 2 1 
\circulararc 180 degrees from 2.6 -.4  center at 3 0 
\put{$\Bbb A_5$} at 0.5 0.5 
\endpicture} at 6 0
\put{\beginpicture
\multiput{$\bullet$} at 0 3  1 2  1 4  2 1 2 3  3 0  3 2  4 1  /
\multiput{$\bullet$} at 0 1  -1 2  -2 1  -3 0  4 3  5 2  6 3  7 4   /
\put{$\ssize P(w)$} at -4 0
\put{$\ssize I(w)$} at 8 4
\plot -3 0  1 4  4 1  3 0  0 3 /
\plot 1 2  2 3 /
\plot 2 1  3 2 /
\plot -1 2  0 1  1 2 /
\plot 4 1  7 4 /
\plot 3 2  4 3  5 2 /
\setdots <.5mm>
\plot 0.6 3.6  1.6 2.6 /
\plot 1.4 4.4  2.4 3.4 /
\circulararc -180 degrees from .6 3.6  center at 1 4 
\circulararc 180 degrees from 1.6 2.6  center at 2 3  

\plot 1.6 .6  2.6 -.4 /
\plot 2.4 1.4  3.4 .4 /
\circulararc -180 degrees from 1.6 0.6  center at 2 1 
\circulararc 180 degrees from 2.6 -.4  center at 3 0 
\put{$\Bbb E_6$} at -3 3 
\endpicture} at 17 0 
\endpicture}
$$
The unique $(2,2)$-subsets are encircled by dotted lines. Again, we have poset 
embeddings:
$$
 H(\Bbb A_2\sqcup\Bbb A_2,w,w') \subset H(\Bbb A_5,w) \subset H(\Bbb E_6,w).
$$
If we want to visualize that the poset $H(\Bbb E_6,w)$ is a subquiver of $\Bbb Z\Bbb E_6$,
we should draw it slightly different, namely as follows:
$$
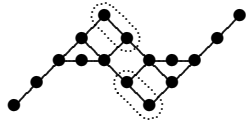
\beginpicture
\setcoordinatesystem units <.25cm,.25cm> 
\multiput{$\bullet$} at 0 3  1 2  1 4  2 1 2 3  3 0  3 2  4 1  /
\multiput{$\bullet$} at 0 2  -1 2  -2 1  -3 0  4 2  5 2  6 3  7 4   /
\plot -3 0  1 4  4 1  3 0  0 3 /
\plot 1 2  2 3 /
\plot 2 1  3 2 /
\plot -1 2   1 2 /
\plot 4 1  7 4 /
\plot 3 2  5 2 /
\setdots <.5mm>
\plot 0.6 3.6  1.6 2.6 /
\plot 1.4 4.4  2.4 3.4 /
\circulararc -180 degrees from .6 3.6  center at 1 4 
\circulararc 180 degrees from 1.6 2.6  center at 2 3  

\plot 1.6 .6  2.6 -.4 /
\plot 2.4 1.4  3.4 .4 /
\circulararc -180 degrees from 1.6 0.6  center at 2 1 
\circulararc 180 degrees from 2.6 -.4  center at 3 0 
\endpicture
$$

Given one of these hammock sets, the representations belonging to the $(2,2)$-subset
will be denoted by $U,V,U',V'$ with $\Hom(U,V) \neq 0,\ \Hom(U',V') \neq 0$.
Of course, in the cases $\Bbb A_2\sqcup\Bbb A_2$ and $\Bbb A_5$, the 
representations $U,V,U',V'$ are thin (as all indecomposable representations of a
quiver of type $\Bbb A$). For $Q$ of type $\Bbb E_6$, two of the representations 
$U,V,U',V'$ are thin, the remaining two are representations $N$ with $\dim N_c = 2$,
where $c$ is the central vertex of $Q$. 
	\medskip
If we consider the quivers of type $\Bbb E_6, \Bbb E_7, \Bbb E_8$ with
subspace orientation, the representations $U,V,U',V'$ are as follows:
$$
{\beginpicture
\setcoordinatesystem units <1cm,.65cm> 
\put{\beginpicture
\put{${\Bbb E}_6$} at -.7 1.3
\put{\beginpicture
  \setcoordinatesystem units <.6cm,.6cm> 
  \multiput{$\circ$} at 0 0  1 0   3 0  4 0 /
  \put{$y$} at 2 1 
\arr{2 0.7}{2 0.3} 
  \put{$x$} at 2 0
  \arr{0.2 0}{0.8 0} 
  \arr{1.2 0}{1.8 0} 
  \arr{2.8 0}{2.2 0} 
  \arr{3.8 0}{3.2 0} 
  \endpicture} at 0 1 
\put{$U$} at -2.5 -.7
\put{$V$} at -2.5 -2.1
\put{$U'$} at -2.5 -3.5
\put{$V'$} at -2.5 -4.9
\put{$\smallmatrix 
                   & & 0 \cr 
         0 & 1 & 0 & 0 & 0\endsmallmatrix$} at 0 -.7
\put{$\smallmatrix  
                    & & 0 \cr 
         1 & 1 & 0 & 0 & 0\endsmallmatrix$} at 0 -2.1
\put{$\smallmatrix  
                    & & 0 \cr 
        0 & 0 & 0 & 1 & 0\endsmallmatrix$} at 0 -3.5
\put{$\smallmatrix  
                    & & 0 \cr 
        0 & 0 & 0 & 1 & 1\endsmallmatrix$} at 0 -4.9
\endpicture} at -.3 0
\put{\beginpicture
\put{${\Bbb E}_7$} at -1 1.3
\put{\beginpicture
  \setcoordinatesystem units <.6cm,.6cm> 
  \multiput{$\circ$} at 1 1  1 0  2 0  3 0  4 0 /
  \put{$x\strut$} at 0 0 
  \put{$y\strut$} at -1 0 
  \arr{0.2 0}{0.8 0} 
  \arr{1.8 0}{1.2 0} 
  \arr{2.8 0}{2.2 0} 
  \arr{3.8 0}{3.2 0} 
  \arr{1 0.8}{1 0.2}
  \arr{-.8 0}{-.2 0} 

  \endpicture} at 0 1
\put{$\smallmatrix   && 1 \cr
                   0 & 0 & 1 & 0 & 0 & 0 \endsmallmatrix$} at 0 -.7
\put{$\smallmatrix   && 1 \cr
                   0 &  0 & 1 & 1 & 0 & 0\endsmallmatrix$} at 0 -2.1
\put{$\smallmatrix   && 0 \cr
                   0 &  0 & 1 & 1 & 1 & 0\endsmallmatrix$} at 0 -3.5
\put{$\smallmatrix   && 0 \cr
                   0 &  0 & 1 & 1 & 1 & 1\endsmallmatrix$} at 0 -4.9
\endpicture} at 4.5 0
\put{\beginpicture
\put{${\Bbb E}_8$} at -1.6 1.3
\put{\beginpicture
  \setcoordinatesystem units <.6cm,.6cm> 
  \multiput{$\circ$} at -1 0  0 0  1 1  1 0  2 0  3 0   /
  \put{$\strut x$} at 4 0 
  \arr{-.8 0}{-.2 0} 
  \arr{0.2 0}{0.8 0} 
  \arr{1.8 0}{1.2 0} 
  \arr{2.8 0}{2.2 0} 
  \arr{3.8 0}{3.2 0} 
  \arr{1 0.8}{1 0.2}
\put{$\strut y$} at 5 0 
  \arr{4.8 0}{4.2 0} 
  \endpicture} at 0 1
\put{$\smallmatrix &  & 0 \cr
                   0 & 1 & 1 & 1 & 1 & 0 & 0  \endsmallmatrix$} at -.4 -.7 
\put{$\smallmatrix &  & 1 \cr
                   1 & 2 & 2 & 1 & 1 & 0 & 0  \endsmallmatrix$} at -.4 -2.1 
\put{$\smallmatrix &  & 1 \cr
                   1 & 1 & 2 & 2 & 1 & 0 & 0  \endsmallmatrix$} at -.4 -3.5
\put{$\smallmatrix &  & 1 \cr
                   0 & 0 & 1 & 1 & 1 & 0 & 0  \endsmallmatrix$} at -.4 -4.9
\endpicture} at 9 0
\endpicture}
$$
	\bigskip
Let us complete the proof of Theorem 3.8. 
There are no problems if
$\Delta$ is of type $\Bbb D_n$, thus we assume that $\Delta$ is of type $\Bbb E_m$.
We are going to discuss the relationship
between the representations of $\Delta'$ and of $\Delta''.$ 
The quiver $\Delta''$ is of the form $\Bbb A_2\sqcup\Bbb A_2,\ \Bbb A_5,\ \Bbb E_6$
and the neighbors in $\Delta''$ of the vertex $x\in\Delta_0$ are just the vertices
$w(1),\dots,w(s)$ mentioned in Proposition 2. 
We write $H(\Delta'')$ instead of $H(\Delta'',w(1),\dots,w(s))$,
thus, for $\Delta''$ of type $\Bbb A_2\sqcup\Bbb A_2$, we have 
$H(\Delta'') = H(\Delta'',w,w')$, and for $\Delta''$ of type $\Bbb A_5$ and $\Bbb E_6$, 
we have $H(\Delta'')= H(\Delta'',w).$

Given a representation $N$ in $H(\Delta'')$, we denote
by $\overline N$ the indecomposable representation of $\Delta'$ with 
$\overline N|\Delta'' = N$
and $\overline N_x \neq 0$ (thus $\dim \overline N_x = 1$). 
Given an antichain pair $N,N'$ in $H(\Delta'')$, we denote
by $\overline {NN'}$ the indecomposable representation of $\Delta'$ with 
$\overline{NN'}|\Delta'' = N\oplus N'$ (and $\dim \overline{NN'}_x = 1)$.
The representations of $\Delta'$ of the form $\overline N$ and $\overline{NN'}$
with $N,N' \in H(\Delta'')$ together with $S(x)$ form the hammock
$H(\Delta',x)$. 

Let us describe inside the hammock $H(\Delta',x)$ the subset $\Cal X$
given by the indecomposable representations $X$ of $\Delta'$ with $X|\Delta''
\in \add\{U,V,U',V'\}.$ This depends on the orientation of the arrows outside
of $\Delta''.$ 
Up to duality, we can assume that the arrow between $x$ and $w$
points to $w$. In the case of $\Delta''$ being of type $\Bbb A_2\sqcup \Bbb A_2$
we have to distinguish the two cases $x \to w'$ (thus $x$ is a source)
and $x \leftarrow w'.$

First, let us assume that $x$ is a source. Then $\Cal X$ looks as follows:
$$
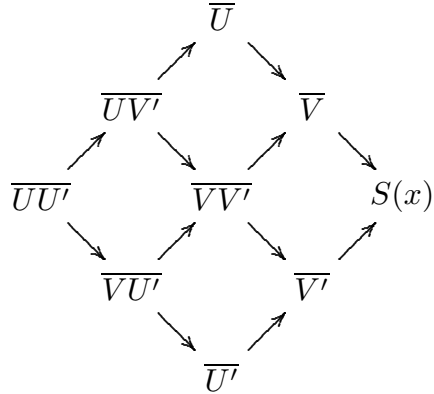
{\beginpicture
\setcoordinatesystem units <1cm,1cm> 
\put{$\overline{UU'}$} at 0 0 
\put{$\overline{UV'}$} at 1 1
\put{$\overline{VU'}$} at 1 -1
\put{$\overline{U}$} at 2 2 
\put{$\overline{VV'}$} at 2 0 
\put{$\overline{U'}$} at 2 -2
\put{$\overline{V}$} at 3 1
\put{$\overline{V'}$} at 3 -1
\put{$S(x)$} at 4 0 
\arr{0.3 0.3}{0.7 0.7}
\arr{1.3 1.3}{1.7 1.7}

\arr{1.3 -.7}{1.7 -.3}
\arr{2.3 .3}{2.7 .7}

\arr{2.3 -1.7}{2.7 -1.3}
\arr{3.3 -.7}{3.7 -.3}

\arr{0.3 -.3}{0.7 -.7}
\arr{1.3 -1.3}{1.7 -1.7}

\arr{1.3 .7}{1.7 .3}
\arr{2.3 -.3}{2.7 -.7}

\arr{2.3 1.7}{2.7 1.3}
\arr{3.3 .7}{3.7 .3}
\endpicture}
$$

Second, we assume that $\Delta''$ is of type $\Bbb A_2\sqcup \Bbb A_2$ and that
there is the subquiver $w \leftarrow x \leftarrow w'$ in $\Delta'.$ Then $\Cal X$
looks as follows:
$$
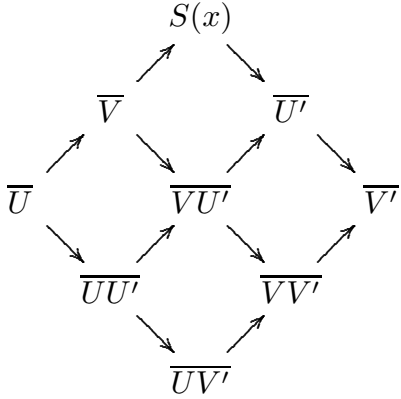
{\beginpicture
\setcoordinatesystem units <1cm,1cm> 

\put{$\overline{U}$} at 0 0
\put{$\overline{V}$} at 1 1
\put{$\overline{UU'}$} at 1 -1 
\put{$S(x)$} at 2 2
\put{$\overline{VU'}$} at 2 0
\put{$\overline{UV'}$} at 2 -2
\put{$\overline{U'}$} at 3 1
\put{$\overline{VV'}$} at 3 -1 
\put{$\overline{V'}$} at 4 0

\arr{0.3 0.3}{0.7 0.7}
\arr{1.3 1.3}{1.7 1.7}

\arr{1.3 -.7}{1.7 -.3}
\arr{2.3 .3}{2.7 .7}

\arr{2.3 -1.7}{2.7 -1.3}
\arr{3.3 -.7}{3.7 -.3}

\arr{0.3 -.3}{0.7 -.7}
\arr{1.3 -1.3}{1.7 -1.7}

\arr{1.3 .7}{1.7 .3}
\arr{2.3 -.3}{2.7 -.7}

\arr{2.3 1.7}{2.7 1.3}
\arr{3.3 .7}{3.7 .3}
\endpicture}
$$

Always we see that the subset $\Cal X$ of $H(\Delta',x)$ 
contains an antichain triple. But as we know $H(\Delta',x)$ has a unique antichain
triple, namely the representations $A(1),A(2),A(3).$ This shows that
$$
 \{A(1),A(2),A(3)\} = \left\{\matrix \{\overline{U},\overline{VV'},\overline{U'}\},
    & \text{if $x$ is a source,}\cr
      \{S(x),\overline{U'V},\overline{UV'}\}
    & \ \text{for $w \leftarrow x \leftarrow w'.$}\endmatrix \right.
$$
Of course, in both cases we have
$$
 A(1)\oplus A(2)\oplus A(3)\, |\,\Delta'' = U\oplus V \oplus U'\oplus V',
$$
this completes the proof of Theorem 3.8.
	\bigskip 
{\bf Proposition 3.} {\it Let $\Delta$ be a quiver of type $\Bbb E_6,\Bbb E_7,\Bbb E_8$
with exceptional vertex $y$ and $x$ the neighbor of $y$.
Then there is a unique thick subcategory of $\rep\Delta$ which
contains $S(x)$ and $S(y)$ such that its
$\Ext$-quiver is of the following form}
$$
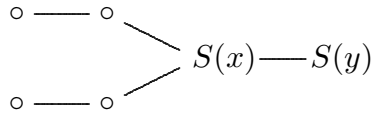
{\beginpicture
\setcoordinatesystem units <1cm,.5cm> 
\multiput{$\circ$} at 0 1  1 1  0 -1  1 -1 /
\put{$S(x)$} at 2.3 0
\put{$S(y)$} at 3.6 0
\plot 2.7 0  3.2 0 /
\plot 0.2 1  0.8 1 /
\plot 0.2 -1  0.8 -1 /
\plot 1.2 0.8 1.8 0.2 /
\plot 1.2 -.8 1.8 -.2 /
\endpicture}
$$
	\medskip
Remark. In case $\Delta$ is of type $\Bbb E_6$, there is only one 
thick subcategory of $\rep\Delta$ with 6 simple objects, namely $\rep\Delta$
itself, thus nothing has to be shown. 
	\medskip 
Proof. It is sufficient to observe that the smallest thick subcategories containing $U,V$
of $U',V'$, respectively, both are of type $\Bbb A_2$. 
	\bigskip 
Let us add a table which provides for any quiver $\Delta$ of Dynkin 
type $\Bbb D_n$ or $\Bbb E_m$
the number of indecomposable representations $X$ of
$\Delta'$ of the various kinds (similar to the table in section 3.5).
$$
{\beginpicture
\setcoordinatesystem units <1cm,.7cm>
\put{} at -4.5 2

\put{$\Delta'$} at -3 2
\put{$\Bbb D_{n-2}$} at -.7 2 
\put{$\Bbb A_5$} at 1 2
\put{$\Bbb D_6$} at 2 2  
\put{$\Bbb E_7$} at 3 2

\put{$\Delta''$} at -3 1
\put{$\Bbb D_{n-3}$} at -.7 1 
\put{$\Bbb A_2\!\sqcup\!\Bbb A_2$} at 1 1
\put{$\Bbb A_5$} at 2 1  
\put{$\Bbb E_6$} at 3 1

\put{$X \in \rep \Delta''$} at -4 0 
\put{$\m(X|\Delta'') = 1,\ X_x\neq 0$} at -4 -1 
\put{$\m(X|\Delta'') = 2$} at -4 -2
\put{$X|\Delta'' = 0$} at -4 -3

\put{sum} at -6.5 -4 

\put{$(n-4)(n-3)$} at -.7 0
\put{$2n-8$} at -.7 -1 
\put{1} at -.7 -2 
\put{1} at -.7 -3 
\put{$(n-3)(n-2)$} at -.7 -4

\put{6} at 1 0 
\put{4} at 1 -1 
\put{4} at 1 -2 
\put{1} at 1 -3 
\put{15} at 1 -4 
  
\put{15} at 2 0 
\put{8} at 2 -1 
\put{6} at 2 -2 
\put{1} at 2 -3 
\put{30} at 2 -4 

\put{36} at 3 0 
\put{16} at 3 -1 
\put{10} at 3 -2 
\put{1} at 3 -3 
\put{63} at 3 -4 
 
\plot -7 .5  3.5 .5 /
\plot -7 -3.5  3.5 -3.5 /
\plot -2 2.4  -2 -4.4 /
\endpicture}
$$
Of course, the number of indecomposable representations $X$ of $\Delta'$
with $\m(X|\Delta'') = 1$ and $X_x \neq 0$ is just the number of vertices
of the hammock set $H(\Delta'')$, whereas the 
number of indecomposable representations $X$ of $\Delta'$
with $\m(X|\Delta'') = 2$ is the number of antichain pairs in 
$H(\Delta'')$.
		       \bigskip
{\bf 3.9. Invariants with value $2,\ 4,\ 8$.} 
Let $\Delta$ be a Dynkin quiver of type $\Bbb E_m$ with exceptional vertex $y$. 
Let $x$ be the neighbor of $y$ and $\Delta' = \Delta\setminus\{y\},
\Delta' = \Delta\setminus\{x,y\}.$ 

Let $P'(x)$ be the indecomposable projective 
representation of $\Delta'$ corresponding
to the vertex $x$, let $I'(x)$ 
be the indecomposable injective 
representation of $\Delta'$ corresponding
to the vertex $x$. Let $\tau'$ be the Auslander-Reiten translation in 
$\rep \Delta'$.
{\it Then $P'(x)$ and $I'(x)$
belong to the same $\tau'$-orbit  and
there is the following $2-4-8$ assertion:}
$$
 P'(x) = \left\{\matrix (\tau')^2 I'(x)  \cr 
                       (\tau')^4 I'(x)  \cr
                       (\tau')^8 I'(x) 
               \endmatrix \right.
\quad\text{for $\Delta'$ of type} \left\{\matrix \Bbb A_5\, , \cr 
                    \Bbb D_6\, , \cr
                    \Bbb E_7\, , \cr  
               \endmatrix \right.
\quad\text{thus for $\Delta$ of type} \left\{\matrix \Bbb E_6\, , \cr 
                    \Bbb E_7\, , \cr
                    \Bbb E_8\, . \cr  
               \endmatrix \right.
$$
Actually, also the number of vertices of the hammock $H(\Delta',x)$ 
is related to this number $r = 2,4,8$, namely it is $3(r+1)$ 
(and, as we have seen, the number of vertices of the hammock coincides with the
number of antichain pairs in the hammock).

Already the quiver $\Delta''$ with its vertices $w(1),\dots,w(s)$ (the neighbors of
$x$ in $\Delta'$) provides a $2-4-8$ assertion. Namely, the number of vertices of
the hammock set $H(\Delta'')$ is equal to $2r,$ with $r=2,4,8$
for $\Delta$ of type $\Bbb E_6,\Bbb E_7,\Bbb E_8$, respectively. 
       \bigskip
{\bf 3.10. Final Remarks for Part 3.} It is a well-known procedure to use one-point extensions
(and vector space categories) in order to determine the representation type of quivers,
in particular to deal with the quivers of type $\Bbb E_6,\Bbb E_7, \Bbb E_8$. A usual
way seems to be to use a sequence of one-point extensions starting with $\Bbb D_5$, 
and adding successively vertices to obtain $\Bbb E_6, \Bbb E_7, \Bbb E_8$. 
Our method presented is slightly different: whereas as usual $\Bbb E_8$ is obtained as 
a one-point extension from $\Bbb E_7$, we draw the attention to a parallel way for obtaining
$\Bbb E_7$ from $\Bbb D_6$ and $\Bbb E_6$ from $\Bbb A_5.$ In this way, the three exceptional
cases $\Bbb E_6, \Bbb E_7, \Bbb E_8$ are considered as being on equal footing.

Let us stress that any Dynkin diagram $\Gamma$ should be seen as being accompanied by
the corresponding Euclidean diagram $\widetilde \Gamma$; many properties of $\Gamma$ can be read off
very nicely by looking at $\widetilde \Gamma.$ In contrast to the inclusion sequence of the
Dynkin diagrams $\Bbb E_6, \Bbb E_7, \Bbb E_8$, the corresponding Euclidean diagrams $\widetilde {\Bbb E}_6,
\widetilde {\Bbb E}_7, \widetilde {\Bbb E}_8$ are pairwise incomparable. Of course, our equal footing approach to
$\Bbb E_6, \Bbb E_7, \Bbb E_8$ is motivated by a look at the further extensions from $\Gamma$
to $\widetilde \Gamma$. 

Henning Krause has drawn my attention to a famous letter included in the paper [A] by F.~Adams,
and purportedly written by $\Bbb E_8$ itself: it regrets the prevailing opinion that the case
$\Bbb E_8$ is {\it remote and unapproachable} and needs {\it an arduous course of preparation with 
$\Bbb E_6$ and $\Bbb E_7$.} The letter adds the following comment:
{\it Any right-thinking mathematician who wishes to
construct the root-system of $\Bbb E_6$ does so as follows: 
first he constructs the root-system of $\Bbb E_8$, and
then inside it he locates the root-system of $\Bbb E_6$. 
In this way he benefits from the great symmetry of
the root-system of $\Bbb E_8$ and its perspicuous nature. 
If this good precedent is not followed in other
researches, one should consider whether to infer a lack of boldness 
in the investigator rather than
a lack of cooperation from the subject-matter.}

In contrast to Adams suggestion to start with $\Bbb E_8$ and 
to consider the root systems of $\Bbb E_6$ and $\Bbb E_7$ just
as subsystems of $\Bbb E_8$, we have the feeling that our attempt to provide
a unified treatment of $\Bbb E_6,\Bbb E_7, \Bbb E_8$ is also 
supported by the subject-matter. 
	  \medskip 
Let us add a remark concerning our use of the representations of a quiver $R$ of type 
$\widetilde{\Bbb D}_4$ in order to deal with the representations of an arbitrary Euclidean quiver. To invoke $R$ means a reduction to the 4-subspace quiver.
In the paper [GP], Gelfand and Ponomarev have drawn the attention to the 4-subspace
quiver and have provided a classification of all the indecomposable representation. 
In the paper, they have stressed the importance of the 4-subspace quiver:
{\it Many tame problems of Linear Algebra can be reduced to the classification} of 
representations of the 4-subspace quiver ({[GP]} p. 163). The results presented in Part 3
confirm this observation. Our approach provides for the Euclidean quivers $\widetilde \Delta$ of type
$\widetilde{\Bbb D}_n, \widetilde{\Bbb E}_6, \widetilde{\Bbb E}_7, \widetilde{\Bbb E}_8$
a thick subcategory of type $\widetilde{\Bbb D}_4$. It recovers not only 
the 1-parameter families of indecomposable representations, but also one 
representative from any $\tau$-orbit of simple regular representations.
Of course, if one is only interested in the 1-parameter families of indecomposable
representations, one will be content with a thick subcategory which is equivalent to the
category of Kronecker modules, say with $\Cal E(M,S(z))$
where $M$ is the maximal indecomposable representation of a corresponding Dynkin subquiver
$\Delta$ and $z$ is the vertex in $\widetilde\Delta \setminus \Delta$. 
	\medskip
The $2-4-8$ assertion should remind the reader on mathematical settings which rely on
an invariant $r$ being equal to 1, 2, 4, or 8. The prototype is the Hurwitz theorem
which asserts that the only Euclidean Hurwitz algebras are the real numbers $\Bbb R$, 
the complex numbers $\Bbb C$, the quaternions $\Bbb H$ and the octonions $\Bbb O$ 
(a {\it Euclidean Hurwitz algebra} is a finite-dimensional 
(not necessarily associative) $\Bbb R$-algebra $A$ with identity 
endowed with a positive quadratic form $q$ such that $q(ab) = q(a)q(b)$).
Is there a direct relationship between the Dynkin quivers $\Bbb E_6,\ \Bbb E_7, \Bbb E_8$
and the division algebras $\Bbb C,\ \Bbb H,\ \Bbb O$? We do not know. But we should
draw the attention of the reader to the magic Freudenthal-Tits  square, see 
[F, T, V]. It is a symmetric $(4\times 4)$-matrix whose entries $D_{AB}$ 
are Dynkin types, with column index $A$ and row index $B$ being one  
$\Bbb R,\ \Bbb C,\ \Bbb H,\ \Bbb O$:
$$  
 {\beginpicture
 \setcoordinatesystem units <1cm,.7cm> 
\multiput{$\Bbb R$} at 1 5  0 4 / 
\multiput{$\Bbb C$} at 2 5  0 3 /
\multiput{$\Bbb H$} at 3 5  0 2 /
\multiput{$\Bbb O$} at 4 5  0 1 /
\plot 0.5 0.6  .5 5.2 /
\plot -.2 4.5  4.3 4.5 /
\multiput{$\Bbb F_4$} at 4 4  1 1 /  
\multiput{$\Bbb  E_6$} at 4 3  2 1 /
\multiput{$\Bbb  E_7$} at 4 2  3 1 /
\multiput{$\Bbb  E_8$} at 4 1  /
\multiput{$\Bbb  A_1$} at 1 4 / 
\multiput{$\Bbb  A_2$} at 2 4  1 3 /
\multiput{$\Bbb  C_3$} at 3 4  1 2 /
\multiput{$\Bbb  A_2\!\sqcup\! \Bbb A_2$} at 2 3 / 
\multiput{$\Bbb  A_5$} at 2 2  3 3 /
\multiput{$\Bbb  D_6$} at 3 2 /
\endpicture}
$$
The magic square
 concerns constructions of the Lie groups or Lie algebras of type $D_{AB}$ 
starting with $A$ and $B$. Since all possible $A$ and $B$ are 
subalgebras of $\Bbb O$, the
magic square asserts that all exceptional Lie groups and Lie algebras can
be constructed starting with the octonions $\Bbb O$ (note that the exceptional
Lie type $\Bbb G_2$ is missing in the magic square, but the Lie group
of type $\Bbb G_2$ is just the automorphism group of $\Bbb O$). 

We are interested in the right lower corner of the magic square, it
is the following matrix:
$$  
 {\beginpicture
 \setcoordinatesystem units <1cm,.7cm> 
\multiput{$\Bbb  E_6$} at 4 3  2 1 /
\multiput{$\Bbb  E_7$} at 4 2  3 1 /
\multiput{$\Bbb  E_8$} at 4 1  /
\multiput{$\Bbb  A_2\!\sqcup\! \Bbb A_2$} at 2 3 / 
\multiput{$\Bbb  A_5$} at 2 2  3 3 /
\multiput{$\Bbb  D_6$} at 3 2 /
\endpicture}
$$
Here, any row (or column) lists the Dynkin types of $\Delta'',\Delta',\Delta$, 
where $\Delta$ is a Dynkin quiver of type $\Bbb E_m$ with exceptional vertex $y$
and where $x$ is the neighbor of $y$, such that $\Delta' = \Delta\setminus\{y\}$ and 
$\Delta'' = \Delta\setminus\{x,y\}.$  

What we have recovered in our presentation are invariants $r$ of $\Delta$ which 
have the value $2,4,8$ respectively. 

One may wonder about the case $r= 1.$
According to the magic square, this concerns the Dynkin type $\Bbb F_4$. Now,
there is no quiver of type $\Bbb F_4.$ There are the species of type $\Bbb F_4$,
but they are outside of the scope of the present investigation.
Actually, it turns out that for the species of type $\Bbb F_4$,
the table in 3.5 can be completed by a column concerning $\Bbb F_4$ with entries
$\Bbb F_4, \Bbb C_3, \Bbb A_2, 3, 6, 6, 6, 2, 1, 24$, and 
$6 = 3\cdot(r+1)$ with $r=1.$ This confirms the magic square philosophy. 
However, if we look at
the hammock in question, we have $(\tau')^t I'(x) = P'(x)$
with $t = 2$ and not $t = 1$.

	 \bigskip\bigskip
{\bf Panoramic View.}
     \medskip
The Dynkin graphs have been exhibited in the Preliminaries. There, at the end, one also finds the
Euclidean graphs $\widetilde {\Bbb A}_n$. The remaining Euclidean graphs $\widetilde {\Bbb D}_n$,
$\widetilde {\Bbb E}_6,
\widetilde {\Bbb E}_7, \widetilde {\Bbb E}_8$ have been discussed in Part 3. 
 
 \medskip 
{\bf (1)} {\it A finite connected quiver $Q$ without a Euclidean subquiver is a Dynkin quiver.}
Proof: Since $Q$ does not contain a subquiver of type $\widetilde{\Bbb A}_n$, we see that 
$Q$ is a tree quiver
Since $Q$ has no subquiver of type $\widetilde{\Bbb D}_4$, any vertex of $Q$ has at most 3 neighbors.
Vertices with 3 neighbors are called branching vertices. 
Since $Q$ has no subquiver of type $\widetilde{\Bbb D}_n$, with $n\ge 5$, there is at most one
branching vertex. If there is no branching vertex, then $Q$ is of type $\Bbb A_n$, thus a
Dynkin quiver. If $Q$ has a branching vertex, then $Q$ is a star quiver of type $(t(1),t(2),t(3))$
with $t(1) \ge t(2) \ge t(3) \ge 2.$ If $t(3) \ge 3,$ then $Q$ contains a subquiver
of type $\widetilde{\Bbb E}_6.$ Thus $t(3) = 2$. If $t(2)\ge 4$, then $Q$ contains a subquiver
of type $\widetilde{\Bbb E}_7.$  Thus $t(2)$ is equal to $2$ or $3$. If $t(2) = 2$, then $Q$
is of Dynkin type $\Bbb D_n$. It remains that $t(2) = 3$. 
If $t(1)\ge 6$, then $Q$ contains a subquiver
of type $\widetilde{\Bbb E}_8.$ Thus, the only remaining cases are $t(1)\in\{3,4,5\}$, these are
the Dynkin cases $\Bbb E_6, \Bbb E_7, \Bbb E_8.$
    \medskip 

{\bf (2)} {\it Dynkin quivers are representation-finite.}
This has been shown above.

     \medskip
{\bf (3)} {\it Euclidean quivers are representation-infinite.}
     \smallskip 
Proof. The quivers of type $\widetilde {\Bbb A}_n$ have been discussed at the end
of the Preliminaries in section 0.3.
Let $Q$ be such a quiver and let us fix one of the arrows, say $\alpha$. For 
any $\lambda\in k$, we define $M(\lambda)$ as follows: $M(\lambda)_x = k$ for all vertices $x$, let
$M(\alpha)_\alpha = \lambda$ and $M(\lambda)_\beta = 1$ for the remaining arrows $\beta$. Then the
representations $M(\lambda)$ with $\lambda\in k$, are indecomposable and pairwise non-isomorphic.
This yields infinitely many isomorphism classes of indecomposable representations, provided $k$
is an infinite field. More generally, for any natural number $n$, 
we define $M(\lambda,n)$ as follows: $M(\lambda,n)_x = k^n$ for all vertices $x$, let
$M(\lambda,n)_\alpha$ be the Jordan $(n\times n)$-matrix with eigenvalue $\lambda\in k$ and $M(\lambda,n)_\beta$ the
identity map, for the remaining arrows $\beta$. Then the
representations $M(\lambda,n)$ with $\lambda\in k$, are indecomposable and pairwise non-isomorphic.

Next, let us consider the 4-subspace quiver $Q$ (its type is $\widetilde{\Bbb D}_4$). 
We consider the quiver $\widetilde {\Bbb A}_0$ (one vertex, one loop) and define a
full exact embedding $\zeta\:\rep \widetilde {\Bbb A}_0 \to \rep Q$ as follows: a representation of 
$\widetilde {\Bbb A}_0$ is a pair $(V,\phi)$, where $V$ is a vector space and $\phi\:V \to V$ a
linear endomorphism (note that $\rep \widetilde {\Bbb A}_0$ can also be considered as the category of $k[T]$-modules,
where $k[T]$ is the polynomial ring in one variable, since a $k[T]$-module is also just a vector space
with a linear endomorphism, namely given by the multiplication with $T$).
If $(V,\phi)$ is a representation of $\widetilde {\Bbb A}_0$, let $\zeta(V,\phi)$
be the following representation of $Q$
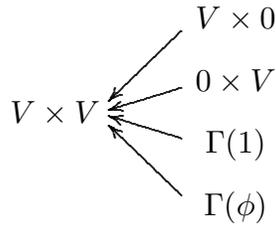
$$  
 {\beginpicture
 \setcoordinatesystem units <2cm,.7cm> 
\put{$V\times V$} at 0 0
\put{$V\times 0$} at 1 1.5
\put{$0\times V$} at 1 0.5
\put{$\Gamma(1)$} at 1 -.5
\put{$\Gamma(\phi)$} at 1 -1.5
\arr{0.7 1.3}{0.3 0.2}
\arr{0.7 0.4}{0.3 0.05}
\arr{0.7 -.4}{0.3 -.05}
\arr{0.7 -1.3}{0.3 -.2}
\endpicture}
$$
where $\Gamma(\phi) = \{(v,\phi(v))\mid v\in V\}$ (and thus $\Gamma(1) = \{(v,v)\mid v\in V\}$).
Since $\widetilde {\Bbb A}_0$ is representation-infinite, also $Q$ is representation-infinite.
Of course, almost all the indecomposable representations of $Q$ are conical, thus any quiver of type
$\widetilde{\Bbb D}_4$ is representation-infinite, by Proposition 1.6. 
		 
Consider now a Euclidean quiver $Q$ of type $\widetilde {\Bbb E}_n$ (thus $n=6,7,8$), but also 
those of type $\widetilde{\Bbb D}_n$ with $n\ge 5.$
As we have seen in Part 3,
there is a full exact embedding $\rep R \to \rep Q$, where $R$ is of type $\widetilde{\Bbb D}_4$.
This shows that the quivers of type $\widetilde {\Bbb E}_6,\ \widetilde {\Bbb E}_7,\ \widetilde {\Bbb E}_8$
are representation-infinite. 

    \bigskip\bigskip
{\bf References.}
     \medskip
\item{[A]} J.~F.~Adams. Finite $H$-spaces and Lie groups. Journal Pure Appl.~Algebra 19 (1980), 
    1--8.

\item{[ARS]} M.~Auslander, I.~Reiten, S.~S.~Smal\o. {\it Representation Theory of  
    Artin Algebras.} Cambridge Studies in Advanced Mathematics 36. Cambridge University 
    Press. 1997.

\item{[BGP]} I.~N.~Bernstein,  I.~M.~Gelfand, V.~A.~Ponomarev.
   Coxeter functors, and Gabriel's theorem. Russian mathematical surveys 28 (2) (1973) 17–-32,

\item{[B]} S.~Brenner. A combinatorial characterization of finite Auslander-Reiten quivers.
   In: {\it Representation Theory I: Finite Dimensional Algebras.} 
   Springer Lecture Notes in Mathematics 1177 (1986), 13--49

\item{[F]} H.~Freudenthal, Lie groups in the foundations of geometry.
   Adv. Math., 1 (1964), 145--190.

\item{[G1]} P.~Gabriel. Unzerlegbare Darstellungen I. Manuscripta Math. 6 (1972), 71--103.

\item{[G2]} P.~Gabriel. Auslander-Reiten sequences and representation-finite algebras.
    In: {\it Representation Theory I.} Springer Lecture Notes in Mathematics 831 (1980), 72--103.

\item{[GP]}  I.~M.~Gelfand, V.~A.~Ponomarev. Problems of  linear  algebra  and  classification
    of  quadruples of subspaces  in a finite-dimensional vector space.  Coll.~Math.~Soc. Janos Bolyai 5.
   {\it Hilbert Space Operators.} Tihany, Hungary (1970), 163--237.

\item{[R1]} C.~M.~Ringel. Representations of K-species and bimodules. J. Algebra 41 (1976),    269--302.

\item{[R2]} C.~M.~Ringel. Tame algebras.
    In: {\it Representation Theory I.} Springer Lecture Notes in Mathematics 831 (1980),  137--287. 

\item{[R3]} C.~M.~Ringel. {\it  Tame algebras and integral quadratic forms.} 
   Springer Lecture Notes in Mathematics 1099 (1984).
   
\item{[R4]} C.~M.~Ringel.  Distinguished bases of exceptional modules.
    In {\it Algebras, quivers and representations}. Proceedings of the Abel symposium 2011. 
    Springer Series Abel Symposia Vol 8. (2013), 253-274. 

\item{[R5]} C.~M.~Ringel. The Catalan combinatorics of the hereditary artin algebras. 
   To appear in: {\it Recent Developments in Representation Theory,} Contemp. Math., 673, Amer.~Math.~Soc., 
   Providence, RI, 2016. 

\item{[RV]} C.~M.~Ringel, D.~Vossieck. Hammocks. Proc.~London Math.~Soc. (3) 54 (1987), 216--246.

\item{[T]} J.~Tits. 
  Algebres alternatives, algebres de Jordan et algebres de Lie exceptionnelles, 
  Indag. Math. 28 (1966) 223-237. 
\item{[V]} E.~B.~Vinberg. A construction of exceptional simple Lie groups (Russian), 
   Tr. Semin. Vektorn. Tensorn. Anal. 13 (1966), 7-9. 

	    		    \bigskip\medskip
\noindent
Fakult\"at f\"ur Mathematik, Universit\"at Bielefeld\newline
POBox 100 131, D-33501 Bielefeld, Germany, and 
\smallskip

\noindent 
Department of Mathematics, Shanghai Jiao Tong University, \newline
Shanghai 200240, P. R. China.
	 \smallskip 

\noindent 
E-mail: {\tt ringel\@math.uni-bielefeld.de}

\bye